\documentclass[12pt,a4paper,twoside]{article}
\usepackage[utf8]{inputenc}
\usepackage[T1]{fontenc}
\usepackage[margin=2.5cm]{geometry}
\usepackage[french,english]{babel}
\usepackage{lmodern}
\usepackage{listings}
\usepackage{bef_aras_nothm}
\usepackage{amsmath,amsthm,amssymb}
\usepackage{fixmath}
\usepackage{upgreek}
\usepackage{enumitem}
\usepackage{amsfonts}
\usepackage{microtype}
\usepackage{float}
\usepackage{mathrsfs}
\usepackage{mathtools}

\usepackage{hyperref}
\hypersetup{
     colorlinks   = true,
     citecolor    = blue
}

\usepackage{fancyhdr}
\theoremstyle{plain}

\allowdisplaybreaks[1]

\numberwithin{equation}{section}

\newtheorem{thm}{Theorem}[section]
\newtheorem{lem}[thm]{Lemma}

\theoremstyle{definition}

\newtheorem{rem}[thm]{Remark}
\newtheorem{exa}[thm]{Example}

\begin{document}

\newpage
\setcounter{page}{1}

\title{Nonsmooth Analysis of Doubly Nonlinear Second-Order Evolution Equations with Non-Convex Energy Functionals} 

\author{Aras Bacho\footnotemark[1]}

\date{}
\maketitle

\footnotetext[1]{Department of Computing and Mathematical Sciences, California Institute of Technology, Pasadena, CA}


\begin{abstract}
We investigate the existence of strong solutions to a general class of doubly multivalued and nonlinear evolution equations of second order. The multivalued operators are generated by the subdifferential of nonsmooth potentials that live in different spaces, $U$ and $V$, where in general $U \nsubseteq V$ and $V \nsubseteq U$. The proof is based on the regularization of the dissipation potential using the generalized Moreau--Yosida regularization and a semi-implicit time-discretization scheme, which demonstrates the existence of strong solutions to the regularized problem. The existence of solutions to the original problem is then shown by letting the regularization parameter converge to zero. Furthermore, we establish an energy-dissipation inequality for the solution. We conclude with applications of this abstract theory.
\end{abstract}
\vspace*{1em} \textbf{Keywords} Evolution equation of second order  $ \cdot $  Nonsmooth analysis $ \cdot $ Variational approximation scheme $ \cdot $ Subdifferential operator $ \cdot $ Calculus of Variations $ \cdot $ Non-convex energy  

\section{Introduction}
\subsection{Problem setting}
In the present work, we aim to prove the existence of strong  solutions to the following abstract \textsc{Cauchy} problem
\begin{align} 
\label{eq:I.2}
\begin{cases}
u''(t)+\partial\Psi(u'(t))+\partial \calE_t(u(t))+B(t,u(t),u'(t)) \ni f(t), &\quad t\in (0,T),\\
u(0)=u_0, \quad u'(0)=v_0,
\end{cases}
\end{align} 
 where $\Psi$ is the dissipation potential, $\calE_t$ the energy functional, $B$ the perturbation, and $f$ the external force. The functionals and operators are defined on suitable spaces, which are specified below. Here, the main assumptions are that the leading part of $\Psi$ is defined by a strongly positive, symmetric, and bounded bilinear form $a$, the energy functional $\calE_t$ is $\lambda$-convex, and the perturbation $B$ is a strongly continuous perturbation of $\partial \Psi$ and $\partial\calE_t$. Within the above-mentioned class of dissipation potentials, we consider the following two cases separately: in the first case (Case \textbf{(a)}), we assume that $\Psi(v)=a(v,v)$ and in the second case (Case \textbf{(b)}), we assume that $\Psi=\Psi_1+\Psi_2$, where $\Psi_1(v)=a(v,v)$ and $\Psi_2$ is a strongly continuous,   convex and nonsmooth perturbation that. Furthermore, we will specifically consider the case when $\calE_t$ is convex. We note that the energy functional and the dissipation potential are, in general, defined on different spaces. An illustrative example that satisfies all the assumptions above is given by
\begin{align*}
&\partial_{tt} u -\nabla\!\!\cdot\! \big(A \nabla \partial_t u \big ) + \nu \left \vert \partial_t u \right \vert^{q-2}\partial_t u+\nu\eta-\nabla \!\! \cdot\! \left (\vert \nabla u\vert^{p-2}\nabla u\right )+\nabla\!\!\cdot\! \xi+W'(u)+b(u,\partial_t u)=f,\\
&\eta(x,t)\in \mathrm{Sgn}(\partial_t u(x,t)), \quad \xi(x,t)\in \mathrm{Sgn}(\nabla u(x,t)) \quad \text{a.e. in } \Omega\times(0,T),
\end{align*} where  $p,q>1$ are to be chosen suitably, $\nu\geq 0$,  $A:\mathbb{R}^d\rightarrow \mathbb{R}^d$ is a linear, symmetric, and elliptic operator, $W:\mathbb{R}\rightarrow \mathbb{R}$ is a double-well potential given by $W(u)=\frac{1}{4}(u^2-1)^2$, $b:\mathbb{R}\rightarrow \mathbb{R}$ a lower order perturbation, and $f:\mathbb{R}\rightarrow \mathbb{R}$ an external force. The energy functional and the dissipation potential are given by 
\begin{align*}
\calE(u)=\int_{\Omega}\left(\frac{1}{p}\vert \nabla u(x)\vert^p+\vert \nabla u(x)\vert +\frac{1}{4}(u^2(x)-1)^2\right)\dd x
\end{align*} and
\begin{align*}
\Psi(v)=\int_{\Omega}\left( A(x)\nabla v(x)\cdot \nabla v(x)+ \frac{\nu}{q}\vert v(x)\vert^q+\nu\vert v(x)\vert\right)\dd x
\end{align*} and the perturbation is (formally) given by 
\begin{align*}
\langle B(u,v),w\rangle_{\rmL^2}=\int_{\Omega}b(u(x),v(x))w(x)\dd x.
\end{align*} We note that both functionals $\calE$ and $\Psi$ are nonsmooth and are, therefore, not eligible for the abstract results known in the literature. Moreover, the evolution equations of the second order, where the operators acting on $u$ and $u'$ are multi-valued, have not been studied before, even for concrete equations. Note that if $\nu=0$, we are in Case \textbf{(a)} and if $\nu>0$, we are in Case \textbf{(b)}. More, in particular, multi-valued applications will be discussed in Section \ref{se:app.DI1A} and \ref{se:app.DI1B}.\\\\
Evolution equations of second order occur in many applications, in particular in physics where many phenomena are modeled as such, e.g., the equations in peridynamics, the equations in visco-elastodynamics, the nonlinear \textsc{Klein--Gordon} equation from quantum mechanics and other nonlinear wave equations arising in mechanics and quantum mechanics,  equations for describing a vibrating membrane, in the see, e.g., \cite{EmmTha10CTDD,EmmSis11FDSO,EmmSis13EESO, RosTho17CDIP, BuMaRa12KVMG,Roub13NPDE}. Multivalued evolution equations (evolution inclusions) occur in many applications, e.g., physical phenomena where rate-independent responses of the body are typical, such as in plasticity \cite{MieRou15RIST}, in ferromagnetic hysteresis \cite{Visi85OLLE, MiRoSa13NADN} occurs or in Elasto-viscoplasticity \cite{RajRou03EDSA}. Applications are also found in optimal control theory \cite{AubCel84DI} or nonsmooth dynamical systems  \cite{Kunz00NSDS}.
\subsection{Literature review}
Evolution inclusions of second order of the form \eqref{eq:I.2} have been studied by \textsc{Rossi} \& \textsc{Thomas} in \cite{RosTho17CDIP} where $\partial \Psi=A:V\rightarrow V^*$ is a single-valued, linear, bounded, strongly positive, and symmetric operator defined on a reflexive and separable \textsc{Banach} space $V$, i.e., possessing a convex and quadratic potential. In their framework $\partial \calE_t$ is the subdifferential of a $\lambda$-convex functional with effective domain in a reflexive and separable \textsc{Banach} space $U\subset V$. In the framework of the \textsc{Gelfand} quintuplet
\begin{align*}
U\overset{d}{\hookrightarrow}V \overset{d}{\hookrightarrow}H\cong H^*\overset{d}{\hookrightarrow} V^*\overset{d}{\hookrightarrow} U^*
\end{align*} and under the assumption that $\calE_t$ satisfies a chain rule and $\partial\calE_t$ satisfies a closedness condition, they showed the existence of a strong solution. While this work is covered with the result presented here, we allow further a strongly continuous non-monotone and non-variational perturbation that depends on $u$ and $u'$ as well as a nonlinear monotone perturbation of $A$ of variational type. We also do not assume the rather restrictive assumption that $U\overset{d}{\hookrightarrow}V$. Furthermore, the strong closedness condition of $\partial \calE_t$ assumed in \cite{RosTho17CDIP} excludes the application to nonlinear elastodynamics where the operator satisfies a so-called \textsc{Andrews--Ball} type condition, see \cite{EmmSis13EESO}. In \textsc{Emmrich} \& \textsc{{\v{S}}i{\v{s}}ka} \cite{EmmSis13EESO}, the authors develop an abstract theory in the smooth setting with the application to nonlinear elastodynamics. They prove the existence of strong solutions for the case where $A:V_A\rightarrow V_A^*$ is linear, bounded, strongly positive and symmetric, and $B:V_B\rightarrow V_B^*$ is supposed to be demicontinuous and a bounded potential operator. In addition, $B$ satisfies an \textsc{Andrews--Ball}-type condition, meaning that $(B+\lambda A):V\rightarrow V^*$ is monotone where $V\vcentcolon=V_A \cap V_B$ is densely and continuously embedded into the separable and reflexive \textsc{Banach} spaces, $V_A$ and $V_B$, for which we assume not that either of the two spaces is continuously embedded in the other one. Since we allow a more general nonsmooth functional $\calE$, this result is also covered by our main result. Another motivation for this work is to complement the results obtained in B. \cite{Bach22NDIS} where the principal part of the operator acting on $u'$ is nonlinear and multi-valued, and the principal part of the operator acting on $u$ is linear, symmetric, and positive, see also \cite{Bach21ONSA}. Hence, our contribution to the literature concerns the following:
\begin{itemize}
\item We allow the functionals that are acting on $u$ and $u'$ to be nonsmooth, hence generating multi-valued subdifferentials in the differential inclusion. This has not even been considered in the non-abstract setting.
\item We allow the multi-valued operators to live on different spaces.
\item We allow non-variational and non-monotone perturbations of the subdifferential operators. 
\item We provide in the applications an existence result to concrete problems for which there are no results known and can not be obtained by other abstract results.
\item The existence result is based on a numerical scheme which can be used to obtain numerical approximations.
\end{itemize}
The case where the principal part of the operator $A$ is nonlinear has also been discussed by several authors, see, e.g., \cite{LioStr65SNLE, EmSiTh15FDNS,Bach22NDIS}.\\
Doubly nonlinear evolution equations where the leading parts of $A$ and $B$ are both nonlinear and contain the same order of spatial derivatives in the applications are unfortunately not covered by our result. However, in the \textsc{Hilbert} space setting, the well-posedness of doubly nonlinear evolution equations of second order under a \textsc{Lipschitz} condition on the operator $B$ has been shown in \cite{Bach22WPFN}. In certain other situations, this can also be demonstrated by exploiting the special structure of the operators; see, e.g., \cite{Puhs15OEVF,FriNec88SNWE,BuMaRa12KVMG,BuKaSt1GKVM}..\\\\
For further results on nonlinear evolution equations, we refer to \textsc{Leray} \cite{Lera53HDE},  \textsc{Dionne} \cite{Dion62SPCH},  \textsc{Emmrich \& Thalhammer} \cite{EmmTha10CTDD,EmmTha11DNEE}, \textsc{Emmrich, \& {\v{S}}i{\v{s}}ka} \cite{EmmSis11FDSO} including stochastic perturbations, \textsc{Emmrich,\textsc{{\v{S}}i{\v{s}}ka} \& Thalhammer} \cite{EmSiTh15FDNS} for a numerical analysis, \textsc{Emmrich, \textsc{{\v{S}}i{\v{s}}ka} \& Wr\'{o}blewska-Kami\'{n}ska} \cite{EmSiWr16ESOO} and \textsc{Ruf} \cite{Ruf17CFDS} for results on \textsc{Orlicz} spaces, and the monographs \textsc{Lions} \cite{Lion69QMRP}, \textsc{Lions \& Magenes} \cite[Chapitre 3.8]{LioMag72NHBV}, \textsc{Barbu} \cite[Chapter V]{Barb76NSDE}, \textsc{Wloka} \cite[Chapter V]{Wlok82PD}, \textsc{Zeidler} \cite[Chapter 33]{Zeid90NFA2b},  \textsc{Roub{\'\i}{\v{c}}ek} \cite[Chapter 11]{Roub13NPDE} and the references therein.\\ 
The list of literature presented in this section is not intended to be exhaustive.
\subsection{Organization of the paper}
The paper is organized as follows. In Section 2, we set the analytical framework and briefly introduce some notions and results from the theory of subdifferential calculus and convex analysis. In Section 3, we present and discuss the assumptions on the dissipation potential $\Psi$, the energy functional $\calE$ and the perturbation $B$ as well as the external force $f$, and we state the main result. Section 4 is devoted to the proof and in Section 5 we show some examples which does not fit into the abstract framework of the cited authors.

\section{Topological assumptions and preliminaries}
\label{se:AssumExistRe.1}
\subsection{Function space setting}
In the following, let $(U,\Vert\cdot\Vert_U), (V,\Vert\cdot\Vert_V)$, $(W,\Vert\cdot\Vert_W)$ and $(\widetilde{W},\Vert\cdot\Vert_{\widetilde{W}})$ be real, separable, and reflexive  \textsc{Banach} spaces such that the dual space $W^*$ is uniformly convex. Furthermore, let $(H,\vert\cdot\vert,(\cdot,\cdot))$ be a \textsc{Hilbert} space with norm $\vert
\cdot\vert$ induced by the inner product $(\cdot,\cdot)$. Then, we assume the dense, continuous, and compact embeddings
\begin{align*}
\begin{cases}
U\cap V \overset{d}{\hookrightarrow}  U\overset{c,d}{\hookrightarrow} \widetilde{W} \overset{d}{\hookrightarrow} H \cong H^*\overset{d}{\hookrightarrow} \widetilde{W}^*\overset{d}{\hookrightarrow} U^* \overset{d}{\hookrightarrow}V^*+U^* \\
 U\cap V \overset{d}{\hookrightarrow} V\overset{c,d}{\hookrightarrow} W \overset{d}{\hookrightarrow} H \cong H^*\overset{d}{\hookrightarrow} W^*\overset{d}{\hookrightarrow} V^*\overset{d}{\hookrightarrow}V^*+U^*
 \end{cases}
\end{align*} and if the perturbation does not explicitly depend on $u$ or $u'$, then we do not need to assume $U\overset{c}{\hookrightarrow} \widetilde{W}$ or $V\overset{c}{\hookrightarrow} W$, respectively, but instead that $V \overset{c}{\hookrightarrow} H$. We stress that we neither assume $U\hookrightarrow V$ nor $V \hookrightarrow U$. The spaces can coincide if a certain embedding is not assumed to be compact. For instance, the cases $V=U$, $\WW=H$ or $W=H$ are admissible. Introducing the spaces $W$ and $\widetilde{W}$ allows us to make use of the finer structure of the spaces, which enables us to treat additional nonlinearities of lower order. As examples for the appearing spaces, we can think of the \textsc{Sobolev} spaces $U=\rmW^{k,p}(\Omega),V=\rmH^{l}(\Omega)$ and the \textsc{Lebesgue} spaces $W=\rmL^q(\Omega)$ and $H= \rmL^2(\Omega)$ or $U=\rmW^{k,p}(\Omega),V=\rmW^{s,p}(\Omega),\, W=\rmH^l(\Omega)$ and $H= \rmL^2(\Omega)$ for suitably chosen numbers $k,l\in \mathbb{N}$ and real values $s ,p>0, q>1$. 
\subsection{Facts from convex analysis}
Before we present the precise assumptions on the functionals and the operators, we recall some facts and results from functional and convex analysis. 
For a proper functional $F:X\rightarrow (-\infty,+\infty]$ defined on a \textsc{Banach} $(X,\Vert \cdot \Vert_X)$ space, the $($\textsc{Fr\'{e}chet}$)$ subdifferential of $F$ is given by the multivalued map $\partial F:X\rightrightarrows X^*$ with
\begin{align*}
\partial F(u):=\left\lbrace \xi\in X^*: \liminf_{v\rightarrow u} \frac{F(v)-F(u)-\langle \xi,v-u\rangle_{X^*\times X}}{\Vert v-u\Vert_X} \geq 0\right \rbrace,
\end{align*} where $\langle\cdot,\cdot\rangle$ denotes the duality pairing between the \textsc{Banach} space $X$ and its topological dual space $X^*$. The elements of the subdifferential are also called subgradients. The effective domain of $F$ and the domain of its subdifferential $\partial F$, we denote by $\DOM(F):=\lbrace v\in X \mid F(v)<+\infty \rbrace$ and $\DOM(\partial F):=\lbrace v\in X : \partial F(v)\neq \emptyset\rbrace$, respectively. If the set of subgradients of $F$ at a given point $u$ is nonempty, we say that $F$ is subdifferentiable at $u$. The following lemma states that the subdifferential of the sum of a subdifferential function and a \textsc{G\^{a}teaux} differentiable function equals the sum of the subdifferentials of both functions. More precisely, there holds
\begin{lem} \label{le:Subdif} Let $F_1:X\rightarrow (-\infty,+\infty]$ and $F_2:X\rightarrow (-\infty,+\infty]$ be subdifferentiable and \textsc{Fr\'{e}chet} differentiable at $u\in \DOM(\partial F_1)\cap \DOM(D F_2)\neq \emptyset$, respectively. Then, there holds
\begin{align*}
\partial(F_1+F_2)(u)=\partial F_1(u)+D F_2(u),
\end{align*} where $D F_2$ denotes the \textsc{Fr\'{e}chet} derivative of $F_2$.
\end{lem}
\begin{proof}
This immediately follows the definition of a subdifferential.
\end{proof}

Since we are dealing with convex and $\lambda$-convex functions, we would like to express the subdifferential in an equivalent way that is easier to handle. The function $F$ is called $\lambda$-convex with parameter $\lambda\in \mathbb{R}$ if 
\begin{align*}
F(tu+(1-t)v)\leq tF(u)+(1-t)F(v)+\lambda t(1-t)\Vert u-v \Vert_X^2
\end{align*} for all $u,v\in D(F)$ and $t\in (0,1)$. If $\lambda<0$ and $\lambda=0$, then the function is called strongly convex and convex, respectively.
Employing the definition of a subdifferential and the $\lambda$-convexity, it is easy to show that for a $\lambda$-convex and proper function $F$, the subdifferential of $F$ is equivalently given by
\begin{align}\label{eq:subdiff.conv}
\partial F(u)=\left\lbrace \xi\in X^*: F(u)\leq F(v)+\langle \xi,u-v\rangle_{X^*\times X}+\lambda \Vert u-v\Vert_X^2 \quad \text{ for all  } v\in X \right\rbrace.
\end{align} It follows immediately that 
\begin{align*}
\langle \xi-\zeta, u-v\rangle_{X^*\times X}\geq -\lambda \Vert u-v\Vert_X^2
\end{align*} for all  $\xi \in \partial F(u)$ and $\zeta \in \partial F(v)$. Hence, for $\lambda<0$, the subdifferential of a $\lambda$-convex function $F$ is strongly monotone.  This gives rise to the following lemma.
\begin{lem}[Variational sum rule] \label{le:Subdif2} Let $F_1:X\rightarrow (-\infty,+\infty]$ and $F_2:X\rightarrow (-\infty,+\infty]$ be proper, lower semicontinuous and convex, and if there is a point $\tilde{u}\in \DOM(F_1)\cap \DOM(F_2)$ where $F_2$ is continuous, we have
\begin{align}
\partial (F_1+F_2)(v)=\partial F_1(v)+\partial F_2(v) \quad \text{for all }v\in X.
\end{align}
If $F_2$ is, in addition \textsc{G\^{a}teaux} differentiable on $V$, there holds $\partial F_2(v)= D_G F_2(v)$ and we have
\begin{align*}
\partial(F_1+F_2)(v)=\partial F_1(v)+D_G F_2(v) \quad \text{for all } v\in X,
\end{align*} where $D_G F_2$ denotes the \textsc{G\^{a}teaux} derivative of $F_2$.
\end{lem}
\begin{proof}
Proposition 5.3. on p. 23 and Proposition 5.6 on p. 26 in \cite{EkeTem76CAVP}.
\end{proof}
Let us now introduce an important tool from the theory of convex analysis. For a
proper, lower semicontinuous and convex function $F: X\rightarrow (-\infty,+\infty]$, we define the so-called convex conjugate (or \textsc{Legendre--Fenchel} transform) $F^*:X^*\rightarrow (-\infty,+\infty]$ by 
\begin{align*}
  F^*(\xi):=\sup_{u\in V}\left \lbrace \langle \xi,u\rangle -
    F(u)\right \rbrace, \quad \xi\in X^*.
\end{align*} By definition, we directly obtain the \textsc{Fenchel--Young}
inequality
\begin{align*}
 \langle \xi,u\rangle \leq  F(u)+F^*(\xi), \quad v\in X, \xi\in X^*.
\end{align*} It can be checked that the convex conjugate itself is proper, lower semicontinuous and convex, see, e.g., \textsc{Ekeland} and \textsc{T\'emam} \cite{EkeTem76CAVP}. If,  in addition, we assume $F(0)=0$, then $F^*(0)=0$ holds as well. We may ask how the convex conjugate of the sum of two functions can explicitly be expressed in terms of the two functions.
\begin{lem}\label{le:Bre.Att}
 Let $F_1:X\rightarrow (-\infty,+\infty]$ and $F_2:X \rightarrow (-\infty,+\infty]$ be a proper, lower semicontinuous and convex functional such that 
 \begin{align*}
 \bigcap_{\lambda\geq 0} \lambda \left(\DOM(F_1)-\DOM(F_2) \right) \quad \text{ is a closed vector space.}
 \end{align*} Moreover, let $F_1^*, F_2^*:X^*\rightarrow (-\infty,+\infty]$ be the associated convex conjugate of $F_1$ and $F_2$, respectively. Then, there holds
\begin{align}\label{eq:brez.formula}
(F_1+F_2)^*(\xi)= \min_{\eta\in X^*}\left( F_1^*(\xi-\eta)+F_2^*(\eta)\right) \quad \text{for all }\xi \in X^*.
\end{align}
\end{lem}
\begin{proof} Theorem 1.1, pp. 126, in \textsc{Br\'{e}zis} u. \textsc{Attouch} \cite{AttBre86DSVF}.
\end{proof} For an illustration of the lemma, we consider
\begin{exa}\label{exa.Brez} Let $(X,\Vert \cdot \Vert_{X})$ and $(Y, \Vert \cdot \Vert_Y)$ be two \textsc{Banach} spaces such that both $X$ and $Y$ are continuously embedded into another \textsc{Banach} space $Z$, and such that $X\cap Y$, equipped with the norm $\Vert \cdot \Vert_{X\cap Y}= \Vert \cdot \Vert_{X}+\Vert \cdot \Vert_Y$, is dense in both $X$ and $Y$. Then, the space $X\cap Y$ becomes a \textsc{Banach} space itself and the dual space can be identified as $X^*+Y^*$ with the dual norm $\Vert \xi \Vert_{X^*+Y^*}=\inf_{\overset{\xi_1\in X^*, \xi_2\in Y^*}{\xi=\xi_1+\xi_2}}\max \lbrace \Vert \xi_1\Vert_{X^*}, \Vert\xi_2\Vert_{Y^*}\rbrace$, see, e.g., Chapter I Section 5 in \cite{GaGrZa74NOOD}.  Now, let for $p,q\in (1,+\infty)$ the functions $F_1, F_2:X\cap Y\rightarrow \mathbb{R}$ be given by
\begin{align*}\label{eq:rock.conv.conj.}
F_1(u)=\frac{1}{p}\Vert u\Vert_{X}^p, \quad F_2(u)=\frac{1}{q}\Vert u\Vert_{Y}^q, \quad u\in X\cap Y.
\end{align*}  Then, according to Lemma \ref{le:Bre.Att}, the convex conjugate $(F_1+F_2)^*: X^*+Y^*\rightarrow (-\infty,+\infty]$ given by
\begin{align}
(F_1+F_2)^*(\xi)=\min_{\overset{\xi_1\in X^*, \xi_2\in Y^*}{\xi=\xi_1+\xi_2}}\left( \frac{1}{p'}\Vert \xi_1 \Vert_{X^*}^{p'}+\frac{1}{q'}\Vert \xi_2 \Vert_{Y^*}^{q'}\right) \quad \text{for all } \xi\in X^*+Y^*,
\end{align} where $p'>1$ and $q'>1$ denote the conjugate exponent of $p$ and $r$, respectively, i.e., fulfilling $1/p+1/p'=1$ and $1/q+1/q'=1$.  Applying \textsc{Young}'s inequality to \eqref{eq:rock.conv.conj.}, we obtain the estimates
\begin{align*}
(F_1+F_2)(u)&\geq  C\Vert u\Vert_{X\cap Y}-C \quad \text{for all } u\in X\cap Y\\
(F_1+F_2)^*(\xi)&\geq C\Vert \xi \Vert_{X^*+Y^*}-C \quad \text{for all } \xi \in X^*+Y^*
\end{align*} for some constant $C>0$. We will make use of the latter estimates later on by choosing particularly $X=\rmL^2(0,T;V)$ and $Y=\rmL^{r}(0,T;W)$.
\end{exa}
 The following lemma shows that there is, in fact, a relation between the subgradient of a function and its convex conjugate.
 
\begin{lem}\label{le:Leg.Fen}
 Let $V$ be a \textsc{Banach} space and let $F:V\rightarrow (-\infty,+\infty]$ be a proper, lower semicontinuous and convex functional and let $F^*:V^*\rightarrow (-\infty,+\infty]$ be the convex conjugate of $F$. Then for all $(u,\xi)\in V\times V^*$, the following assertions are equivalent:
\begin{itemize}
\item[$i)$]$\xi\in \partial F(u) \quad \text{in } V^*;$
\item[$ii)$] $u\in \partial F^*(\xi)\quad \text{in } V;$
\item[$iii)$]$\langle \xi, u\rangle=F(u)+F^*(\xi) \quad \text{in
  } \mathbb{ R}.$
\end{itemize}
\end{lem}
\begin{proof} Proposition 5.1 and Corollary 5.2 on pp. 21 in \cite{EkeTem76CAVP}.
\end{proof}
The next result shows a certain chain rule for the subdifferential. 
\begin{lem} \label{le:chain.rule} Let $U$ and $V$ be \textsc{Banach} spaces. Let $\Lambda:U\rightarrow V$ be a linear, bounded operator and $f:U\rightarrow (-\infty,+\infty]$ be a proper, lower semicontinuous, and convex functional. If there exists a point $\Lambda \tilde{u}\in V$ with $\tilde{u}\in U$, where $f$ is finite and continuous, then for all $u\in U$, there holds 
\begin{align*}
\partial (f\circ \Lambda)(u) = \Lambda^* \partial f(\Lambda u) \quad \text{for all }u\in U,
\end{align*} where $\Lambda^*:V^*\rightarrow U^*$ denotes the adjoint operator of $\Lambda$.
\end{lem}
\begin{proof}[Proof]
This has been proven in \textsc{Ekeland \& Temam} \cite[Proposition 5.7]{EkeTem76CAVP}.
\end{proof}
Finally, we recall that the space $U\cap V$ equipped with the norm $\Vert \cdot \Vert_{U\cap V}=\Vert \cdot \Vert_{U}+\Vert \cdot \Vert_{V}$ is a separable and reflexive \textsc{Banach} space and the dual space is given by $(U\cap V)^*=U^*+V^*$ with the norm $\Vert \xi \Vert_{U^*+V^*}=\inf_{\overset{\xi_1\in U^*, \xi_2\in V^*}{\xi=\xi_1+\xi_2}}\max \lbrace \Vert \xi_1\Vert_{U^*}, \Vert\xi_2\Vert_{V^*}\rbrace$, see Example \ref{exa.Brez}. Furthermore, the duality pairing between $U\cap V$ and $U^*+V^*$ is given by
\begin{align*}
\langle f,v\rangle_{(U^*+ V^*)\times (U\cap V)}=\langle f_1,u\rangle_{U^*\times U}+\langle f_2,u\rangle_{V^*\times V}, \quad  u\in U\cap V,
\end{align*} for all $v\in U\cap V$ and any partition $f=f_1+f_2$ with $f_1\in U$ and $f_2\in V$. Second, for any $p\in [1,+\infty]$, there holds $\rmL^p(0,T;U)\cap \rmL^p(0,T;V)=\rmL^p(0,T;U\cap V)$, where the measurability immediately follows from the \textsc{Pettis} theorem, see, e.g., \textsc{Diestel \& Uhl} \cite[Theorem 2, p. 42]{DieUhl77VEME}. And third, for the \textsc{Banach} spaces $X,Y\in \lbrace U\cap V,U,V,W,\WW,H\rbrace$ satisfying the embedding $X\hookrightarrow Y$, there holds
\begin{align*}
\langle f,v\rangle_{X^*\times X}=\langle f,v\rangle_{Y^*\times Y} \quad \text{if } v\in X \text{ and } f\in Y^*.
\end{align*} see, e.g, \textsc{Br\'{e}zis} \cite[Remark 3, pp. 136]{Brez11FASS} and \textsc{Gajewski} et al. \cite[Kapitel 1, \S  5]{GaGrZa74NOOD}.
\section{Main result}
\subsection{Assumptions on the functionals and operators}

In this section, we summarize all the assumptions regarding the functionals $\Psi$ and $\mathcal{E}$, and the operators $B$ and $f$. Since the subdifferential of the main part of $\Psi$ is linear, we refer henceforth to the inclusion \eqref{eq:I.2} in the given framework as \textit{linearly damped inertial system} $(U,V,W,\WW,H,\calE,\Psi,B,f)$. The assumptions we make for the linearly damped inertial system resemble the structure of those we made for the perturbed gradient system in \cite{BaEmMi19EREG} where the same evolution inclusion has been investigated after neglecting the inertial term  $u''(t)$. Involving inertia makes the situation much more delicate. As a consequence, we will impose, in general, stronger conditions on the functionals and operators in order to ensure the solvability of the problem. Hereinafter, we collect the assumptions for the dissipation potential $\Psi$  and remind the reader that we distinguish two cases \textbf{(a)} and \textbf{(b)}. 

\begin{enumerate}[label=(\thesection.$\Uppsi$), leftmargin=3.2em] \label{eq:Psi1}
\item \label{eq:Psi1.1} \textbf{Dissipation potential.} \\
Case \textbf{(a)}: we assume that there exists a strongly positive, symmetric, and continuous bilinear form  $a: V\times V\rightarrow \mathbb{R}$ such that $\Psi(v)=\frac{1}{2}a(v,v)$, i.e., there is a constant $\mu>0$ such that
\begin{align} \label{eq:growth.Psi}
\mu\Vert v\Vert_V^2   \leq \Psi(v)  \quad \text{for all }v\in V.
\end{align}
Case \textbf{(b)}: we assume that $\Psi=\Psi_1+\Psi_2$, where $\Psi_1(v)=\frac{1}{2}a(v,v)$ with the bilinear form $a:V\times V\rightarrow \mathbb{R}$ as above and $\Psi_2:\Wt\rightarrow \mathbb{R}$ to be a lower semicontinuous and convex functional with $\Psi_2(0)=0$ satisfying the following growth condition: there exists a positive number $q>1$ and  constants $\hat{c},\hat{C}>0$ such that 
\begin{align}\label{eq:Psi.growth}
\hat{c}(\Vert v\Vert^q_{\Wt}-1)\leq \Psi_2(v)\leq \hat{C}(\Vert v\Vert^q_{\Wt}+1) \quad \text{for all }v\in \Wt.
\end{align} In addition, we assume for the subgradients of $\Psi_2$ the following growth condition: for all $R>0$, there exists a positive real constant $C_R>0$ such that 
\begin{align}\label{eq:Psi.Gateaux}
\begin{split}
\Vert \eta \Vert_{W^*}&\leq C_R(1+\Vert v\Vert_{\Wt}^{q-1}) \quad \text{for all }\eta\in \partial \Psi_2(v)\text{ and } v\in W \text{ with }\vert v\vert\leq R.
\end{split}
\end{align} 
\end{enumerate}

\begin{rem}\label{re:Assump.Psi1} \mbox{} \vspace{-0.6em}
\begin{itemize}
\item[$i)$] Assumption \ref{eq:Psi1.1} yields the convexity and continuity of the mapping $v\mapsto \Psi(v)$. Furthermore, $\Psi$ is \textsc{G\^{a}teaux} differentiable with the \textsc{G\^{a}teaux} derivative given by a positive, linear bounded and symmetric operator $A:V\rightarrow V^*$ such that $ $ $\partial \Psi(v)=\lbrace A v\rbrace $ and the potential can be expressed by $\Psi(v)=\frac{1}{2}\langle Av,v\rangle$. Assumption \ref{eq:Psi1.1} implies that the
  \textsc{Legendre--Fenchel} transform $\Psi^*$ of $\Psi$ is convex, continuous, finite everywhere,
  i.e., $\DOM(\Psi^*)=V^*$, and can be explicitly expressed by $\Psi^*(\xi)=\frac{1}{2}\langle \xi,A^{-1}\xi\rangle$, where $A^{-1}:V^*\rightarrow V$ is also continuous, symmetric and positive, which follows from the \textsc{Lax--Milgram} theorem, see, e.g., \textsc{Br\'{e}zis} \cite[Corollary 5.8, p. 140]{Brez11FASS}.
\item[$ii)$] From the properties of the conjugate, we obtain from \eqref{eq:growth.Psi} and \eqref{eq:Psi.growth} the following growth condition for the conjugates $\Psi_1^*:V^*\rightarrow \mathbb{R}$ and $\Psi_2^*:V^*\rightarrow (-\infty,+\infty]$: there exist positive constants $\bar{c},\bar{C}>$ such that 
\begin{align}\label{eq:Psi*.growth}
\begin{split}
\bar{c} \Vert \xi \Vert_{V^*}^2\leq &\Psi_1^*(\xi)\leq \bar{C}\Vert \xi \Vert_{V^*}^2\quad \text{for all }\xi\in V^*,\\
\begin{rcases}
\bar{c}(\Vert \xi \Vert^{q^*}_{\Wt^*}-1)\\
+\infty 
\end{rcases} \leq &\Psi_2^*(\xi)\leq 
\begin{cases} 
 \bar{C}(\Vert \xi\Vert^{q^*}_{\Wt^*}+1) \quad &\text{if } \xi\in \Wt^*\\
+\infty \quad &\text{otherwise },
\end{cases}
\end{split}
\end{align}
 where $q^*>1$ denotes the conjugate exponent of $q$. To justify the formula \eqref{eq:Psi*.growth}, it is not restrictive to show it for $\Psi_2(v)=\frac{1}{q}\Vert v\Vert^q_{\Wt}, v\in V$. To do so, we remark that the conjugate function $f^*$ of any proper, convex, and lower semicontinuous function $f:V\rightarrow \overline{\mathbb{R}}$ is also proper, convex, and lower semicontinuous, and that $f^{**}=f$, see, e.g., \textsc{Ekeland \& Temam} \cite[Section 4.1, pp. 16]{EkeTem76CAVP}. Thus, defining $\widetilde{\Psi}:V^*\rightarrow \overline{\mathbb{R}}$ via
\begin{align*}
\widetilde{\Psi}(\xi)=
\begin{cases} 
 \frac{1}{q^*}\Vert \xi\Vert^{q^*}_{\Wt^*} \quad &\text{if } \xi\in \Wt^*\\
+\infty \quad &\text{otherwise },
\end{cases}
\end{align*} it follows that $\widetilde{\Psi}$ is a proper, convex, and lower semicontinuous function on $V^*$, which easily follows from the fact that a function is convex and lower semicontinuous if and only if its epigraph is convex and closed \cite{Zali02CAGV}. Then, we show that $\widetilde{\Psi}^*=\Psi_2$ which in turn implies $\Psi_2^{*}=\widetilde{\Psi}=\widetilde{\Psi}^{**}$ where the first equality follows from 
\begin{align*}
\widetilde{\Psi}^{*}(v)&=\sup_{\xi \in V^*}\left\lbrace \langle \xi,v\rangle_{V^*\times V}-\widetilde{\Psi}^{*}(\xi)\right\rbrace\\
&=\sup_{\xi \in \Wt^*}\left \lbrace \langle \xi,v\rangle_{\Wt^*\times \Wt}-\frac{1}{q^*}\Vert \xi\Vert^{q^*}_{\Wt^*}\right \rbrace\\
&=\frac{1}{q}\Vert v\Vert^q_{\Wt}\\
&=\Psi_2(v) \quad \text{for all }v\in \Wt,
\end{align*} where we have used that $(\frac{1}{q^*}\Vert \cdot\Vert^{q^*}_{\Wt})^*=\frac{1}{q}\Vert v\Vert^q_{\Wt}$ on $W$ \cite[Remark 4.1., pp. 19]{EkeTem76CAVP}
\item[$iii)$] We remark that we could also allow for a time-dependent dissipation potential $\Psi_t=\Psi^1_t+\Psi^2_t$ when we assume that $t\mapsto a(t,u,v)\in \rmC([0,T])\cap\rmC^1(0,T)$ for all $u,v\in V$ and a strong monotonicity and boundedness of $A(t):V\rightarrow V^*$ uniformly in time as well as a slight modification of Assumption \ref{eq:cond.E1.3} and \ref{eq:B1.2}, whereas for $\Psi_t^2$ we would assume that for all $t\in [0,T]$, the functional $\Psi^2_t$ is lower semicontinuous and convex satisfying the Conditions \eqref{eq:Psi.growth} and \eqref{eq:Psi.Gateaux} uniformly in time. For simplicity, we will not consider this case here.
\end{itemize}
\end{rem} \noindent
We proceed with collecting the assumptions for the energy functional $\calE$. To do so, we define $V_\lambda=U$ if $\lambda=0$ and $V_\lambda=U\cap V$ if $\lambda>0$. We make this distinction because for the convex case, i.e., when $\lambda=0$, we will obtain a stronger result meaning that the initial value $u_0$ can be chosen to be in $\DOM(\calE_t)$ instead of $\DOM(\calE_t)\cap V$ as in the $\lambda$-convex case with $\lambda\neq 0$, and that the subgradient of $\calE_t$ is in $U^*$ instead of $U^*+V^*$, see Theorem \ref{th:MainExist1}.

\begin{enumerate}[label=(\thesection .E\alph*),
 leftmargin=3.2em]  \label{eq:cond.E1}
\item \label{eq:cond.E1.1} 
\textbf{Lower semicontinuity.} For all $t\in[0,T]$, the
  functional $\calE_t:U \rightarrow (-\infty,+\infty]$ is proper and sequentially weakly lower semicontinuous with time-independent effective domain $D:=\DOM(\calE_t)\subset U$ for all $t\in [0,T]$. Furthermore, the set $D\cap V$ is dense in $D$ in the topology of $U$, and if $\calE_t$ is convex, the interior of $D$ is non-empty.
\item \label{eq:cond.E1.2}
 \textbf{Bounded from below.} $\calE_t$ is bounded from below uniformly in time, i.e., there exists a constant $C_0\in \mathbb{R}$ such that 
  \begin{align*}
  \calE_t(u)\geq C_0 \quad \text{for all } u\in U\text{ and }t\in[0,T].
  \end{align*} Since the potential of a function is, if it exists, unique up to a constant, we assume without loss of generality $C_0=0$.
\item \label{eq:cond.E1.3} \textbf{Coercivity.}  For every $t\in[0,T]$, $\calE_t$ has bounded sublevel sets in $U$.

\item \label{eq:cond.E1.4} \textbf{Control of the time derivative.} For all $u\in U$,
  the mapping $t\mapsto \calE_t(u)$ is in $\rmC([0,T])\cap \rmC^1(0,T)$ and its derivative $\partial_t\calE_t$ is
  controlled by the function $\calE_t$, i.e., there exists $C_1>0$ such that
\begin{align*}
  \vert \partial_t \calE_t(u)\vert \leq C_1 \calE_t(u)\quad \text{for all } t\in
  (0,T) \text{ and } u\in U.
\end{align*} 
\item \label{eq:cond.E1.7}
 \textbf{$\lambda$-convexity.} 
 There exists $\lambda \geq 0$ such that for every $t\in[0,T]$, the energy functional $\calE_t$ is $\lambda$-convex on $V$ (by extending $\calE$ on $V$), i.e., for all $u,v\in U\cap V $ and $\vartheta\in (0,1)$, there holds
\begin{align}\label{eq:lambda.convex.2}
\calE_t(\vartheta v +(1-\vartheta)u)\leq \vartheta\calE_t(v)+(1-\vartheta)\calE_t(u)+\lambda \vartheta(1-\vartheta)\Vert v-u \Vert_{V}^2.
\end{align}

\item \label{eq:cond.E1.6} \textbf{Closedness of $\mathrm{Gr}(\partial\calE)$.} For all sequences of measurable functions $(\tf_n)_{n\in \mathbb{N}}$ with $\tf_n:[0,T]\rightarrow [0,T],\,n\in \mathbb{N}$, $ (u_n)_{n\in \mathbb{N}}$, $(\xi_n)_\nN$, and measurable functions $u,\xi$ satisfying 
\begin{enumerate}[label=\alph*)]
\item $\tf_n(t)\rightarrow t \, \text{ for a.e. t}\in(0,T),$ as $n\rightarrow \infty$,
\item $\exists \, C_2>0: \, \sup_{\nN,t\in[0,T]} \calE_t (u_n(t))\leq C_2,$
\item $\xi_n(t)\in \partial_{V_\lambda} \calE_{\tf_n(t)}(u_n(t)) \text{ a.e. in } (0,T), n\in \mathbb{N}$,
\item $u_n-\tilde{u}_0\overset{*}{\rightharpoonup} u-\tilde{u}_0 \text{ in } \rmL^{\infty} (0,T;V_\lambda)$ and $u_n-\tilde{u}_0 \rightarrow  u-\tilde{u}_0 \text{ in } \rmL^{2} (0,T;V)$ for any $\tilde{u}_0\in D$ and $\xi_n \rightharpoonup \xi \text{ in } \rmL^2(0,T;V_\lambda^*)$ as $n\rightarrow \infty$. Additionally, there exists a constant $C_3>0$ such that for sufficiently small $h>0$, there holds\\\\
Case \textbf{(a)}: 
\begin{align}\label{LAD.lemma.a}
\sup_{n\in \mathbb{N}}\Vert \sigma_h u_n-u_n\Vert_{\rmL^2(0,T-h,V)}\leq C_3 h  
\end{align}
Case \textbf{(b):}
\begin{align}\label{LAD.lemma.b}
\sup_{n\in \mathbb{N}}\Vert \sigma_h u_n-u_n\Vert_{\rmL^2(0,T-h,V)\cap \rmL^{r}(0,T-h,W)}\leq  C_3 h,
\end{align} 
 where $\sigma_h v:=\chi_{[0,T-h]}v(\cdot+h)$ for any function $v:[0,T]\rightarrow V$,
\item $\limsup_{n\rightarrow \infty} \int_0^T \langle \xi_n(t)-\xi(t),u_n(t)-u(t)\rangle_{V_\lambda^*\times V_\lambda} \dd t \leq 0$,
\end{enumerate}
 we have the relations
\begin{align*} 
\xi(t) \in \partial_{V_\lambda} \calE_t(u(t))\subset V_\lambda^*, \quad \calE_{\tf_n(t)}(u_n(t))\rightarrow \calE_t(u(t)) \quad \text{as } n\rightarrow \infty\\
\text{ and}\,\,\,\, \limsup_{n\rightarrow \infty} \partial_t  \calE_{t}(u_n(t))\leq \partial_t\calE_t(u(t))\quad \text{for a.e. }t\in (0,T).
\end{align*} 
\item \label{eq:cond.E1.8} \textbf{Control of the subgradient.}
  There exist constants $C_4>0$ and $\sigma>0$ such that
\begin{align*}
\Vert \xi \Vert_{V_\lambda^*}^\sigma \leq C_4 (1+\calE_t(u)+\Vert u\Vert_{V_\lambda}) \quad \forall t\in [0,T], u \in D(\partial_{V_\lambda} \calE_t),\, \xi \in \partial_{V_\lambda} \calE_t(u).
\end{align*} 
\end{enumerate}
Henceforth, we write (3.E) when we refer to the Conditions \ref{eq:cond.E1.1}-\ref{eq:cond.E1.8}. Before we proceed with the conditions for the perturbation $B$ and the external force $f$, we first give a few relevant comments on these assumptions that will be important later on. 

\begin{rem}\label{re:Assump.E1} \mbox{}\vspace{-0.6em}

\begin{itemize}

\item[$i)$] From Assumption \ref{eq:cond.E1.4}, we deduce with \textsc{Gronwall}'s lemma, see, e.g., \cite[Lemma 7.3.1, pp.180]{Emmr04GUOD}, the chain of inequalities
\begin{align*}
\ee^{-C_1\vert t-s\vert}\calE_s(u)\leq \calE_t(u)\leq \ee^{C_1\vert
  t-s\vert}\calE_s(u) \quad \text{ for all }s,t\in [0,T],\, u\in D.
\end{align*} 
In particular, there holds
\begin{align*}
 \sup_{t\in [0,T]}\calE_t(u)\leq \ee^{C_1T} \inf_{t\in[0,T]} \calE_t(u)
  \quad \text{for all } u\in D.
\end{align*} 
\item[$ii)$] In Case \textbf{(b)} it is possible to improve the assumption of $\lambda$-convexity in the following way: there exist positive constants $\lambda_1, \lambda_2> 0$ such that
\begin{align*}
\calE_t(\vartheta u+(1-\vartheta )v)\leq &\vartheta \calE_t(u)+(1- \vartheta)\calE_t(v)\\
&+\vartheta(1-\vartheta)\left(\lambda_1\Vert u-v \Vert_V^2+\lambda_2\Psi(u-v)^{\frac{1}{q}}\vert u-v\vert\right)
\end{align*} for all $u\in D, v\in V$ and $\vartheta\in (0,1)$, where $q>1$ comes from Assumption \ref{eq:Psi1.1}
\end{itemize} 
\end{rem}

Finally, we present the assumptions on the non-variational non-monotone perturbation $B$ and the external force $f$.

\begin{enumerate}[label= (\thesection.B\alph*), leftmargin=3.2em]  
\item \label{eq:B1.1} \textbf{Continuity.} The mapping $(t,u,v) \mapsto
  B(t,u,v):[0,T]\times \WW \times W \rightarrow V^*$ is continuous on the sublevels of $\calE$, i.e., for every sequence
$(t_n,u_n,v_n)\to (t,u,v)$ in $[0,T]\ti \WW \ti W$ with $\sup_{t\in [0,T], n\in\mathbb{N}}\calE_t(u_n)<+\infty$, there holds $B(t_n,u_n,v_n)\to B(t,u,v)$ in $V^*$ as $n\rightarrow \infty$.
\item \label{eq:B1.2} \textbf{Control of the growth.} There exist positive constants $\beta>0$ and $c,\nu \in(0,1)$  such that
\begin{align*}
c\,\Psi^*\left(\frac{-B(t,u,v)}{c}\right)\leq \beta (1+\calE_t(u)+\vert v\vert ^2+ \Psi(v)^{\nu})
\end{align*} for all $u\in D\cap V,v\in V, \, t\in[0,T]$.
\item[(\thesection.f)]\label{eq:f1} \textbf{External force.} There holds $f\in \rmL^2(0,T;H)$.
\end{enumerate}
Similarly, we write (3.B) when we refer to the Conditions \ref{eq:B1.1}-\ref{eq:B1.1}. Several remarks are in order.

\begin{rem}\label{re:Assump.B1} 
\begin{itemize}
\item[$i)$]  In fact, the continuity of the perturbation with an image in $V^*$ is only needed to show the energy-dissipation inequality \eqref{sol:EDI}. If we only address the existence of solutions to the inclusion \eqref{sol:IC} without the energy-dissipation inequality, then it is sufficient to suppose that $B:[0,T]\times \WW \times W \rightarrow U^*+V^*$ is a mapping with values in $U^*+V^*$ which is continuous on sublevel sets of the energy.
\item[$ii)$] The condition \ref{eq:f1} can be relaxed to $f\in \rmL^1(0,T;H)+\rmL^2(0,T;V^*)$ in the Case \textbf{(a)} and to $f\in \rmL^1(0,T;H)+\rmL^2(0,T;V^*)+\rmL^{q^*}(0,T;\Wt^*)$ in the case \textbf{(b)}, where $q^*>1$ is again the conjugate exponent to $q>1$ from Assumption \ref{eq:Psi1.1}. 
  \end{itemize}
  \end{rem} 
  
\subsection{Discussion of the assumptions}\label{su:Diss.LDS}
Having collected the assumptions on the system $(U,V,W,H,\calE,\Psi,B,f)$ system, we want to discuss several conditions in more detail apart from the assertions and implications made in the remarks. In particular, we want to discuss the practical meaning of the assumptions and provide sufficient conditions for them to hold true. \\\\
As we already mentioned before, evolution equations of second order are, in general, more delicate than evolution equations of first order because of the nonsmoothing effect in time caused by the term $\partial_{tt} u$. This leads to a formation of discontinuities or a blow-up of a solution in finite time despite having smooth initial data, which makes it more difficult to prove strong solutions, see, e.g., \textsc{Zeidler} \cite[Section 33.7]{Zeid90NFA2b} for a discussion of these phenomena in connection with problems arising in physics. Therefore, we need well-adjusted assumptions.\\\\
Ad (\thesection.$\Uppsi$). Here, the dissipation potential is (in Case \textbf{(b)}) the sum of a leading part $\Psi_1$, which is defined by a strongly positive and bounded bilinear form, and a strongly continuous perturbation $\Psi_2$ of polynomial growth. As mentioned in Remark \ref{re:Assump.Psi1} $i)$,  the subdifferential $\partial \Psi= A+\partial \Psi_2$ is given by the sum of a positive, linear, bounded, and symmetric operator $A$ and the subdifferential of a strongly continuous perturbation $\Psi_2$. The essential assumption here is that $A$ is a positive, linear, bounded, and symmetric operator, which is crucial in identifying the limit of a sequence of approximate solutions stemming from a discretization scheme, see Section \ref{se:VarApprox.1}.

Admissible examples of dissipation potentials are, e.g., 
\begin{align*}
\Psi_1(v)=\frac{1}{2}\int_\Omega \vert \nabla^n v(\xx)\vert^2\dd \xx  \quad \text{and}\quad \Psi_2(v)=\frac{1}{p}\int_\Omega \left(\vert \nabla^{m} v(\xx)\vert^p+\vert \nabla^{k} v(\xx)\vert\right)\dd \xx 
\end{align*} on the \textsc{Sobolev} space $\rmH^n(\Omega)$ for any $k,m,n\in \mathbb{N}$ with $k\leq m<n$ and $p\in (1,+\infty)$ such that the compact embedding $\rmH^n(\Omega)\overset{c}{\hookrightarrow }\rmW^{m,p}(\Omega)$ holds.
\\\\
Ad (\thesection.E). In Condition \ref{eq:cond.E1.1}, we assume that the effective domain of $\calE_t$ is assumed to be time-independent.  The assumption that $D\cap V$ is dense in $D$ ensures the existence of an approximating sequence in $D\cap V$ to any initial value $u_0\in D$, which is needed in order to obtain a priori estimates for the interpolations in Lemma \ref{le:DUEE1}. The non-emptiness of $D$ in the case that $\calE_t$ is convex ensures the existence of a continuity point for $\calE_t$, see \textsc{Ekeland \& Temam} \cite[Corollary 2.5., p. 13]{EkeTem76CAVP}, which in turn allows us to use the variational sum rule in Lemma \ref{le:Subdif2}. \\ The Assumption \ref{eq:cond.E1.6} imposes the closedness of the graph of $\mathrm{Gr}(\partial \calE)$ on \textsc{Bochner} spaces which mimics the closedness result for maximal monotone operators \cite{Barb10NDMT}.

The Condition \ref{eq:cond.E1.8} is necessary to obtain appropriate a priori estimates for the subgradients of $\calE_t$, which in turn is needed to obtain a priori estimates for $u''$ since we solely obtain a priori estimates for the sum of the subgradient of $\calE_t$ and $u''$. Condition \ref{eq:cond.E1.8} could be replaced by the more general condition that $\partial\calE_t$ is a bounded operator.
\\\\ Ad (\thesection.B). The continuity condition \ref{eq:B1.1} means that $B$ is a continuous perturbation of $\partial \calE_t$. The term $B$  contains in the applications non-variational and non-monotone terms of lower order in terms of only containing lower order derivatives as well as obeying a growth condition of lower order. This is reflected by Condition \ref{eq:B1.2}, where $B$ satisfies a growth condition in terms of the conjugate dissipation potential $\Psi^*$ and the energy functional $\calE_t$. In fact, the growth condition shows that the higher the order of the growths of $\Psi$ and $\calE_t$ are, the more we can allow for the growth of the perturbation. Condition \ref{eq:B1.2} ensures that we are able to derive a priori estimates. Both conditions can be generalized in a framework that instead of a point-wise continuity and a pointwise growth condition, a continuity on suitable \textsc{Bochner} spaces can be imposed, as well as a growth condition on the level of time integrals. Furthermore, it would be sufficient to define the perturbation on the domain of the subdifferential of $\calE_t$. We note that the perturbation is nonlinear in both $u$ and $u'$.
\subsection{Statement of the main result}
Now, we are in the position to state the main result.
 
\begin{thm}[Main result]
  \label{th:MainExist1} Let the linearly damped inertial system \\
  $(U,V,W,H,\calE,\Psi,B,f)$ be given and fulfill Assumptions  \textnormal{(\thesection.E), (\thesection.$\Uppsi$)} and \textnormal{(\thesection.B)} as well as Assumption \textnormal{(\thesection.f)}. Let $\lambda\geq0$ be from Condition \textnormal{\ref{eq:cond.E1.7}}. Furthermore, let the initial values $u_0\in D\cap V_\lambda$ and $v_0\in H$ be given and assume that there exists a sequence $(u_0^k)_{k\in \mathbb{N}}, u_0^k \in D\cap V,$ such that
  \begin{align*}
  u_0^k\rightarrow u_0 \quad \text{in $V_\lambda$ as $ k\rightarrow\infty$ } \quad\text{ and}\quad \sup_{k\in \mathbb{N}}\calE_0(u_0^k)<+\infty.
  \end{align*} Then, there exists a strong solution to \eqref{eq:I.2}, i.e., there exist functions
  \begin{align}\label{sol:SP}
  u\in\rmC_{w}([0,T];U) \cap\rmW^{1,\infty}(0,T;H) \text{ with $u'\in \rmL^2(0,T;V)$,}\quad \xi\in \rmL^\infty(0,T;V_\lambda^*),
\end{align} additionally satisfying $u\in \rmH^{2}(0,T;U^*+V^*)$ in Case \textbf{(a)} and $u\in \rmW^{1,q}(0,T;\Wt)\cap \rmW^{2,\min\lbrace 2,q^*\rbrace}(0,T;U^*+V^*)$ in Case \textbf{(b)} such that the initial conditions $u(0)=u_0$ in $V_\lambda$ and $u'(0)=v_0$ in $H$ as well as the inclusions
  \begin{align}\label{sol:IC}
\xi(t)\in \partial_{V_\lambda} \calE_t(u(t)), \,\, f(t)\in u''(t)+\partial\Psi(u'(t))+\xi(t)+B(t,u(t),u'(t)) \quad \text{in }U^*+V^*
\end{align} are fulfilled for almost every $t \in(0,T)$. Furthermore, the energy-dissipation inequality 
\begin{align} \label{sol:EDI}
&\frac{1}{2}\vert u'(t) \vert^2 +\calE_t(u(t)) + \int_s^t \left(\Psi(u'(r))+\Psi^*(S(r)-\xi(r)-u''(r)) \right)\dd r \notag \\
&\leq \frac{1}{2}\vert u'(s) \vert^2 +\calE_s(u(s)) + \int_s^t  \partial_r \calE_r(u(r)) \dd r+\int_s^t \langle S(r),u'(r) \rangle_{V^*\times V}\dd r,
\end{align} holds for all $0 <t\leq T$ for $s=0$, and almost every $s\in (0,t)$, where $S(r):=f(r)-B(r,u(r),u'(r)),\, r\in[0,T]$, and $V_\lambda=U$ if $\calE_t$ is convex, i.e., $\lambda=0$, and $V_\lambda=U\cap V$ if $\calE_t$ is $\lambda$-convex with $\lambda>0$. 
\end{thm} 

\section{Proof of the main result}
In order to prove the main result, we proceed as follows:
\begin{itemize}
\item[\textbf{1)}] We regularize the dissipation potential $\Psi_2$ with regularization parameter $\varepsilon>0$.
\item[\textbf{2)}] We discretize the inclusion in time via a semi-implicit \textsc{Euler} scheme with time step $\tau>0$ and show solvability of the discrete problem.
\item[\textbf{3)}] We define interpolations functions and show a priori estimates for them.
\item[\textbf{4)}] We show the compactness of the interpolation functions in suitable spaces. 
\item[\textbf{5)}] We pass to the limit with $\tau\searrow 0$ and show existence of solutions to the regularized problem.
\item[\textbf{6)}] We pass to the limit with $\varepsilon\searrow 0$ and show existence of solutions to the original problem.
\end{itemize}

In the following, each step will be carried out in a subsection.

\subsection{Regularized problem} \label{se:regularized.problem}
In an intermediate step of our main proof, we consider the regularized problem
\begin{align} 
\label{eq:I.2.reg}
\begin{cases}
u''(t)+\partial \Psi_1(u'(t))+\rmD_G\Psi^\varepsilon_2(u'(t))+\partial \calE_t(u(t))+B(t,u(t),u'(t)) \ni f(t), &\quad t\in (0,T),\\
u(0)=u_0, \quad u'(0)=v_0,
\end{cases}
\end{align} where $\Psi_2^\varepsilon :W\rightarrow \mathbb{R}$ denotes the $q$-\textsc{Moreau--Yosida} regularization defined by
\begin{align*}
\Psi_2^\varepsilon(v)=\inf_{w\in W}\left \lbrace \frac{\varepsilon}{q}\left \Vert \frac{v-w}{\varepsilon}\right\Vert^p+\Psi_2(w)\right \rbrace,
\end{align*} where $q>1$ is from Condition $(3.\Uppsi)$. It is easily checked that $\Psi_2^\varepsilon$ satisfies the growth condition \eqref{eq:Psi.growth} for constants independent of $\varepsilon$ for sufficiently small $\varepsilon$ (and hence also for it's convex conjugate $\Psi_2^{\varepsilon,*}$). Furthermore, since $W$ is a reflexive Banach space such that $W^*$ is uniformly convex, by B. \cite[Theorem 2.2]{Bach22GMYR}, the functional $\Psi_2^\varepsilon$ \textsc{G\^{a}teaux} differentiable with continuous \textsc{G\^{a}teaux} derivative denoted by $\rmD_G \Psi^\varepsilon_2$. The \textsc{G\^{a}teaux} derivative satisfies in addition the bound
\begin{align*}
\Vert \rmD_G \Psi^\varepsilon_2(v)\Vert_{W^*}\leq C(\varepsilon)(1+\Vert v\Vert_W^{q-1})
\end{align*} for a constant $C=C(\varepsilon)$ dependent on $\varepsilon$, see, e.g., \cite{AkaMel18ERNP}. In fact, as $\varepsilon \searrow 0$ the constant $C(\varepsilon)$ tends to $+\infty$. Hence, the dissipation potential $\Psi=\Psi_1+\Psi^\varepsilon_2$ complies with the conditions $(3.\Uppsi)$.
Moreover, it has been shown in B. \cite[Theorem 2.3]{Bach22GMYR} that $\Psi_2^\varepsilon \Mto \Psi_2$ as $\varepsilon\searrow 0$ in the sense of \textsc{Mosco} convergence. In the following sections, we suppress the dependency of $\Psi_2^\varepsilon$ from $\varepsilon$ for readability.

\subsection[Variational approximation scheme]{Variational approximation scheme}\label{se:VarApprox.1}
As mentioned above, the proof of Theorem \ref{th:MainExist1} is based on the construction of strong solutions to \eqref{eq:I.2} via a
semi-implicit time discretization scheme. More specifically, we will employ a semi-implicit \textsc{Euler} method where all terms will be discretized implicitly, except for the non-variational perturbation term $B$ in order to obtain a variational approximation scheme to inclusion \eqref{eq:I.2}. Therefore, let for $N\in
\mathbb{N}\backslash \lbrace 0\rbrace$
\begin{align*}
 I_\tau=\lbrace 0=t_0<t_1<\cdots< t_n=n\tau<\cdots <t_N=T \rbrace
\end{align*} be an equidistant partition of the time interval $[0,T]$ with step size $\tau:=T/N$, where we omit the dependence of the nodes of the partition on the step size for simplicity. The discretization of \eqref{eq:I.2} is then given by
\begin{align}
\label{eq:EuLa}
\frac{\va-\vb}{\tau}+\partial_{V_\lambda}\Psi \left(\va\right)+\partial_{V_\lambda} \calE_{\ta}(\ua) +B\left(\ta,\ub,\vb\right) \ni f_\tau^n \quad \text{in }U^*+V^*
 \end{align} for $n=1,\dots,N$, where $V_\tau^0=v_0$, $\va:=\frac{\ua-\ub}{\tau}$, and $f_\tau^n :=\frac{1}{\tau}\int_{\tb}^{\ta} f(\sigma)\dd\sigma, n=1,\dots,N$, where $\partial_{V_\lambda}$ denotes the subdifferential operator with respect to the strong topology of $V_\lambda$. The inclusion \eqref{eq:EuLa} is equivalent to saying that there exists a subgradient $\xi_\tau^n \in \partial_{V_\lambda}\calE_t(\ua)$ such that 
 \begin{align*}
\frac{\va-\vb}{\tau}+\rmD_G \Psi(\va)+\xi_\tau^n +B\left(\ta,\ub,\vb\right) = f_\tau^n \quad \text{in }U^*+V^*,
 \end{align*} where $\rmD_G \Psi(\va)=A\va$ in Case \textbf{(a)} and   $\rmD_G \Psi(\va)=A\va+\rmD_G \Psi_2(\va)$ in Case \textbf{(b)}.
 
 The values $\ua \approx u(t_n)$ and $\va \approx u'(\ta)$ shall approximate the exact solution and its time derivative, and are to be determined successively from \eqref{eq:EuLa}. By Lemma \ref{le:Subdif2}, it follows that the approximate value $\ua$ is characterized as the solution to the \textsc{Euler--Lagrange} equation
associated with the mapping
\begin{align*}
u\mapsto \Upphi(\tau,\tb,\ub,\uc,B(\ta,\ub,\vb)-f_\tau^n;u),
\end{align*} where 
\begin{align*}
  \Upphi(r,t,v,w,\eta;u)= \frac{1}{2r^2} \vert u-2v-w\vert^2+
  r\Psi \left(\frac{u-v}{r}\right)+\calE_{t+r}(u)+\langle \eta,u\rangle_{V^*\times V}
\end{align*}
for $ r\in \mathbb{R}^{>0},t\in[0,T)$ with $r+t\in[0,T]$,  $u\in D,v\in V, w \in H$,
and $ \eta\in V^*$.\\
 We end up with the recursive scheme
\begin{align}\label{eq:ApproxSc}
\begin{cases}
\un \in D\cap V \text{ and }\vn \in V \text{ are given; whenever $\ue,\dots,\ub \in D\cap V$ are known,}\\
\text{find } \ua \in J_{\tau,\tb}(\ub,\uc;B(\ta,\ub,\vb)-f_\tau^n)
\end{cases}
\end{align} for $n=1,\dots, N$, where $J_{r,t}(v,w;\eta): =\mathrm{argmin}_{u\in U\cap V}\Upphi(r,t,v,w,\eta;u)$ and $U^{-1}_\tau=\un-\vn \tau$. 

The following lemma ensures the solvability of the variational scheme \eqref{eq:ApproxSc}. 

\begin{lem}\label{le:Exist.Min1} Let the linearly damped inertial system  $(U,V,W, \Wt,H,\calE,\Psi)$ be given and let the Conditions \textnormal{ \ref{eq:cond.E1.1}-\ref{eq:cond.E1.3}, \ref{eq:cond.E1.7}}, and \textnormal{\ref{eq:Psi1.1}} be fulfilled. Furthermore, let $r\in (0,T)$ and $t\in [0,T)$ with $r+t\leq T$ as well as $v\in V, w\in H$ and $\eta \in V^*$. Then, the set $J_{r,t}(v,w;\eta)$ is non-empty and single valued if $r\leq \frac{\mu}{4\lambda}$, where $\mu$ and $\lambda$ are from \textnormal{\ref{eq:Psi1.1} and \ref{eq:cond.E1.7}}, respectively. Furthermore, to every $u\in J_{r,t}(v,w;\eta)$ there exists $\xi\in \partial_{V_\lambda}\calE_t(u)\subset V_\lambda^*$ such that
  \begin{align*}
 \frac{u-2v-w}{r^2}+\rmD_G\Psi\left( \frac{u-v}{r}\right)+\xi+\eta=0 \quad \text{ in }U^*+V^*.
 \end{align*}
  \end{lem}
\begin{proof}[Proof] Since the proof is similar for the cases \textbf{(a)} and \textbf{(b)}, we restrict the proof by showing the assertion for the case \textbf{(b)}. Let $u\in D\cap V,v\in V,w\in H, \eta \in V^*$, and $r\in (0,T), t\in [0,T)$ with $r+t\leq T$ be given. Employing the \textsc{Fenchel--Young} inequality, we obtain 
\begin{align}
\label{eq:II.11}
   \Upphi(r,t,v,w,\eta;u)&= \frac{1}{2r^2} \vert u-2v+w\vert^2+
  r\Psi\left(\frac{u-v}{r}\right)+\calE_{t+r}(u)+\langle \eta,u\rangle_{V^*\times V}\notag \\
  &= \frac{1}{2r^2} \vert u-2v+w\vert^2+
   r\Psi_1\left(\frac{u-v}{r}\right)+ r\Psi_2\left(\frac{u-v}{r}\right)+\calE_{t+r}(u)\notag \\
   &\quad +\langle \eta,u\rangle_{V^*\times V} \notag\\
   &\geq  \frac{1}{2r^2} \vert u-2v+w\vert^2+
   \frac{1}{r}\Psi_1\left(u-v\right)+ r\hat{c}\left(\left\Vert\frac{u-v}{r}\right\Vert_\Wt^q-1\right)+\calE_{t+r}(u)\notag \\
   &\quad +\langle \eta,u-v\rangle_{V^*\times V}+\langle \eta,v\rangle_{V^*\times V} \notag \\
    &\geq   \frac{1}{2r^2} \vert u-2v+w\vert^2+
  \left(\frac{1}{r}-\varepsilon\right)\Psi_1\left( u-v \right)-r\hat{c}+\calE_{t+r}(u) \notag\\
  &\quad -\varepsilon\Psi_1^*\left(-\frac{\eta}{\varepsilon}\right)-\langle \eta,v\rangle_{V^*\times V}  \\
    &\geq   \frac{1}{2r^2} \vert u-2v+w\vert^2+
  \left(\frac{1}{r}-\varepsilon\right)\Psi_1\left( u-v \right)-\varepsilon\Psi_1^*\left(\frac{\eta}{\varepsilon}\right)-r\hat{c}\notag \\
  &\quad-\langle \eta,v\rangle_{V^*\times V}  \notag,
\end{align} 
for all $0<\varepsilon<\frac{1}{r}$, where in the two penultimate inequalities we used the inequality \eqref{eq:Psi.growth}. In regard of Conditions \ref{eq:Psi1.1},  \ref{eq:cond.E1.2}, and \ref{eq:cond.E1.3}, this implies $\inf_{u\in U\cap V} \Upphi(r,t,v,w,\eta;u)>-\infty$. On the other hand, we observe that
\begin{align}
\label{eq:II.12}
  \inf_{u\in U\cap V} \Upphi(r,t,v,w,\eta;u) \leq \frac{1}{2r^2} \vert \bar{u}-2v+w\vert^2+
  r\Psi\left(\frac{
 \bar{u}-v}{r}\right)+\calE_{t+r}(\bar{u})+\langle \eta,\bar{u}\rangle_{V^*\times V}
\end{align} for any $\bar{u} \in D\cap V$,
so that $\inf_{u\in U\cap V} \Upphi(r,t,v,w,\eta;u)<+\infty$ holds as well. It remains to show that the global minimum is attained by an element in $D\cap V$. In order to show that, let $(u_n)_{n \in
  \mathbb{N}}\subset U\cap V$ be a minimizing sequence for
$\Upphi(r,t,v,w,\eta;\cdot)$. From \eqref{eq:II.11} and the coercivity of $\Psi_1$ and $\calE$, we deduce that
$(u_n)_{n \in \mathbb{N}}\subset U\cap V$ is contained in a sublevel set of $\Psi_1$ and $\calE$, and thus bounded in $U\cap V$. Hence, by the reflexivity of $U\cap V$, there exists a subsequence (not relabeled) that converges weakly in $U\cap V$ to a
limit $\tilde{u}\in U\cap V$. We note that by Conditions \ref{eq:Psi1.1} and \ref{eq:cond.E1.1} as well as the continuity of the \textsc{Hilbert} space norm and the duality pairing, the sequential weak lower semicontinuity of the mapping
$u\mapsto \Upphi(r,t,v,w,\eta;u)$ on $U\cap V$ follows.  Then, by the sequential weak lower semicontinuity of $\Upphi$, we have 
\begin{align*}
  \Upphi(r,t,v,w,\eta;\tilde{u})\leq \liminf_{n\rightarrow
    \infty}\Upphi(r,t,v,w,\eta;u_n)=\inf_{\tilde{v}\in U\cap V}
  \Upphi(r,t,v,w,\eta;\tilde{v}),
\end{align*} and therefore, $u\in J_{r,t}(v,w;\eta)\neq \emptyset$ and $u\in D\cap V$. If $r>0$ is sufficiently small, then there is a unique global minimizer. Indeed, assuming there are two different global minimizer $\tilde{u}_0,\tilde{u}_1\in D\cap V$, then in view of the $\lambda$-convexity of $\calE_t$, the convexity of $\Psi_2$ and the fact that $\vert \cdot\vert^2$ and $\Psi_1$ fulfil a parallelogram identity, we obtain for every $s\in (0,1)$
\begin{align*}
 &\Upphi(r,t,v,w,\eta;s\tilde{u}_0+(1-s)\tilde{u}_1) \\
 &\leq s \Upphi(r,t,v,w,\eta;\tilde{u}_0)+(1-s)\Upphi(r,t,v,w,\eta; \tilde{u}_1)-s(1-s)\vert \tilde{u}_0-\tilde{u}_1\vert^2+\lambda\Vert \tilde{u}_0-\tilde{u}_1\Vert_V^2\\ 
 &\quad-\frac{s(1-s)}{r}\Psi_1(\tilde{u}_0-\tilde{u}_1)\\
 &= \min_{\tilde{v}\in U\cap V}\Upphi(r,t,v,w,\eta;\tilde{v})-s(1-s)\vert \tilde{u}_0-\tilde{u}_1\vert^2-\frac{s(1-s)}{r}\Psi_1(\tilde{u}_0-\tilde{u}_1)+\lambda\Vert \tilde{u}_0-\tilde{u}_1\Vert_V^2\\
&\leq \min_{\tilde{v}\in U\cap V}\Upphi(r,t,v,w,\eta;\tilde{v})-s(1-s)\vert \tilde{u}_0-\tilde{u}_1\vert^2-\left(\frac{s(1-s)}{r}-\frac{\lambda}{\mu}\right)\Psi_1(\tilde{u}_0-\tilde{u}_1),
\end{align*} where we also used the strong positivity of $\Psi_1$ with constant $\mu>0$. Choosing $s=\frac{1}{2}$ and $r\leq \frac{\mu}{4\lambda}$, the uniqueness follows. In order to prove the last assertion, we first assume that the energy functional is $\lambda$-convex with $\lambda>0$. Then, from \textsc{Fermat}'s theorem, we know that for any minimizer $u\in J_{r,t}(v,w;\eta)$, the functional $\Upphi(r,t,v,w,\eta;\cdot)$ is subdifferentiable in $u$ and there holds
\begin{align*}
0&\in \partial_{U\cap V} \Upphi(r,t,v,w,\eta;u)\\
&=\partial_{U\cap V}\left(\frac{1}{2r^2} \vert u-2v+w\vert^2+ r\Psi\left(\frac{u-v}{r}\right)+\calE_{t+r}(u)+\langle \eta,u\rangle_{V^*\times V}\right).
\end{align*} Since all terms except the energy functional are convex and \textsc{G\^{a}teaux} differentiable on the space $U\cap V$, we obtain with Lemma \ref{le:Subdif2} that $\calE_t$ is subdifferentiable in $u$ and there holds
\begin{align*}
\frac{u-2v-w}{r^2}+\rmD_G\Psi\left( \frac{u-v}{r}\right)+\eta \in \partial_{U\cap V}\calE_t(u).
\end{align*} Thus, we define $\xi:=\frac{u-2v-w}{r^2}+\rmD_G\Psi\left( \frac{u-v}{r}\right)+\eta \in U^*+V^*$. Now, we consider the case when $\lambda=0$. Then, we define the functionals $\widetilde{\Psi}:U\rightarrow (-\infty,+\infty]$ and $h:U\rightarrow (-\infty,+\infty]$ by
\begin{align*}
\widetilde{\Psi}(\bar{v})=
\begin{cases}
&\Psi(\bar{v}) \quad \text{if }\bar{v}\in V\cap U \\
& +\infty \quad \text{otherwise}.
\end{cases}\quad \text{ and } \quad h(\bar{v})=
\begin{cases}
&\langle \eta ,\bar{v}\rangle_{V^*\times V} \quad \text{if }\bar{v}\in V\cap U \\
& +\infty \quad \text{otherwise}.
\end{cases}
\end{align*} It can be shown that $\widetilde{\Psi}+h$ is proper, convex, and lower semicontinuous on $U$. The first two properties are readily seen. For the lower semicontinuity, we make use of the equivalent characterization of the lower semicontinuity, which states that all sublevel sets are closed in the strong topology of $U$. Thus, let $\alpha\in \mathbb{R}$ and $(u_n)_{n\in \mathbb{N}}\subset J_\alpha:=\lbrace \tilde{v}\in U : \widetilde{\Psi}(\tilde{v})+h(\tilde{v})\leq \alpha\rbrace$ such that $u_n\rightarrow u\in U$ as $n\rightarrow \infty$. We want to show that $u\in J_\alpha$. From the definition of $\widetilde{\Psi}$ and $h$, there holds $J_\alpha =\lbrace \tilde{v}\in U : \Psi(\tilde{v})+\langle \eta,\tilde{v}\rangle_{V^*\times V}\leq \alpha\rbrace$ and that $(u_n)_{n\in \mathbb{N}}$ is bounded in $V$ by the coercivity of $\Psi$ on $V$. Hence, there exists a weakly convergent subsequence (labeled as before) such that $u_n\rightharpoonup \tilde{u}\in V$ as $n\rightarrow \infty$. Therefore, $(u_n)_{n\in \mathbb{N}}$ is bounded in $U\cap V$ and from the reflexivity of $U\cap V$, we can extract a further weakly convergent subsequence (labeled as before) such that $ u_n\rightharpoonup \hat{u}\in U\cap V$ as $n\rightarrow \infty$. We obtain,
\begin{align*}
\langle f,\hat{u}\rangle_{(U^*+V^*)\times (U\cap V)}&=\langle f_1,\hat{u} \rangle_{U^*\times U}+\langle f_2,\hat{u}\rangle_{V^*\times V}\\
&= \lim_{n\rightarrow \infty}\left( \langle f_1,u_n \rangle_{U^*\times U}+\langle f_2,u_n \rangle_{V^*\times V}\right)\\
&= \langle f_1,u\rangle_{U^*\times U}+\langle f_2, \tilde{u}\rangle_{V^*\times V} 
\end{align*}  for all $f=f_1+f_2\in U^*+V^*$ and in particular for all $f \in U^*$ and $f\in V^*$ whence $u=\bar{u}=\hat{u}$ in $U\cap V$. From the weak lower semicontinuity of $\Psi_1$ on $V$, we obtain
\begin{align*}
\Psi(u)+\langle \eta,u\rangle_{V^*\times V}\leq \liminf_{n\rightarrow \infty}\left(\Psi(u_n)+\langle \eta,u\rangle_{V^*\times V}\right) \leq \alpha,
\end{align*} and thus, $u\in J_\alpha$, from which the lower semicontinuity on $U$ follows. Noting that 
\begin{align*}
\min_{\tilde{v}\in U\cap V } \Upphi(r,t,v,w,\eta;\tilde{v})&=\min_{\tilde{v}\in U\cap V } \widetilde{\Upphi}(r,t,v,w,\eta;\tilde{v})\\
&=\min_{\tilde{v}\in U}\left( \frac{1}{2r^2} \vert \tilde{v}-2v+w\vert^2+ r\widetilde{\Psi}\left(\frac{\tilde{v}-v}{r}\right)+\calE_{t+r}(u)+h(\tilde{v})\right),\\
\end{align*} we obtain again by \textsc{Fermat}'s theorem
\begin{align*}
0&\in \partial_U \widetilde{\Upphi}(r,t,v,w,\eta;u)\\
&=\partial_{U}\left(\frac{1}{2r^2} \vert u-2v+w\vert^2+ r\widetilde{\Psi}\left(\frac{u-v}{r}\right)+\calE_{t+r}(u)+h(u)\right)
\end{align*} for any global minimizer $u\in J_{r,t}(v,w;\eta)$. To decompose the elements of the subdifferential of the sum of the functionals in terms of the subgradients of each functional, we employ Lemma \ref{le:Subdif2} and note that with Remark \ref{re:Assump.E1} $i)$ all assumptions of that lemma are satisfied. Hence, there holds
\begin{align*}
&\partial_{U}\left(\frac{1}{2r^2} \vert u-2v+w\vert^2+ r\widetilde{\Psi}\left(\frac{u-v}{r}\right)+\calE_{t+r}(u)+h(u)\right)\\
&=\partial_{U}\left(\frac{1}{2r^2} \vert u-2v+w\vert^2\right)+ \partial_{U}\left(r\widetilde{\Psi}\left(\frac{u-v}{r}\right)+h(u)\right)+\partial_{U}\calE_{t+r}(u)
\end{align*} and therefore there exists a subgradient $\xi\in \partial_U\calE_{t+r}(u)$ such that
\begin{align*}
-\frac{u-2v-w}{r^2}-\xi \in \partial_U \left( r\widetilde{\Psi}\left(\frac{u-v}{r}\right)+h(u)\right).
\end{align*} Unfortunately, we are not allowed to decompose the right subdifferential further since the functional $h$ is, in general, not lower semicontinuous on $U$. However, since the sum is proper, convex, and lower semicontinuous on $U$, we can make use of the equivalent description of the subdifferential by the inequality
\begin{align*}
r\widetilde{\Psi}\left(\frac{u-v}{r}\right)+h(u)-r\widetilde{\Psi}\left(\frac{\tilde{v}-v}{r}\right)+h(\tilde{v})\leq \left \langle -\frac{u-2v-w}{r^2}-\xi,\tilde{v}-u \right \rangle_{U^*\times U}
\end{align*} for all $\tilde{v}\in U$ and in particular
\begin{align*}
&r\Psi\left(\frac{u-v}{r}\right)+\langle \eta,u \rangle_{V^*\times V}-r\Psi \left(\frac{\tilde{v}-v}{r}\right)-\langle \eta,\tilde{v}\rangle_{V^*\times V} \\
&\leq \left \langle -\frac{u-2v-w}{r^2}-\xi,\tilde{v}-u \right \rangle_{(U^*+V^*)\times (U\cap V)}
\end{align*} for all $\tilde{v}\in U\cap V$, which in turn implies 
\begin{align*}
-\frac{u-2v-w}{r^2}-\xi &\in \partial_{U\cap V} \left(\Psi\left(\frac{u-v}{r}\right)+\langle \eta,u \rangle_{V^*\times V}\right)\\
&=\rmD_G \Psi\left(\frac{u-v}{r}\right)+\eta\in U^*+V^*,
\end{align*} which means that we can decompose the elements of the subdifferential in the weaker space $U^*+V^*$. We finally obtain
\begin{align*}
-\frac{u-2v-w}{r^2} +\rmD_G \Psi\left(\frac{u-v}{r}\right)+\xi+\eta \quad \text{in } U^*+V^*,
\end{align*} and hence the completion of the proof.

\end{proof}
\subsection{Discrete Energy-Dissipation inequality and a priori estimates}\label{se:TimeDiscret.1}
Since the previous lemma ensures the solvability of the approximation scheme \eqref{eq:ApproxSc}, we are now able to define piecewise linear and constant interpolations which will interpolate the values $(\ua)_{n=0}^N$ and $(\va)_{n=0}^N$ for every $\tau>0$, respectively, and we will derive a priori estimates for them. The interpolations shall approximate the desired solution to \eqref{eq:I.2} and its derivative, and are therefore also referred to as approximate solutions to \eqref{eq:I.2}. In order to define the approximate solutions, we assume for the moment that $u_0\in D\cap V$ and $v_0 \in V$. In the main existence proof, we will then approximate the initial values from $D\cap V_\lambda$ and $H$ by sequences from $V\cap V$ and $V$, respectively. For $\tau>0$, let $(U_\tau^n)_{n=1}^N \subset D\cap V$ be the sequence of
approximate values obtained from the variational approximation scheme \eqref{eq:ApproxSc} for $\un:=u_0$ and $\vn:= v_0$. Moreover, let $(\xi_\tau^n)_{n=1}^N \subset V_\lambda^*$ be a sequence of subgradients of the energy determined by the preceding lemma and satisfying $\xi_\tau^n \in \partial_{V_\lambda} \calE_{t_n}(U_\tau^n),\, i=1,\dots,N$ and \eqref{eq:EuLa}.
The piecewise constant and linear interpolations are defined by
\begin{align}
\label{eq:Approx.U}
  &\overline{U}_\tau(0)=\underline{U}_\tau(0)=\widehat{U}_\tau(0):=U_\tau^0=u_0\, \text{ and } \notag \\
  &\underline{U}_\tau(t):=U_\tau^{n-1}, \quad \widehat{U}_\tau(t): =\frac{t_n-t}{\tau}U_\tau^{n-1}+\frac{t-t_{n-1}}{\tau}U_\tau^{n} \quad \text{for } t\in[t_{n-1},t_n),\\
  &\overline{U}_\tau(t):=U_\tau^n \quad \text{ for } t\in(t_{n-1},t_n]
  \quad \text{and}\quad \underline{U}_\tau(T)=U_\tau^N,\, n=1,\dots,N,\notag
\end{align} 
as well as
\begin{align}
\label{eq:Approx.V}
  &\overline{V}_\tau(0)=\underline{V}_\tau(0)=\widehat{V}_\tau(0):=V_\tau^0=v_0 \,\text{ and } \notag \\
  &\underline{V}_\tau(t):=V_\tau^{n-1}, \quad \widehat{V}_\tau(t): =\frac{t_n-t}{\tau}V_\tau^{n-1}+\frac{t-t_{n-1}}{\tau}V_\tau^{n} \quad \text{for } t\in[t_{n-1},t_n),\\
  &\overline{V}_\tau(t):=V_\tau^n \quad \text{ for } t\in(t_{n-1},t_n]
  \quad \text{and} \quad \underline{V}_\tau(T)=V_\tau^N,\, \ n=1,\dots,N, \notag
\end{align}  where $\va=\frac{\ua-\ub}{\tau}$ for $n=1,\dots,N$. We note that $\widehat{U}_\tau'=\overline{V}_\tau$ in the weak sense.
Furthermore, we define the functions $\xi_\tau:[0,T]\rightarrow V_\lambda^*$ and $f_\tau:[0,T]\rightarrow H$ by

\begin{align} \label{eq:Approx.xif}
&\xi_\tau(t)= \xi^n_{\tau}, \quad f_\tau (t)=f_\tau^n=\frac{1}{\tau}\int_{t_{n-1}}^{t_n}f(\sigma)\dd \sigma \quad \quad \text{ for } t\in(t_{n-1},t_n],\, n=1,\dots,N, \\
& \xi_\tau(T)= \xi^N_{\tau} \quad \text{and} \quad f_\tau (T)=f_\tau^N.\notag
\end{align}

For notational convenience, we also introduce the piecewise constant functions $\overline{\mathbf{t}}_\tau : [0,T]
\rightarrow [0,T]$ and $\teu_\tau:[0,T]\rightarrow [0,T]$ given by
\begin{align}\label{eq:Approx.t}
\begin{split}
&\teo_{{\tau}}(0): = 0\,\,\,\, \text{ and } \, \teo_{{\tau}}(t): =
t_n \quad \text{ for } t\in (t_{n-1},t_n], \\
&\teu_{{\tau}}(T): = T \, \text{ and } \, \teu_{{\tau}}(t): =
t_{n} \quad \text{ for } t\in [t_{n-1},t_n), \quad n=1,\dots,N.
\end{split}
\end{align} Obviously, there holds $\teo_\tau(t)\rightarrow t$ and
$\teu_\tau(t)\rightarrow t$ as $\tau\rightarrow 0$. 

Finally, we are in a position to show useful priori estimates.

\begin{lem}[A priori estimates]\label{le:DUEE1} Let the system LDS $(U,V,W,\Wt,H,\calE,\Psi,B,f)$ be given and satisfy the Assumptions
  \textnormal{(3.E), (3.$\mathrm{\Psi}$)}, 
  \textnormal{(3.B)} as well as Assumption \textnormal{(3.f)}.  Furthermore, let $\overline{U}_\tau,\underline{U}_\tau, \widehat{U}_\tau,\overline{V}_\tau,\underline{V}_\tau, \widehat{V}_\tau, \xi_\tau$ and $f_\tau$ be the interpolations defined in
  \eqref{eq:Approx.U}-\eqref{eq:Approx.xif} associated with the given initial values $u_0\in D\cap V, v_0 \in V$ and the
  step size $\tau>0$. Then, the discrete energy-dissipation inequality
\begin{align}
\label{eq: DUEE1}
 &\frac{1}{2} \left \vert \Vo(t)\right \vert^2+ \calE_{\teo_\tau(t)}(\overline{U}_\tau(t)) + \int_{\teo_\tau(s)}^{\teo_\tau(t)}\left(
    \Psi ( \Vo (r) ) +
    \Psi^*\left(S_\tau(r)-\Vh'(r)-
      \xi_\tau(r)\right) \right) \dd r
   \notag \\ 
  &\leq \frac{1}{2} \left \vert \Vo(s) \right \vert^2+
  \calE_{\teo_\tau(s)}(\overline{U}_\tau(s))+\int_{\teo_\tau(s)}^{\teo_\tau(t)} \partial_r
  \calE_r(\Uu (r))\dd r+\int_{\teo_\tau(s)}^{\teo_\tau(t)}
  \langle S_\tau(r), \Vo(r) \rangle_{V^*\times V}\dd r \notag \\
 &\quad +\tau \lambda\int_{\teo_\tau(s)}^{\teo_\tau(t)}\Vert \Vo(r) \Vert_V^2\dd r 
\end{align} holds for all $0\leq s< t\leq T$, where $S_\tau(r): = f_\tau(r)- B(\teo_\tau(t),\Uu (t),\Vu(t)), r\in[0,T]$. Moreover, there exist positive constants $M,\tau^*>0$ such that the estimates 
\begin{align}
\label{eq:II.29.2}
\sup_{t\in [0,T]}  \left \vert \Vo(t)\right \vert \leq M, \quad \sup_{t\in [0,T]} \calE_t(\Uo(t)) \leq M, \quad \sup_{t\in [0,T]}\vert \partial_t \calE_t(\Uu(t))\vert \leq M, \\ 
\label{eq:II.30.2}
\int_0^T \left( \Psi\left( \Vo(r)\right) + \Psi^*\left( S_\tau(r)-\Vh'(r)- \xi_\tau(r)\right) \right) \dd r\leq M
\end{align} hold for all $0<\tau\leq \tau^*$. In particular, the families of functions
\begin{subequations}
\label{eq:allbounds}
\begin{align}
\label{eq:bUo}
&(\Uo)_{0<\tau\leq \tau^*} \subset \rmL^\infty(0,T;U),\\
\label{eq:bxi}
&(\xi_\tau)_{0<\tau\leq \tau^*} \subset \rmL^\infty(0,T;V_\lambda^*),\\
\text{in Case \textbf{(a)}}\notag \\
\label{eq:bVo}
&(\Vo)_{0<\tau\leq \tau^*} \subset \rmL^2(0,T;V)\cap \rmL^\infty(0,T;H),\\
\label{eq:bVh'}
&(\Vh')_{0<\tau\leq \tau^*}\subset \rmL^{2}(0,T;U^*+V^*),\\
\label{eq:bB}
&(B_\tau)_{0<\tau\leq \tau^*}\subset \rmL^\frac{2}{\nu}(0,T;V^*),\\
\text{in Case \textbf{(b)}}\notag \\
\label{eq:bVob}
&(\Vo)_{0<\tau\leq \tau^*}\subset \rmL^{2}(0,T;V)\cap \rmL^q(0,T;\Wt)\cap \rmL^\infty(0,T;H),\\
\label{eq:bVh'b}
&(\Vh')_{0<\tau\leq \tau^*}\subset \rmL^{\min \lbrace 2,q\rbrace}(0,T;U^*+V^*),\\
\label{eq:bBb}
&(B_\tau)_{0<\tau\leq \tau^*}\subset \rmL^\frac{2}{\nu}(0,T;V^*)+\rmL^{\frac{q^*}{\nu}}(0,T;W^*),
\end{align} 
\end{subequations}
are uniformly bounded with respect to $\tau$ in the respective spaces, where we abbreviate $B_\tau (t):=B(\teo_\tau(t),\Uu (t),\Vu(t)),\ t\in[0,T]$, $q^*>0$ is the conjugate exponent of $q>1$, and $\nu\in (0,1)$ stemming from Assumption \textnormal{\ref{eq:B1.2}}. Finally, there holds
\begin{align}
\begin{split}
\label{eq:II.31.2}
\sup_{t\in [0,T]}\left( \Vert \underline{U}_\tau(t)-\overline{U}_\tau(t)\Vert_V+\Vert \widehat{U}_\tau(t)-\overline{U}_\tau(t)\Vert_V \right)\rightarrow 0 \\
\sup_{t\in [0,T]}\left(\Vert \Vu(t)-\Vo(t)\Vert_{U^*+V^*}+\Vert \Vh(t)-\Vo(t)\Vert_{U^*+V^*}\right)
\rightarrow 0
\end{split}
\end{align} as $\tau \rightarrow 0$.
\end{lem}
\begin{proof}[Proof]
Let $(U_\tau^n)_{n=1}^N \subset D\cap V$ be the approximative values obtained from the variational approximation scheme \eqref{eq:ApproxSc} and let $(\xi_\tau^n)_{n=1}^N \subset U^*+V^*$ be the associated subgradients. Then, by Lemma \ref{le:Subdif2}, the approximate value $\ua$ solves the \textsc{Euler--Lagrange} equation \eqref{eq:EuLa}, i.e., 
\begin{align}\label{eq:dis.inc}
S_\tau^n -\frac{\va-\vb}{\tau}-\xi^n_\tau \in \partial_{V\cap U} \Psi(\va)=\lbrace D_G\Psi(\va)\rbrace \quad \text{ and }\quad \xi^n_\tau\in \partial_{V_\lambda} \calE_{\ta}(\ua),
\end{align} where $S_\tau^n:=f_\tau^n -B(\ta,\ub,\vb)$. Due to Lemma \ref{le:Leg.Fen}, the first inclusion is equivalent to
\begin{align*}
\Psi (\va)+\Psi^*\left(S_\tau^n -\frac{\va-\vb}{\tau}-\xi_\tau^n\right)=\left \langle S_\tau^n -\frac{\va-\vb}{\tau}-\xi_\tau^n,\va \right \rangle_{V^*\times V}
\end{align*} and the second one implies 
\begin{align*}
-\left \langle \xi_\tau^n,\ua-\ub \right \rangle_{V_\lambda^*\times V_\lambda}&\leq \calE_{\ta}(\ub)-\calE_{\ta}(\ua)+\lambda \Vert \ua-\ub \Vert_V^2 \notag \\
&=\calE_{\tb}(\ub)-\calE_{\ta}(\ua)+\int_{\tb}^{\ta} \partial_r E_r(\ub) \dd r \\
&\quad+\lambda \Vert \ua-\ub \Vert_V^2
\end{align*} for all $n=1,\dots,N$.
Using the identity 
\begin{align} \label{rule}
\left( u-v,u \right)= \frac{1}{2} \left( \vert u \vert^2-\vert v \vert^2+ \vert u-v \vert^2\right) \quad \text{for all }u,v\in H
\end{align} and the fact that $\langle w,v\rangle_{V^*\times V} =(w,v)$ for $v\in V$ and $w\in H$, we obtain
\begin{align}\label{eq:pre.DUEE}
&\frac{1}{2} \vert \va \vert^2 +\calE_{\ta}(\ua)+ \tau \Psi(\va)+\tau\Psi^*\left(S_\tau^n -\frac{\va-\vb}{\tau}-\xi_\tau^n\right)-\tau\left \langle S_\tau^n,\va \right \rangle_{V^*\times V}\\
 &\leq \frac{1}{2} \vert \vb \vert^2 +\calE_{\tb}(\ub)+\int_{\tb}^{\ta} \partial_r \calE_r(\ub) \dd r +\lambda\Vert \ua-\ub \Vert_V^2 \notag
\end{align} for all $n=1,\dots,N$, which, by summing up the inequalities, implies \eqref{eq: DUEE1}.
In order to show the bounds \eqref{eq:II.29.2} and \eqref{eq:II.30.2}, we make use of the following estimates: first, from Assumption \ref{eq:B1.2} and the \textsc{Fenchel--Young} inequality, we obtain
\begin{align*}
\tau \langle S_\tau^n,\va \rangle&= \tau \langle -B(\ta,\ub,\vb)+f_\tau^n ,\va \rangle_{V^*\times V}\notag \\
&= \tau \langle -B(\ta,\ub,\vb),\va \rangle_{V^*\times V}+\tau \langle f_\tau^n,\va \rangle_{V^*\times V}\notag \\
&\leq  c \tau \Psi(\va)+c\tau\Psi^*\left(\frac{-B(\ta,\ub,\vb)}{c}\right)+ \frac{\tau}{2} (\vert f_\tau^n \vert^2+\vert \va \vert^2)\notag \\
&\leq c \tau \Psi(\va) + \tau \beta (1+\calE_{\ta}(\ub)+\vert \vb \vert^2 +\Psi(\vb)^\nu)\\
&+\frac{\tau}{2} (\vert f_\tau^n \vert^2+\vert \va \vert^2),\notag\\
&\leq c \tau \Psi(\va) + \tau \beta (1+\calE_{\ta}(\ub)+\vert \vb \vert^2)+\tau \varepsilon \Psi(\vb)+\tau C\\
&+\frac{\tau}{2} (\vert f_\tau^n \vert^2+\vert \va \vert^2),\notag
\end{align*} for positive constants $\varepsilon,C=C(\varepsilon,\beta)>0$ such that $\varepsilon<\frac{1-c}{2}$ and $C=\frac{\beta^{\frac{1}{1-\nu}}}{\varepsilon^\frac{\nu}{1-\nu}}$. Second, by the  strong positivity of $\Psi_1$ and the growth condition for $\Psi_2$, we have in Case \textbf{(b)}
\begin{align*}
\mu\Vert \ua-\ub\Vert^2_V=\mu \tau^2 \Vert \va\Vert_V^2\leq \tau^2 \Psi_1\left(\va\right)+\tau^2 \Psi_2\left(\va\right)+\tau^2\tilde{c}=\tau^2 \Psi\left(\va\right)+\tau^2\tilde{c}
\end{align*} where $\tilde{c}>0$ is from Condition \ref{eq:Psi.growth}. In the Case \textbf{(a)}, we only employ the strong positivity of $\Psi$ obtaining
\begin{align*}
\mu\Vert \ua-\ub\Vert^2_V=\mu \tau^2 \Vert \va\Vert_V^2\leq \tau^2 \Psi\left(\va\right).
\end{align*}

 Finally, we use the estimate following from Condition \ref{eq:cond.E1.4},
\begin{align*}
\int_{\tb}^{\ta} \partial_r \calE_r(\ub) \dd r\leq \int_{\tb}^{\ta} C_1 \calE_r(\ub) \dd r \leq C_1 \int_{\tb}^{\ta}  \sup_{t\in [0,T]}\calE_t(\ub) \dd r.
\end{align*} 
Inserting all preceding inequalities in \eqref{eq:pre.DUEE} and summing up all inequalities from $1$ to $n$, we find a positive constant $C>0$ such that
\begin{align}
&\frac{1}{2} \vert \va \vert^2 +\frac{1}{C_1}\sup_{t\in [0,T]} \calE_t(\ua)+ \int_{0}^{\ta}\left((1- \alpha(\tau)) \Psi(\Vo(r))+\Psi^*\left(S_\tau (r) -\Vh'(r)-\xi_{\tau}(r) \right) \right) \dd r \notag \\
 &\leq C\left(\vert v_0 \vert^2 +\calE_{0}(u_0)+T+\Vert f \Vert^2_{\rmL^2(0,T;H)}+\Psi(v_0)\right) \label{eq:vo} \\
  &+C\int_{0}^{\ta} \left( \vert \Vo(r)\vert^2+ \sup_{t\in [0,T]}\calE_t(\Uo(r))\right)\dd r,\notag
\end{align}  where $\alpha(\tau):=c+\ct+\tau\frac{\lambda}{\mu}<1$ for all $\tau<\tau^*: =\min\lbrace\frac{\mu}{\lambda}(1-c-\tilde{c}),1\rbrace$ and $\alpha(\tau)$ is decreasing for decreasing $\tau$. In the step \eqref{eq:vo}, we made use of the estimate for the interpolation $f_\tau$
\begin{align}\label{eq:f.est}
\Vert f_\tau\Vert^2_{\rmL^2(0,T;H)}&=\sum_{k=1}^n \tau \vert f_\tau^k\vert^2\notag\\
&=\sum_{k=1}^n \frac{1}{\tau}\vert \int_{t_{k-1}}^{t_k}f(\sigma)\dd \sigma \vert^2\notag \\
&\leq  \sum_{k=1}^n \int_{t_{k-1}}^{t_k} \vert f(\sigma)\vert^2\dd \sigma= \int_0^{t_n}\vert f(\sigma)\vert^2\dd \sigma\leq \Vert f\Vert^2_{\rmL^2(0,T;H)}.
\end{align} Then, by the discrete version of \textsc{Gronwall}'s lemma, see, e.g., \cite[Lemma 3.2.4, p. 68]{AmGiSa05GFMS}, there exists a constant $M>0$ such that \eqref{eq:II.29.2} and \eqref{eq:II.30.2} are satisfied. Now, we seek to show the bounds in \eqref{eq:allbounds} by distinguishing the two cases \textbf{(a)} and \textbf{(b)}.\\
Ad \textbf{(a)}. Due to the coercivity of $\Psi$ and $\Psi^*$, the uniform boundedness of $(\Vo)_{0<\tau\leq \tau^*}\subset \rmL^2(0,T;V)$ and $(S_\tau-\Vh'- \xi_\tau)_{0<\tau\leq \tau^*}\subset \rmL^2(0,T;V^*)\subset \rmL^2(0,T;V^*)$ follow immediately from the a priori estimate \eqref{eq:II.30.2}. The boundedness of $(B_\tau)_{0<\tau\leq \tau^*}\subset \rmL^2(0,T;V^*)$ uniformly in $\tau$ is a consequence of Assumption \ref{eq:B1.2} and the coercivity of $\Psi^*$: there holds
\begin{align}\label{eq:II.38.1}
\bar{c}\int_0^T \Vert B_\tau(r))\Vert_*^\frac{2}{\nu} \dd r&\leq \int_0^T \Psi^*\left( B(\teo_{\tau}(r),\Uu(r),\Vu(r))\right)^\frac{1}{\nu}\dd r \notag\\
&\leq \int_0^T c\Psi^*\left(\frac{ B(\teo_{\tau}(r),\Uu(r),\Vu(r))}{c}\right)^\frac{1}{\nu}\dd r \notag\\
&\leq \int_0^T \left( C((1+\calE_{\teu_{\tau}(r)}(\Uu(r))^\frac{1}{\nu}+\vert \Vu(r)\vert^2)^\frac{1}{\nu}+ \Psi\left( \Vu(r) \right)\right)\dd r \notag\\
&\leq N 
\end{align} for positive constants $C,N>0$ independent of $\tau$, where $c\in (0,1)$ is from Assumption \ref{eq:B1.2} and where we have used the fact that for all $\zeta\in V^*$ the mapping $r\mapsto r\Psi^*(\zeta/r)$ is monotonically decreasing on $(0,+\infty)$ which follows from the convexity of $\Psi^*$ and $\Psi^*(0)=0$. Since $(f_\tau)_{0<\tau\leq \tau^*}$ is uniformly bounded in $\rmL^2(0,T;H)$, we infer that $(\Vh'+\xi_\tau)_{0<\tau\leq \tau^*}$ is uniformly bounded in $\rmL^2(0,T;V^*)$ with respect to $\tau$ as well. Finally, Assumption \ref{eq:cond.E1.8} implies a uniform bound for $(\xi_\tau)_{0<\tau\leq \tau^*}$ in $\rmL^\infty(0,T;V_\lambda^*)$. Since all previous families of functions are bound in the common space $\rmL^2(0,T;U^*+V^*)$, we deduce that $(\Vh')_{0<\tau\leq \tau^*}$ is uniformly bounded in $\rmL^\infty(0,T;U^*+V^*)$ with respect to $\tau$.\\
Ad \textbf{(b)}. Again, the coercivity of the dissipation potential $\Psi$ leads to the boundedness of the sequence of discrete derivatives $(\Vo)$ in $\rmL^{2}(0,T;V)\cap \rmL^q(0,T;\Wt)$ uniformly in $\tau\in (0,\tau^*$). In order to show that $(S_\tau-\Vh'- \xi_\tau)_{0<\tau\leq \tau^*}\subset \rmL^2(0,T;V^*)+\rmL^{q^*}(0,T;\Wt^*)$ is uniformly bounded with respect to $\tau$, we make the following observation: let $\zeta:[0,T]\rightarrow V^*$ be any measurable function such that 
\begin{align*}
\int_0^T \Psi^*(\zeta(t))\dd t\leq M.
\end{align*} We want to show that there exists a positive constant $\widetilde{M}>0$ such that
\begin{align*}
\tilde{M}\geq \Vert \zeta \Vert_{\rmL^2(0,T;V^*)+\rmL^{q^*}(0,T;W^*)}.
\end{align*} First, by the formula \eqref{eq:brez.formula} from Lemma \ref{le:Bre.Att} and the growth conditions for the conjugate in Remark \ref{re:Assump.Psi1}, there exists a constant $M_1>0$ such that
\begin{align*}
&M_1\geq \int_0^T \min_{\eta\in W^*}\left( \Vert\zeta(t)-\eta\Vert^2_{V^*}+\Vert \eta\Vert^{q^*}_{W^*}\right) \dd t.
\end{align*} Second, the mapping $t\mapsto \alpha(t): =\min_{\eta\in W^*}\left( \Vert\zeta(t)-\eta\Vert^2_{V^*}+\Vert \eta\Vert^{q^*}_{W^*}\right) $ is \textsc{Lebesgue} measurable: since $V^*$ is separable, there exists a countable dense subset $(\eta_n)_{n\in \mathbb{N}}\subset V^*$. Then, there holds $\alpha(t)=\inf_{n\in \mathbb{N}}\left( \Vert\zeta(t)-\eta_n\Vert^2_{V^*}+\Vert \eta_n\Vert^{q^*}_{W^*}\right) $ and from the measurability of the function $\alpha_n(t): = \left( \Vert\zeta(t)-\eta_n\Vert^2_{V^*}+\Vert \eta_n\Vert^{q^*}_{W^*}\right)$ for each $n\in \mathbb{N}$, the measurability of $\alpha$ follows. Further, we note that the mapping $g:[0,T]\times W^*\rightarrow \mathbb{R}, \, (t,\eta)\mapsto g(t,\eta)=\Vert\zeta(t)-\eta \Vert^2_{V^*}+\Vert \eta \Vert^{q^*}_{W^*}$ is a \textsc{Carath\'{e}odory} function and therefore, by the Inverse Image Theorem, see, e.g., \textsc{Aubin \& Frankowska} \cite[Theorem 8.2.9, p. 315]{AubFra90SVA}, the set-valued map 
\begin{align*}
H(t): =\lbrace \eta\in W^* : g(t,\eta)=\alpha(t)\rbrace
\end{align*} is measurable and there exists a measurable selection $\omega:[0,T]\rightarrow W^*$ with $\omega(t)\in H(t)$ and $g(t,\omega(t))=\alpha(t)$ for all $t\in [0,T]$. We obtain
\begin{align*}
M_1&\geq \int_0^T \min_{\eta\in W^*}\left( \Vert\zeta(t)-\eta\Vert^2_*+\Vert \eta\Vert^{q^*}_{W^*}\right) \dd t\\
&=\int_0^T \left( \Vert\zeta(t)-\omega(t)\Vert^2_*+\Vert \omega(t)\Vert^{q^*}_{W^*}\right) \dd t,
\end{align*} whence $\omega \in \rmL^{q^*}(0,T;W^*)$ and $\zeta-\omega \in \rmL^2(0,T;V^*)$. It follows
\begin{align*}
M_1&\geq \int_0^T \left( \Vert\zeta(t)-\omega(t)\Vert^2_*+\Vert \omega(t)\Vert^{q^*}_{W^*}\right) \dd t,\\
&\geq  \inf_{\overset{\xi_1\in \rmL^2(0,T;V^*),\xi_2\in \rmL^{q^*}(0,T;W^*)}{\zeta=\xi_1+\xi_2}} \int_0^T \left( \Vert\xi_1(t)\Vert^2_*+\Vert \xi_2(t)\Vert^{q^*}_{W^*}\right) \dd t\\ 
 &\geq \inf_{\overset{\xi_1\in \rmL^2(0,T;V^*),\xi_2\in \rmL^{q^*}(0,T;W^*)}{\zeta=\xi_1+\xi_2}} \left(\Vert \xi_1\Vert_{\rmL^2(0,T;V^*)}+\Vert\xi_2\Vert_{\rmL^{q^*}(0,T;W^*)}\right)-M_2\\
 &\geq \inf_{\overset{\xi_1\in \rmL^2(0,T;V^*),\xi_2\in \rmL^{q^*}(0,T;W^*)}{\zeta=\xi_1+\xi_2}}\max \lbrace \Vert \xi_1\Vert_{\rmL^2(0,T;V^*)},\Vert \xi_2\Vert_{\rmL^{q^*}(0,T;W^*)}\rbrace-M_2\\
&= \Vert \zeta\Vert_{\rmL^2(0,T;V^*)+\rmL^{q^*}(0,T;W^*)}-M_2
\end{align*} for a constant $M_2>0$ coming from \textsc{Young}'s inequality. Since the constant $\widetilde{M}:=M_1+M_2>0$ was obtained independently of the function $\zeta$, the uniform bound of the sequence $(S_\tau-\Vh'- \xi_\tau)_{0<\tau\leq \tau^*}$ in $\rmL^2(0,T;V^*)+\rmL^{q^*}(0,T;W^*)$ follows. Employing Condition \ref{eq:B1.2} for the perturbation $B$ as in \eqref{eq:II.38.1} and noting that
\begin{align*}
\int_0^T \left( \min_{\eta\in W^*}\left( \Vert\zeta(t)-\eta\Vert^2_*+\Vert \eta\Vert^{q^*}_{W^*}\right)\right)^\frac{1}{\nu} \dd t&= \int_0^T \min_{\eta\in W^*}\left( \Vert\zeta(t)-\eta\Vert^2_*+\Vert \eta\Vert^{q^*}_{W^*}\right)^\frac{1}{\nu} \dd t\\
&\geq \int_0^T \min_{\eta\in W^*}\left( \Vert\zeta(t)-\eta\Vert^\frac{2}{\nu}_*+\Vert \eta\Vert^\frac{q^*}{\nu}_{W^*}\right) \dd t,
\end{align*} we obtain the uniform boundedness of $(B_\tau)_{0<\tau\leq \tau^*}$ in $\rmL^\frac{2}{\nu}(0,T;V^*)+\rmL^\frac{q^*}{\nu}(0,T;W^*)$ by arguing in the same way as for Case \textbf{(a)}. This, together with the uniform bounds of $(f_\tau)_{0<\tau\leq \tau^*} $ in $\rmL^2(0,T;H)$ and $(\xi_\tau)_{0<\tau\leq \tau^*}$ in $\rmL^\infty(0,T;U^*+V^*)$ yields the uniform bound of $(\Vh')_{0<\tau\leq \tau^*}$ in $\rmL^{\min\lbrace 2,q^*\rbrace}(0,T;U^*+V^*)$ with respect to $\tau$. It remains to show the uniform convergences \eqref{eq:II.31.2}, which follow immediately from the uniform bounds of $(\Vh')_{0<\tau\leq \tau^*}\subset \rmL^{\min\lbrace 2,q^*\rbrace}(0,T;U^*+V^*)$ and $(\Vo)_{0<\tau\leq \tau^*}\subset \rmL^2(0,T;V)$ in the respective spaces together with the estimates
\begin{align*}
&\Vert \Uh(t)-\Uo(t)\Vert_V \leq \Vert \Uu(t)-\Uo(t)\Vert_V = \int_{\teu(t)}^{\teo(t)} \Vert \Vo(r)\Vert_V \dd r \quad \text{ and } \\
&\Vert \Vh(t)-\Vo(t)\Vert_{U^*+V^*} \leq \Vert \Vu(t)-\Vo(t)\Vert_{U^*+V^*}= \int_{\teu(t)}^{\teo(t)} \Vert \Vh'(r)\Vert_{U^*+V^*}\dd r 
\end{align*}
 for all $t\in[0,T]$.
\end{proof} 
\subsection{Compactness of the interpolations}
\label{se:compactness.1}
This section is devoted to the existence of convergent subsequences of the sequence of approximate solutions in some proper \textsc{Bochner} spaces in order to pass to the limit in the discrete inclusion \eqref{eq:EuLa} as the step size vanishes. As we will see, we will indeed obtain in the limit a solution to the \textsc{Cauchy} problem \eqref{eq:I.2}. For this purpose, we will make use of the compactness properties of bounded sets in reflexive and separable spaces with respect to the weak topology. We will elaborate on this in the next result.

\begin{lem}[Compactness] \label{le:LimitPass.1} Under the assumptions of Lemma
  \ref{le:DUEE1}, 
 let $(\tau_n)_{n \in \mathbb{N}}$ be a
  vanishing sequence of positive numbers and let $u_0\in D\cap V$ and $v_0\in V$. Then, there exists a
  subsequence, still denoted by $(\tau_n)_{n \in \mathbb{N}}$, a pair of functions $(u,\xi)$ with \begin{align*}
  u\in\rmC_{w}([0,T];U)\cap \rmH^1(0,T;V) \cap\rmW^{1,\infty}(0,T;H) \text{ and } \xi\in \rmL^\infty(0,T;U^*+V^*)
\end{align*} that satisfies $u\in \rmH^{2}(0,T;U^*+V^*)$ in the case \textbf{(a)} and $u\in \rmW^{1,q}(0,T;\Wt)\cap \rmW^{2,\min\lbrace 2,q^*\rbrace}(0,T;U^*+V^*)$ in the case \textbf{(b)} while fulfilling the initial values $u(0)=u_0$ in $U$ and $u'(0)=v_0$ in $H$ such that the following convergences hold
\begin{subequations}
\label{eq:LP.all.1}
\begin{align}
\label{eq:LP.uhu.weak}
\Uun,\Uon,\Uhn \overset{*}{\rightharpoonup} u \quad \text{in } &\rmL^{\infty}(0,T;U\cap V),\\
\label{eq:uhu.weak}
\Uhn(t),\Uun(t),\Uon(t) \rightharpoonup u(t) \quad \text{in } &U \,\, \text{ for all }t\in[0,T],\\
\label{eq:uhu.weakinV}
\Uun(t) \rightharpoonup u(t) \quad \text{in } &V \,\, \text{ for all }t\in[0,T],\\
\label{eq:B.uhh}
\Uun \rightarrow u \quad \text{in } &\rmL^r(0,T;\WW) \quad \text{for any } r\geq 1,\\
\label{eq:B.uhh.ptw}
\Uhn(t),\Uun(t),\Uon(t) \rightarrow u(t) \quad \text{in } &\WW \,\, \text{ for all }t\in[0,T],\\
\label{eq:LP.vo}
\Von, \Vun \overset{*}{\rightharpoonup} u' \quad \text{in } &\rmL^2(0,T;V)\cap \rmL^{\infty}(0,T;H),\\
\label{eq:LP.vo.strong}
\Von,\Vun \rightarrow u' \quad \text{in } &\rmL^p(0,T;H) \quad \text{for all }p\geq 1,\\
\label{eq:vo.ptw}
\Von(t), \Vun(t) \rightarrow u'(t) \quad \text{in } & H \,\, \text{ for a.e. }t\in(0,T),\\
\label{eq:vh.ptw.weak}
\Vun(t), \Von(t) \rightharpoonup u'(t) \quad \text{in } & H \,\, \text{ for all }t\in[0,T],\\
\label{eq:LP.xi.1}
\Xion \overset{*}{\rightharpoonup} \xi \quad \text{in } &\rmL^\infty(0,T;V_\lambda^*),\\
\label{eq:LP.f}
f_{\tau_n}\rightarrow f \quad \text{in } &\rmL^{2}(0,T;H),\\ \text{and in Case  \textbf{(a)}}\notag\\
\label{eq:LP.vhn'.(a)}
\Vhn' \rightharpoonup u'' \quad \text{in } &\rmL^2(0,T;U^*+V^*),
\\
\label{eq:LP.B.(a)}
B_{\tau_n}\rightarrow B(\cdot,u(\cdot),u'(\cdot))\quad \text{in } &\rmL^{2}(0,T;V^*),\\
 \text{and in Case  \textbf{(b)}}\notag\\
 \label{eq:LP.weak.von.(b)}
\Von \rightharpoonup u' \quad \text{in } &\rmL^q(0,T;W),\\
\label{eq:LP.strong.von.(b)}
\Von \rightarrow u' \quad \text{in } &\rmL^{\max\lbrace2, r\rbrace}(0,T;W) \quad \text{for any }r\in [1,q),\\
\label{eq:L1.strong.DPsi2}
D_G\Psi_2(\Von) \rightarrow D_G\Psi_2(u') \quad \text{in } &\rmL^{r}(0,T;W^*) \quad \text{for any }r\in [1,q^*),\\
\label{eq:LP.vhn'.(b)}
\Vhn' \rightharpoonup u'' \quad \text{in } &\rmL^{\min\lbrace 2,q^*\rbrace}(0,T;U^*+V^*),\\
\label{eq:LP.B.(b)}
B_{\tau_n}\rightarrow B(\cdot,u(\cdot),u'(\cdot))\quad \text{in } &\rmL^{2}(0,T;V^*)+\rmL^{q^*}(0,T;W^*),
\end{align} where  $B_\tau (t):=B(\teo_\tau(t),\Uu (t),\Vu(t)),\ t\in[0,T]$.
\end{subequations} 
\end{lem}
\begin{proof}[Proof]
Let $\overline{U}_\tau,\underline{U}_\tau, \widehat{U}_\tau,\overline{V}_\tau,\underline{V}_\tau, \widehat{V}_\tau, \xi_\tau$ as well as $f_\tau$ be the interpolations with the initial values $u_0\in D\cap V$ and $v_0 \in V$ as defined in
  \eqref{eq:Approx.U}-\eqref{eq:Approx.xif}. Since all spaces are supposed to be separable and reflexive, we note that if a \textsc{Banach} space $X$ is separable and reflexive, the spaces $\rmL^p(0,T;X)$ for $1<p<\infty$ are also separable and reflexive, whereas $\rmL^\infty(0,T;X)$ is the dual of the separable space $\rmL^1(0,T;X^*)$. So, as a consequence of the \textsc{Banach--Alaoglu} theorem, bounded sets in $\rmL^p(0,T;X),\,1<p<\infty$, and $\rmL^\infty(0,T;X)$ are relatively compact with respect to the weak and weak* topology, respectively. In view of the a priori estimates \eqref{eq:II.29.2} and \eqref{eq:II.30.2}, Assumption \ref{eq:cond.E1.3} implies that the sequence $(\overline{U}_{\tau_n})_{n\in \mathbb{N}}$ is bounded in $\rmL^\infty(0,T;U)$. Together with the bounds \eqref{eq:allbounds}, this already yields the existence of converging subsequences (denoted as before) fulfilling \eqref{eq:LP.uhu.weak}, \eqref{eq:LP.vo}, and \eqref{eq:LP.xi.1}. We remark that the limit functions can be identified with $u$ and $u'$ by standard arguments. In order to show \eqref{eq:B.uhh}, we make use of the \textsc{Lions--Aubin--Dubinski\v{i}} lemma, see, e.g., \cite[Theorem 1]{DreJun12CFPC}. The boundedness of the sequence of piecewise linear interpolations $(\underline{U}_{\tau_n})_{n\in \mathbb{N}}$ and the discrete derivatives $(\underline{V}_{\tau_n})_{n\in \mathbb{N}}$ uniformly in $\rmL^\infty(0,T;U)$ and $\rmL^\infty(0,T;H)$, respectively, yields directly the relative compactness in $\rmL^r(0,T;\WW)$ for all $r\geq 1$. In view of \eqref{eq:uhu.weak}, this implies the convergence \eqref{eq:B.uhh.ptw}. The just proven convergence is indeed only needed to deduce the strong convergence of the perturbation $B$, i.e., convergence \eqref{eq:LP.B.(a)} and \eqref{eq:LP.B.(b)}. Before showing this convergence, we first proceed with proving the pointwise weak convergence as stated in \eqref{eq:vh.ptw.weak}. First, we note that from $\Vhn\in \rmW^{1,1}(0,T;U^*+V^*)\hookrightarrow \rmC([0,T];U^*+V^*)$ and \eqref{eq:LP.vhn'.(a)} or \eqref{eq:LP.vhn'.(b)}, there holds $\Vhn(t)\rightharpoonup u'(t)$ in $U^*+V^*$ as $n\rightarrow \infty$ for all $t\in [0,T]$. Since $\Vhn(t)$ is uniformly bounded in $H$ for all $t\in [0,T]$, it is (up to a subsequence) weakly convergent in $H$ to $u'(t)$. Since the weak limit is unique in $U^*+V^*$, we obtain, with the subsequence principle, the convergence of the whole sequence. Together with the strong convergence in \eqref{eq:II.31.2}, this implies \eqref{eq:vh.ptw.weak}. With the same argument, we deduce the pointwise weak convergences \eqref{eq:uhu.weak} and \eqref{eq:uhu.weakinV} where in the latter convergence, we use the fact that $u_0\in D\cap V$. Further, we recall that $\rmL^{\infty}(0,T;X)\cap \rmC_w([0,T];Y)=\rmC_w([0,T];X)$ for two \textsc{Banach} spaces $X$ and $Y$ with $X$ being reflexive and such that the continuous and dense embedding $X\overset{d}{\hookrightarrow} Y$ holds, see, e.g., in \textsc{Lions \& Magenes} \cite[Lemma 8.1, p. 275]{LioMag72NHBV}. Applying the latter result to $X=U$ and $Y=H$, there holds $u\in \rmC_w([0,T];U)$. Now, we seek to apply the \textsc{Lions--Aubin--Dubinski\v{i}} lemma to the sequence $(\Von)_{n\in \mathbb{N}}$ with $X=V,\, B=H$ and $Y=U^*+V^*$ in order to show the strong convergence in $\rmL^2(0,T;H)$ to the limit $u'\in \rmL^2(0,T;H)$. The Assumption of the \textsc{Lions--Aubin--Dubinski\v{i}} lemma \cite[Theorem 1]{DreJun12CFPC} follows for $p=2$ and $r=\min \lbrace 2,q^*\rbrace>1$ directly from the a priori estimate \eqref{eq:II.29.2} and the following estimate
\begin{align*}
 \Vert \sigma_{\tau_n}\Von-\Von\Vert_{\rmL^q(0,T-{\tau_n};U^*+V^*)}=\tau_n \Vert \Vhn'\Vert_{\rmL^q(0,T;U^*+V^*)}\leq \tau_nM \quad \text{for all }n\in \mathbb{N},
\end{align*} from which the strong convergence $\Von\rightarrow u'$ in $\rmL^2(0,T;H)$ as $n\rightarrow \infty$ follows. Taking into account the boundedness of the very same sequence in $\rmL^\infty(0,T;H)$, we obtain by a well-known interpolation inequality the strong convergence in $\rmL^r(0,T;H)$ for all $r\geq 1$, i.e. \eqref{eq:LP.vo.strong}. This, in turn, implies pointwise convergence of the very sequence almost everywhere in $(0,T)$, i.e., \eqref{eq:vo.ptw}. The assertion for $\Vhn$ can be shown analogously. Recalling the fact that the space of continuous functions $\rmC([0,T];H)$ is dense in $\rmL^2(0,T;H)$, for every $\epsilon>0$ there exists a function $f^\varepsilon\in \rmC([0,T];H)$ such that $\Vert f^\varepsilon-f\Vert_{\rmL^2(0,T;H)}<\varepsilon/3$. In view of this approximation property and defining $f_{\tau_n}^\varepsilon(t)=\frac{1}{\tau_n}\int_{t_{n_1}}^{n}f^\varepsilon(\sigma)\dd \sigma, \, t\in [t_{n-1},t_n), n=1,\dots,N$, we find 
\begin{align*}
\Vert f_{\tau_n}-f\Vert_{\rmL^2(0,T;H)}&\leq \Vert f_{\tau_n}-f^\varepsilon_{\tau_n} \Vert_{\rmL^2(0,T;H)}+\Vert f^\varepsilon_{\tau_n}-f^\varepsilon\Vert_{\rmL^2(0,T;H)}+\Vert f^\varepsilon-f\Vert_{\rmL^2(0,T;H)}\\
&\leq\Vert f-f^\varepsilon \Vert_{\rmL^2(0,T;H)}+\Vert f^\varepsilon_{\tau_n}-f^\varepsilon\Vert_{\rmL^2(0,T;H)}+\Vert f^\varepsilon-f\Vert_{\rmL^2(0,T;H)}\\
&\leq \varepsilon/3+\varepsilon/3+\varepsilon/3=\varepsilon
\end{align*} for sufficiently small step sizes $\tau_n$, where we also used the estimate \eqref{eq:f.est} for the first term, and where we made the second term smaller than $\varepsilon/3$ for sufficiently small step sizes which follows from the uniform continuity of $f^{\varepsilon}$. We proceed by showing the convergences which differ from each other in Case \textbf{(a)} and in Case \textbf{(b)}.\\
Ad case \textbf{(a)}. The weak convergence $\Vhn'\rightarrow u''$ as $n\rightarrow \infty$ in $\rmL^2(0,T;U^*+V^*)$ follows immediately from the reflexivity of the space $\rmL^2(0,T;U^*+V^*)$ and the uniform bound of the sequence $(\Vhn')_{n\in \mathbb{N}}$ in the very same space with respect to $n\in \mathbb{N}$. Further, we denote by $\calB(u)(t)=B(t,u(t),u'(t)),t\in[0,T]$, the associated \textsc{Nemitski\v{i}} operator and recall that $B_{\tau_n}(t)=B(\teon(t),\Uun(t),\Vun(t)), t\in[0,T]$. In order to show the strong convergence of the perturbation, we first note that from the uniform convergence \eqref{eq:B.uhh.ptw} and the pointwise convergence \eqref{eq:vo.ptw} together with the continuity condition \ref{eq:B1.1} implies
\begin{align}\label{eq:B.conv}
\Vert \calB_{\tau_n}(t)- \calB(u)(t) \Vert{V^*} \rightarrow 0 \quad \text{a.e. in }(0,T)
\end{align} as $n\rightarrow \infty$. By the growth condition \ref{eq:B1.2}, we also obtain $\calB(u)\in \rmL^\frac{2}{\nu}(0,T;V^*)$ so that we have $B_{\tau_n}-\calB(u)\in \rmL^\frac{2}{\nu}(0,T;V^*)$ being uniformly bounded with respect to $n\in \mathbb{N}$. Using \textsc{Egorov}'s theorem, it is easy to deduce the strong convergence of $B_{\tau_n}\rightarrow \calB(u)$ in $\rmL^{q}(0,T;V^*)$ as $n\rightarrow \infty$ for all $0<q<\frac{2}{\nu}$, and since $2<\frac{2}{\nu}$, we can choose $q=2$, i.e., \eqref{eq:LP.B.(a)}.\\
Ad case \textbf{(b)}. From the boundedness of the sequences of $(\Von)_{n\in \mathbb{N}}$ and $(\Vhn')_{n\in \mathbb{N}}$ in $\rmL^r(0,T;W)$ and $\rmL^{\min\lbrace 2,q^*\rbrace}(0,T;U^*+V^*)$, respectively, we obtain the weak convergences \eqref{eq:LP.weak.von.(b)} and \eqref{eq:LP.vhn'.(b)}. Applying again the \textsc{Lions--Aubin--Dubinski\v{i}} lemma to the sequence $(\Von)_{n\in \mathbb{N}}$ with the choices $X=V,\, B=W$ and $Y=U^*$ with $p=2$ and $r=1$ yields compactness of the sequence in $\rmL^2(0,T;W)$, and if $q> 2$, we obtain compactness of the sequence in every intermediate space $\rmL^s(0,T;W)$ with $2\leq s<q$ between $\rmL^2(0,T;W)$ and $\rmL^q(0,T;W)$ by an interpolation inequality, and hence \eqref{eq:LP.strong.von.(b)}. With the same reasoning as for the perturbation, the latter convergence yields \eqref{eq:L1.strong.DPsi2} employing \textsc{Egorov}'s theorem, the continuity, and the growth condition for $D_G\Psi_2$ on $W$.
The strong convergence of the perturbation in the space $\rmL^2(0,T;V^*)+\rmL^{q^*}(0,T;W^*)$ is more delicate and is established as follows: first, by the a priori estimate \eqref{eq:bBb}, the sequence $B_{\tau_n}$ is bounded in $\rmL^\frac{2}{\nu}(0,T;V^*)+\rmL^\frac{q^*}{\nu}(0,T;W^*)\subset \rmL^2(0,T;V^*)+\rmL^{q^*}(0,T;W^*)$ by a constant denoted by $\tilde{M}>0$.  With the same reasoning as for the first case, we obtain the convergence \eqref{eq:B.conv} and $\calB(u) \in \rmL^\frac{2}{\nu}(0,T;V^*)+\rmL^\frac{q^*}{\nu}(0,T;W^*)$. We choose the constant $\tilde{M}>0$ such that $\tilde{M}\geq \Vert \calB(u)\Vert_{\rmL^\frac{2}{\nu}(0,T;V^*)+\rmL^\frac{q^*}{\nu}(0,T;W^*)}$. Second, defining the set 
\begin{align*}
G_n:=\lbrace \eta\in \rmL^\frac{q^*}{\nu}(0,T;W^*): &\Vert B_{\tau_n}-\calB(u)-\eta \Vert_{\rmL^\frac{2}{\nu}(0,T;V^*)}\leq 2\tilde{M},\\
&\Vert \eta \Vert_{\rmL^\frac{q^*}{\nu}(0,T;W^*)}\leq 2\tilde{M}\rbrace,
\end{align*} there holds 
\begin{align*}
2\tilde{M}\geq &\Vert B_{\tau_n}-\calB(u)\Vert_{\rmL^\frac{2}{\nu}(0,T;V^*)+\rmL^\frac{q^*}{\nu}(0,T;W^*)}\\
&= \inf_{\eta\in \rmL^\frac{q^*}{\nu}(0,T;W^*)}\max \lbrace \Vert B_{\tau_n}-\calB(u)-\eta \Vert_{\rmL^\frac{2}{\nu}(0,T;V^*)},\Vert \eta \Vert_{\rmL^\frac{q^*}{\nu}(0,T;W^*)}\rbrace\\
&=\inf_{ \eta\in G_n } \max \lbrace \Vert B_{\tau_n}-\calB(u)-\eta \Vert_{\rmL^\frac{2}{\nu}(0,T;V^*)},\Vert \eta \Vert_{\rmL^\frac{q^*}{\nu}(0,T;W^*)}\rbrace \quad\text{for all }n\in \mathbb{N},
\end{align*} which restricts the set of functions where the infimum is taken over.
Then, by \textsc{Egorov}'s theorem, for every $\varepsilon>0$ there exists a subset $E\subset [0,T]$ with measure $\mu(E)<\varepsilon$ such that the uniform convergence 
\begin{align*}
\lim_{n\rightarrow \infty}\sup_{t\in [0,T]\backslash E}\Vert B_{\tau_n}(t)- \calB(u)(t) \Vert_*=0.
\end{align*} holds. Now, let $\eta:[0,T] \rightarrow V^*$ be any measurable function chosen to be fixed. Then, for every $\varepsilon>0$ there exists an index $N\in \mathbb{N}$ such that for all $n\geq N$, there holds
\begin{align*}
\Vert B_{\tau_n}(t)- \calB(u)(t)-\eta(t)\Vert_*\leq \varepsilon +\Vert \eta(t)\Vert_* \quad \text{for all } t\in [0,T]\backslash E.
\end{align*} Invoking the latter estimate, we obtain
\begin{align*}
&\Vert B_{\tau_n}-\calB(u)-\eta\Vert_{\rmL^2(0,T;V^*)}\\
&\leq \left(  \int_E \Vert B_{\tau_n}(t)-\calB(u)(t)-\eta(t) \Vert_*^2\dd t\right)^{\frac{1}{2}}+\left(\int_{[0,T]\backslash E} \Vert B_{\tau_n}(t)-\calB(u)(t)-\eta(t) \Vert_*^2\dd t\right)^{\frac{1}{2}}\\
&\leq \mu(E)^{1-\nu} \left(  \int_E \Vert B_{\tau_n}(t)-\calB(u)(t)-\eta(t) \Vert_*^\frac{2}{\nu}\dd t\right)^\frac{\nu}{2}+\left(\int_{[0,T]\backslash E}  (\varepsilon+\Vert \eta(t)\Vert_*)^2\dd t\right)^\frac{1}{2}\\
&\leq \varepsilon^{1-\nu}2\tilde{M}+(2T)^\frac{1}{2}\varepsilon +\left(\int_{[0,T]\backslash E}  2\Vert \eta(t)\Vert_*^2\dd t\right)^\frac{1}{2}\\
&\leq \varepsilon^{1-\nu}2\tilde{M}+(2T)^\frac{1}{2}\varepsilon +2\Vert \eta\Vert_{\rmL^{2}(0,T;V^*)}
\end{align*} for all $\eta\in \rmL^\frac{q^*}{\nu}(0,T;W^*)\subset \rmL^{q^*}(0,T;W^*)$ with $\Vert B_{\tau_n}-\calB(u)-\eta \Vert_{\rmL^\frac{2}{\nu}(0,T;V^*)}\leq 2\tilde{M}$. Finally, we end up with
\begin{align*}
&\Vert B_{\tau_n}-\calB(u)\Vert_{\rmL^2(0,T;V^*)+\rmL^{q^*}(0,T;W^*)}\\
&=\inf_{\eta\in \rmL^{q^*}(0,T;W^*)}\max \lbrace \Vert B_{\tau_n}-\calB(u)-\eta \Vert_{\rmL^2(0,T;V^*)},\Vert \eta \Vert_{\rmL^{q^*}(0,T;W^*)}\rbrace\\
&=\inf_{ \eta \in G_n}\max \lbrace \Vert B_{\tau_n}-\calB(u)-\eta \Vert_{\rmL^2(0,T;V^*)},\Vert \eta \Vert_{\rmL^{q^*}(0,T;W^*)}\rbrace\\
&\leq \inf_{ \eta \in G_n} \max \lbrace \varepsilon^{1-\nu}2\tilde{M}+(2T)^\frac{1}{2}\varepsilon +2\Vert \eta\Vert_{\rmL^{2}(0,T;V^*)},\Vert \eta \Vert_{\rmL^{q^*}(0,T;W^*)}\rbrace\\
&\leq \varepsilon^{1-\nu}2\tilde{M}+(2T)^\frac{1}{2}\varepsilon \quad \text{for all }n\geq N,
\end{align*} and hence \eqref{eq:LP.B.(b)}. Finally, thanks to \eqref{eq:uhu.weak} and \eqref{eq:vh.ptw.weak}, the initial conditions are also fulfilled by $u$ and $u'$, and since $u_0\in D\cap V$, there holds $u\in \rmH^1(0,T;V)$,  which shows the existence of solutions to the regularized problem \eqref{eq:I.2.reg}
\end{proof}

\subsection{Proof of Theorem 3.4}
\label{se:Proof.1}

\begin{proof}[\unskip\nopunct]
We first show that the limit function obtained from the previous lemma is indeed a solution to the regularized \textsc{Cauchy} problem \eqref{eq:I.2.reg}. We note that we still keep the regularization parameter $\varepsilon>0$ fixed and suppress the dependence of the solution as well as the dissipation potential $\Psi_2$ on $\varepsilon$. Let $u_0\in D\cap V_\lambda,$ $v_0\in H$, and a vanishing sequence of step sizes $(\tau_n)_{n\in N}$ be given. We remark that for the estimate \eqref{eq:vo} and the solvability of the variational approximation scheme, we necessarily needed the initial data $u_0$ and $v_0$ to be in $U\cap V$ in order to solve the variational approximation scheme \eqref{eq:ApproxSc} and to make use of the growth condition of $B$ in  \ref{eq:B1.2} for the a priori estimates since the energy functional and the dissipation potential are defined on different spaces. We circumvent this problem via approximating $u_0\in D\cap V_\lambda$ and $v_0\in H$ by approximating sequences $(u_0^{k})_{k\in \mathbb{N}}\subset D\cap V$ and $(v_0^{k})_{k\in \mathbb{N}}\subset V$ such that $u_0^{k}\rightarrow u_0$ in $U$ and $v_0^{k}\rightarrow v_0$ in $H$ as $k\rightarrow \infty$, which exists by Condition \ref{eq:cond.E1.1}.  Henceforth, we assume $k\in \mathbb{N}$ to be fixed, and we define the interpolations associated with the initial values $u^k_0$ and $v_0^{k}$ as in the previous lemma while omitting the dependence on $k$ for notational convenience. Then, again by the previous lemma, we obtain after selecting a subsequence (not relabeled) of the interpolations, the existence of a limit function $u\in \rmL^{\infty}(0,T;U)\cap \rmH^1(0,T;V)\cap \rmW^{1,\infty}(0,T;H)$ with $u(0)=u_0$ in $U$  and $u'(0)=v_0^{k}$ in $H$ that satisfies $u\in \rmH^{2}(0,T;U^*+V^*)$ in Case \textbf{(a)} and $u\in \rmW^{1,q}(0,T;W)\cap \rmW^{\min\lbrace 2,q^*\rbrace}(0,T;U^*+V^*)$ in Case \textbf{(b)}, where again we omit the dependence of the limit function on $k$. Now, the inclusion \eqref{eq:dis.inc} fulfilled by the interpolations reads in the weak formulation
\begin{align*}
\int_0^T \langle f_{\ton}(r)-B_{\tau_n}(r)-\Vhn'(r)-\xi_{\ton}(r) -D_G\Psi(\Von(r)),w(r) \rangle_{(U^*+V^*)\times (U\cap V)} \dd r=0 \quad 
\end{align*} for all $w\in \rmL^r(0,T;U\cap V)$ with $r=2$ in Case \textbf{(a)} and $r=\max\lbrace 2,1+(q-1)/(1-\delta(q-1)) \rbrace$ for a fixed $\delta\in (0,q^*-1)$ in Case \textbf{(b)}, where again $B_{\ton}(r)=B(\teo(r),\Uon(r),\Von(r))$, $r\in [0,T]$ and $D_G\Psi(\Von(r))=A \Von(r)$ in Case \textbf{(a)} and $D_G\Psi(\Von(r))=A \Von(r)+D_G\Psi_2(\Von(r))$ in Case \textbf{(b)}. 
For the readers convenience, we confine ourselves to Case \textbf{(b)}, but remark that Case \textbf{(a)} can be treated in the exact same manner.\\
Ad case \textbf{(b)}. Since $\Psi_1(v)=a(v,v)$ is defined by a strongly positive quadratic form, the \textsc{Fr\'{e}chet} derivative is a linear bounded and strongly positive operator $A:V\rightarrow V^*$, which implies that the associated \textsc{Nemitski\v{i}} operator $\mathcal{A}:\rmL^2(0,T;V)\rightarrow \rmL^2(0,T;V^*)\hookrightarrow \rmL^{\min\lbrace 2,q^*\rbrace}(0,T;U^*)$ is well defined, linear, bounded, and strongly positive. Therefore, the \textsc{Nemitski\v{i}} operator is weak-to-weak continuous so that we can pass with $\tau_n\searrow 0$ to the limit as $n\rightarrow \infty$. The \textsc{G\^{a}teaux} derivative $D_G \Psi_2(\Von)$ is strongly convergent to $D_G \Psi_2(u')$ in $\rmL^{\min\lbrace 2, q^*-\delta\rbrace}(0,T;U^*+V^*)$ so that passing to the limit is also justified in this term. We are also allowed to pass to the limit as the step size vanishes in the terms $f_\ton$ and $B_\ton$ which, according to the previous lemma, converge to $f$ and $B(\cdot,u(\cdot),u'(\cdot))$ strongly in $\rmL^2(0,T;H)\hookrightarrow \rmL^{\min\lbrace 2,q^*\rbrace}(0,T;V^*)$ and  $\rmL^2(0,T;V^*)+\rmL^{q^*}(0,T;W^*)\hookrightarrow \rmL^{\min\lbrace 2,q^*\rbrace}(0,T;V^*)$ as $n\rightarrow \infty$, respectively. Also by the previous lemma, there holds $\Vhn'\rightharpoonup u''$ in $\rmL^{\min\lbrace 2,q^*\rbrace}(0,T;U^*+V^*)$ and $\Xion\rightharpoonup \xi$ in $\rmL^\infty(0,T;U^*+V^*)$. Thus, we are allowed to pass to the limit in the weak formulation in these terms as well. Then, by a well-known density argument and by the fundamental lemma of calculus of variations, we deduce
\begin{align}\label{eq:LDIS.Inclusion}
u''(t)+D_G\Psi(u'(t))+\xi(t)+B(t,u(t),u'(t))=f(t) \quad \text{ in }U^*+V^* \text{ a.e. in }(0,T). 
\end{align} We proceed by showing that $\xi(t)\in \partial_{U\cap V}\calE_t(u(t))$ in $U^*+V^*$ for almost every $t\in(0,T)$. To do so, we employ the closedness condition \ref{eq:cond.E1.6}. Since we have already shown that the conditions a)-c) are satisfied, it remains to show the conditions d) and e). Condition d) follows immediately from
\begin{align*}
 \Vert \sigma_{\tau_n}\Uon-\Uon\Vert_{\rmL^2(0,T-{\tau_n};V)}=\tau_n \Vert \Uhn'\Vert_{\rmL^2(0,T;V)}\leq \tau_nM,
\end{align*} within Case \textbf{(a)} and 
\begin{align*}
& \Vert \sigma_{\tau_n}\Uon-\Uon\Vert_{\rmL^2(0,T-{\tau_n};V)\cap \rmL^r(0,T-{\tau_n};W)}\\
 &= \Vert \sigma_{\tau_n}\Uon-\Uon\Vert_{\rmL^2(0,T-{\tau_n};V)}+\Vert \sigma_{\tau_n}\Uon-\Uon\Vert_{\rmL^r(0,T-{\tau_n};W)}\\
&\leq \tau_n \Vert \Uhn'\Vert_{\rmL^2(0,T-{\tau_n};V)\cap \rmL^r(0,T-{\tau_n};W)}\leq \tau_nM,
\end{align*} in Case \textbf{(b)}. Condition e) in turn, is verified by the following calculations: let $t\in [0,T]$, then we have
\begin{align*}
&\int_0^{\teon (t)} \langle \Xion(r),\Uon(r)\rangle_{V_\lambda^*\times V_\lambda} \dd r\\
&=\int_0^{\teon (t)} \langle S_{\ton}(r)-A\Von(r)-\Vhn'(r)-D_G\Psi(\Von(r)),\Uon(r)\rangle_{V_\lambda^*\times V_\lambda} \dd r\\
&=\int_0^{\teon (t)} \langle S_{\ton}(r),\Uon(r)\rangle_{V^*\times V} \dd r\\
&-\int_0^{\teon (t)} \langle\Vhn'(r),\Uon(r)\rangle_{(U^*+V^*)\times (U \cap V)} \dd r\\
&-\int_0^{\teon (t)}\langle A\Von(r),\Uon(r)\rangle_{V^*\times V} \dd r\\
&-\int_0^{\teon (t)}\langle D_G\Psi_2(\Von(r)),\Uon(r)\rangle_{V^*\times V} \dd r\\
&=:I^n_1(t)+I^n_2(t)+I^n_3(t)+I^n_4(t).
\end{align*} The convergence of the first integral is due to the strong convergence of $S_{\ton}=f_\ton-B_\ton$ to $f-\calB(u)$ in $\rmL^{\min \lbrace 2,q^* \rbrace}(0,T;V^*)$ and the weak* convergence of $\Uon\rightharpoonup u$ in $\rmL^{\infty}(0,T;U\cap V)$ as $n\rightarrow \infty$. For the second integral, we recall the discrete integration by parts formula: let $n\in \mathbb{N}$ and $v^k,u^k\in H, k=0,\dots,N$. Then, there holds
\begin{align*}
\sum_{k=1}^n (v^k-v^{k-1},u^k)=(v^n,u^n)-(v^0,u^0)-\sum_{k=1}^n(v^{k-1},u^k-u^{k-1}).
\end{align*} 
Employing the discrete integration by parts formula, we obtain
\begin{align}\label{disc.int}
&-\int_0^{\teon (t)} \langle \Vhn'(r),\Uon(r)\rangle_{(U^*+V^*)\times (U\cap V)} \dd r \notag\\
&=\int_{0}^{\teon (t)} (\Vun(r),\Von(r))\dd r-(\Von(t),\Uon(t))+(v_0,u_0).
\end{align} Thus, by \eqref{eq:II.31.2}, \eqref{eq:uhu.weak}, \eqref{eq:LP.vo.strong} and \eqref{eq:vh.ptw.weak}
\begin{align*}
\lim_{n\rightarrow \infty}I_2^n(t)=\int_0^t(u'(r),u'(r))\dd r-(u'(t),u(t))+(v_0,u_0) \quad \text{for all }t\in [0,T].
\end{align*} Employing the more general integration by parts formula for \textsc{Bochner} spaces from Lemma A.1 in \textsc{Emmrich \& {\v{S}}i{\v{s}}ka} \cite{EmmSis13EESO} with $a=u$ and $b=u'$, we obtain
\begin{align*}
\int_0^t(u'(r),u'(r))\dd r-(u'(t),u(t))+(v_0,u_0)=-\int_0^t\langle u''(r),u(r)\rangle_{(U^*+V^*)\times (U\cap V)}\dd r 
\end{align*} for all $t\in [0,T]$.
We proceed with showing the convergence of the third integral $I_3^n(t)$. To do so, we use the symmetry of $A$ and the convexity of $\Psi_1$, to obtain
\begin{align}\label{eq:limit.Psi1}
-\int_0^{\teon (t)}\langle A \Von(r),\Uon(r)\rangle_{V^*\times V} \dd r&=-\int_0^{\teon (t)}\langle A \Uon(r),\Von(r)\rangle_{V^*\times V} \dd r \notag \\
&=-\sum_{k=1}^m \langle AU^k_{\tau_n},U^k_{\tau_n}-U^{k-1}_{\tau_{n}} \rangle_{V^*\times V}\notag\\
&\leq- \sum_{k=1}^m\left( \Psi_1(U^k_{\tau_n})-\Psi_1(U^{k-1}_{\tau_n})\right)\notag\\
&=\Psi_1(u_0)-\Psi_1(U^m_{\tau_n})\notag\\
&=\Psi_1(u_0)-\Psi_1(\Uon(t)).
\end{align} for $m\in \lbrace 1,\dots,N\rbrace$. Furthermore, we observe that
\begin{align}\label{eq:limit.Psi2}
\frac{\dd }{\dd t}\Psi_1(u(t))=\frac{\dd }{\dd t}\frac{1}{2} a(u (t),u(t))=\langle Au (t),u'(t)\rangle_{V^*\times V} \quad \text{for a.e. }t\in (0,T),
\end{align} which follows from the properties of $A$ and the fact that $u\in \rmH^1(0,T;V)$.
Then, taking into account \eqref{eq:limit.Psi1}, \eqref{eq:limit.Psi2}, the weak lower semicontinuity of $\Psi_1$, the pointwise weak convergence \eqref{eq:uhu.weakinV} as well as the symmetry of $A$, we obtain
\begin{align*}
\limsup_{n\rightarrow \infty}I^n_3(t)&\leq \limsup_{n\rightarrow \infty}\left(\Psi_1(u_0)-\Psi_1(\Uon(t))\right)\notag \\
&=-\liminf_{n\rightarrow \infty}\left(\Psi_1(\Uon(t))-\Psi_1(u_0)\right)\notag\\
&\leq \Psi_1(u_0)-\Psi_1(u(t))\notag\\
&= \int_0^t \langle Au(r),u'(r)\rangle_{V^*\times V}\dd r\notag\\
&= \int_0^t \langle Au'(r),u(r)\rangle_{V^*\times V}\dd r.
\end{align*} In view of \eqref{eq:LP.uhu.weak} and \eqref{eq:L1.strong.DPsi2}, we obtain for the last integral
\begin{align*}
\lim_{n\rightarrow \infty} I^n_4(t)&=-\lim_{n\rightarrow \infty}\int_0^{\teon (t)}\langle D_G\Psi_2(\Von(r)),\Uon(r)\rangle_{(U^*+V^*)\times (U\cap V)} \dd r\\
&=\int_0^t \langle D_G \Psi_2(u'(r)),u(r)\rangle_{(U^*+V^*)\times (U\cap V)} \dd r.
\end{align*} We end up with 
\begin{align*}
\limsup_{n\rightarrow \infty}\int_0^{\teon (t)} \langle \Xion(r),\Uon(r)\rangle_{V_\lambda^*\times V_\lambda} \dd r\leq \int_0^{t} \langle \xi(r),u(r)\rangle_{V_\lambda^*\times V_\lambda} \dd r
\end{align*} and thus
\begin{align*}
\limsup_{n\rightarrow \infty}\int_0^{T} \langle \Xion(r)-\xi(r),\Uon(r)-u(r)\rangle_{V_\lambda^*\times V_\lambda} \dd r \leq 0.
\end{align*} It remains to show the strong convergence $\Uon-u_0^n\rightarrow u-u_0$ in $\rmL^2(0,T;V)$ as $n\rightarrow \infty$ in order to obtain the conclusions of Assumption \ref{eq:cond.E1.6}. We show equivalently that $(\Uon-u_0^n)_{n\in \mathbb{N}}$ is a \textsc{Cauchy} sequence in $\rmL^2(0,T;V)$. To do so, we follow the idea of the proof of Lemma 4.6 in \textsc{Emmrich} \& \textsc{{\v{S}}i{\v{s}}ka} \cite{EmmSis13EESO} and consider in the first step
\begin{align*}
&\frac{\dd}{\dd t}\Psi_1\left(\hat{U}_{\tau_l}(t)-\hat{U}_{\tau_m}(t)\right)\\
&=\langle A(\hat{U}_{\tau_l}(t)-\hat{U}_{\tau_m}(t)),\overline{V}_{\tau_l}(t)-\overline{V}_{\tau_m}(t)\rangle_{V^*\times V}\\
&=\langle A(\overline{V}_{\tau_m}(t)-\overline{V}_{\tau_l}(t)),\hat{U}_{\tau_l}(t)-\hat{U}_{\tau_m}(t)\rangle_{V^*\times V}\\
&=\langle A(\overline{V}_{\tau_m}(t)-\overline{V}_{\tau_l}(t)),\overline{U}_{\tau_l}(t)-\overline{U}_{\tau_m}(t)\rangle_{V^*\times V}\\
&\quad+\langle A(\overline{V}_{\tau_m}(t)-\overline{V}_{\tau_l}(t)),\hat{U}_{\tau_l}(t)-\overline{U}_{\tau_l}(t)-\hat{U}_{\tau_m}(t)+\overline{U}_{\tau_m}(t)\rangle_{V^*\times V}\\
&=\langle \xi_{\tau_m}(t)-\xi_{\tau_l}(t)+\hat{V}_{\tau_m}'(t)-\hat{V}_{\tau_l}'(t)+S_{\tau_m}(t)-S_{\tau_l}(t)\\
&\quad -D_G\Psi_2(\overline{V}_{\tau_l}(t))+D_G\Psi_2(\overline{V}_{\tau_m}(t)),\overline{U}_{\tau_l}(t)-\overline{U}_{\tau_m}(t)\rangle_{V^*\times V}+b_{l,m}(t)\\
&=\langle \xi_{\tau_m}(t)-\xi_{\tau_l}(t),\overline{U}_{\tau_l}(t)-\overline{U}_{\tau_m}(t)\rangle_{V^*\times V}+ \langle \hat{V}_{\tau_m}'(t)-\hat{V}_{\tau_l}'(t),\overline{U}_{\tau_l}(t)-\overline{U}_{\tau_m}(t)\rangle_{V^*\times V}\\
&\quad+\langle -D_G\Psi_2(\overline{V}_{\tau_l}(t))+D_G\Psi_2(\overline{V}_{\tau_m}(t)),\overline{U}_{\tau_l}(t)-\overline{U}_{\tau_m}(t)\rangle_{V^*\times V}\\
&\quad+\langle S_{\tau_m}(t)-S_{\tau_l}(t),\overline{U}_{\tau_l}(t)-\overline{U}_{\tau_m}(t)\rangle_{V^*\times V}+b_{l,m}(t)\\
&\leq \lambda\Vert \overline{U}_{\tau_l}(t)-\overline{U}_{\tau_m}(t)\Vert^2_V + \langle \hat{V}_{\tau_m}'(t)-\hat{V}_{\tau_l}'(t),\overline{U}_{\tau_l}(t)-\overline{U}_{\tau_m}(t)\rangle_{V^*\times V}\\
&\quad+\langle -D_G\Psi_2(\overline{V}_{\tau_l}(t))+D_G\Psi_2(\overline{V}_{\tau_m}(t)),\overline{U}_{\tau_l}(t)-\overline{U}_{\tau_m}(t)\rangle_{V^*\times V}\\
&\quad+\langle S_{\tau_m}(t)-S_{\tau_l}(t),\overline{U}_{\tau_l}(t)-\overline{U}_{\tau_l}(t)\rangle_{V^*\times V}+b_{l,m}(t)\\
&\leq 2\lambda\Vert \hat{U}_{\tau_l}(t)-\hat{U}_{\tau_m}(t)\Vert^2_V+2\lambda\Vert \overline{U}_{\tau_l}(t)-\hat{U}_{\tau_l}(t)-\overline{U}_{\tau_m}(t)+\hat{U}_{\tau_m}(t)\Vert^2_V\\
&\quad+ \langle \hat{V}_{\tau_m}'(t)-\hat{V}_{\tau_l}'(t),\overline{U}_{\tau_l}(t)-\overline{U}_{\tau_m}(t)\rangle_{V^*\times V}\\
&\quad+\langle -D_G\Psi_2(\overline{V}_{\tau_l}(t))+D_G\Psi_2(\overline{V}_{\tau_m}(t)),\overline{U}_{\tau_l}(t)-\overline{U}_{\tau_m}(t)\rangle_{V^*\times V}\\
&\quad+\langle S_{\tau_m}(t)-S_{\tau_l}(t),\overline{U}_{\tau_l}(t)-\overline{U}_{\tau_m}(t)\rangle_{V^*\times V}+b_{l,m}(t)\\
&\leq \frac{2\lambda}{\mu}\Psi_1(\hat{U}_{\tau_l}(t)-\hat{U}_{\tau_m}(t))+2\lambda\Vert \overline{U}_{\tau_l}(t)-\hat{U}_{\tau_l}(t)-\overline{U}_{\tau_m}(t)+\hat{U}_{\tau_m}(t)\Vert^2_V\\
&\quad+ \langle \hat{V}_{\tau_m}'(t)-\hat{V}_{\tau_l}'(t),\overline{U}_{\tau_l}(t)-\overline{U}_{\tau_m}(t)\rangle_{V^*\times V}\\
&\quad+\langle -D_G\Psi_2(\overline{V}_{\tau_l}(t))+D_G\Psi_2(\overline{V}_{\tau_m}(t)),\overline{U}_{\tau_l}(t)-\overline{U}_{\tau_m}(t)\rangle_{V^*\times V}\\
&\quad+\langle S_{\tau_m}(t)-S_{\tau_l}(t),\overline{U}_{\tau_l}(t)-\overline{U}_{\tau_m}(t)\rangle_{V^*\times V}+b_{l,m}(t)
\\&=\frac{2\lambda}{\mu}\Psi_1\left(\hat{U}_{\tau_l}(t)-\hat{U}_{\tau_m}(t)\right)+c_{l,m}(t)
\end{align*} for almost every $t\in (0,T)$, where we have used the symmetry and strong positivity of $A$, the $\lambda$-convexity of $\calE$, and that \eqref{eq:dis.inc} is fulfilled. Then, by \textsc{Gronwall}'s lemma, there holds
\begin{align*}
\Psi_1(\hat{U}_{\tau_l}(t)-\hat{U}_{\tau_m}(t))\leq c_{l,m}(t)+\int_0^t \frac{2\lambda}{\mu} c_{l,m}(r)e^{\frac{2\lambda}{\mu}(t-r)}\dd r.
\end{align*} Integrating the latter inequality from $t=0$ to $t=T$ and using the strong positivity of $\Psi$ yields
\begin{align*}
\mu\int_0^T\Vert \hat{U}_{\tau_l}(t)-\hat{U}_{\tau_m}(t)\Vert_V^2 \dd t\leq \int_0^T c_{l,m}(t)\dd t+\int_0^T \int_0^t \frac{2\lambda}{\mu} c_{l,m}(r)e^{\frac{2\lambda}{\mu}(t-r)}\dd r\dd t.
\end{align*} Employing again the convergences \eqref{eq:II.31.2}, \eqref{eq:LP.uhu.weak}, \eqref{eq:LP.vo.strong}, \eqref{eq:vh.ptw.weak}, \eqref{eq:LP.f}, and \eqref{eq:L1.strong.DPsi2}-\eqref{eq:LP.B.(b)}, as well as the discrete integration by parts formula \eqref{disc.int}, we obtain 
\begin{align*}
\lim_{l,m\rightarrow \infty}\int_0^t c_{l,m}(r)\dd r&=0 \quad \, \text{for all }t\in [0,T],\\
\int_0^t c_{l,m}(r)\dd r&\leq C \quad \text{for all }l,m\in \mathbb{N}.
\end{align*}
Therefore, by the dominated convergence theorem, $(\Uhn-u_0^n)_{n\in \mathbb{N}}$ is a \textsc{Cauchy} sequence in $\rmL^2(0,T;V)$. By the convergence \eqref{eq:II.31.2}, we obtain that $(\Uon-u_0^n)_{n\in \mathbb{N}}$ is a \textsc{Cauchy} sequence in $\rmL^2(0,T;V)$ as well and thus convergent.
Hence, by the closedness condition \ref{eq:cond.E1.6}, there holds $\xi(t)\in \partial_{V_\lambda} \calE(u(t))$ as well as 
\begin{align}\label{eq:limsup.timed}
\calE_{\teon(t)}(\Uon(t))\rightarrow \calE_t(u(t)) \quad \text{and} \quad \limsup_{n\rightarrow \infty} \partial_t  \calE_{\tf_n(t)}(
\Uon(t))\leq \partial_t\calE_t(u(t))
\end{align} for a.e. $t\in(0,T)$. Now, we show that the energy-dissipation inequality holds. Let $t\in [0,T]$ and $\mathcal{N}\subset (0,T]$ be a set of measure zero such that $\calE_{\teon(s)}(\Uon(s))\rightarrow \calE_t(u(s))$ and $\Vun(s)\rightarrow u'(s)$ for each $s\in [0,T]\backslash \mathcal{N}$. Then, exploiting the convergences \eqref{eq:LP.all.1} and \eqref{eq:limsup.timed} as well as the condition \ref{eq:cond.E1.4}, we obtain from the discrete energy-dissipation inequality,
\begin{align*}
&\frac{1}{2}\vert u'(t) \vert^2 +\calE_t(u(t)) + \int_s^t \left(\Psi(u'(r))+\Psi^*(S(r)-\xi(r)-u''(r)) \right)\dd r \\
&\leq \liminf_{n\rightarrow \infty}\bigg (
 \frac{1}{2} \left \vert \Von(t)\right \vert^2+ \calE_{\teon(t)}(\overline{U}_\ton(t)) \\
 &\quad\quad \quad \quad \quad +\int_{s}^{t}\left( \Psi ( \Von (r) ) + \Psi^*\left(S_\ton(r)-\Vhn'(r)- \Xion(r)\right) \right) \dd r \bigg)
   \notag \\ 
  &\leq \limsup_{n\rightarrow \infty}\bigg ( \frac{1}{2} \left \vert \Von(s) \right \vert^2+
  \calE_{\teon(s)}(\Uon(s))+\int_{\teon(s)}^{\teon(t)} \partial_r
  \calE_r(\Uun (r))\dd r\\
  &\quad \quad \quad\quad\quad +\int_{\teon(s)}^{\teon(t)}
  \langle S_\ton(r), \Von(r) \rangle_{V^*\times V}\dd r
  +\tau \lambda\int_{\teo_\tau(s)}^{\teo_\tau(t)}\Vert \Vo(r) \Vert_V^2\dd r\bigg ) \notag\\
 &= \frac{1}{2}\vert u'(s) \vert^2 +\calE_s(u(s)) + \int_s^t  \partial_r \calE_r(u(r)) \dd r+\int_s^t \langle S(r),u'(r) \rangle_{V^*\times V}\dd r,
\end{align*} for all $t\in [0,T]$ if $s=0$ and almost every $s\in(0,t)$, where $S(r)=f(r)-B(r,u(r),u'(r))$. This shows that $u$ is a strong solution to \eqref{eq:I.2} satisfying the initial conditions $u_{k}(0)=u_0^k\in D\cap V$ and $u'_{k}(0)=v_0^{k}\in V, k\in \mathbb{N}$. We denote with $(u_{k})_{k\in \mathbb{N}}$ and $(\xi_{k})_{k\in \mathbb{N}}$ the associated solutions and subgradients of $\calE_t$ which satisfy \eqref{sol:SP}-\eqref{sol:EDI}. We recall that $u_0^{k}\rightarrow u_0$ in $U\cap V_\lambda$ and $v_0^{k}\rightarrow v_0$ in $H$ as $k\rightarrow \infty$.  The next steps are the same as before: 
\begin{itemize}
\item[1.] We derive a priori estimates based on the energy-dissipation inequality \eqref{sol:EDI},
\item[2.] We show compactness of the sequences $(u_k)_{k\in \mathbb{N}}$ and $(\xi_k)_{k\in \mathbb{N}}$ in appropriate spaces,
\item[3.] We pass to the limit in the inclusion \ref{sol:IC} as $k\rightarrow \infty$.
\end{itemize}
Ad 1. From the energy-dissipation inequality \eqref{sol:EDI} for $t\in [0,T]$ and $s=0$ while using the \textsc{Fenchel--Young} inequality, Condition \ref{eq:B1.2} and \ref{eq:cond.E1.4}, we obtain
\begin{align*}
&\frac{1}{2}\vert u_{k}'(t) \vert^2 +\calE_t(u_{k}(t)) + \int_0^t \left(\Psi(u_{k}'(r))+\Psi^*(S_{k}(r)-\xi_{k}(r)-u_{k}''(r)) \right)\dd r \\
&\leq \frac{1}{2}\vert v_0^{k} \vert^2 +\calE_0(u^k_0) + \int_0^t  \partial_r \calE_r(u_{k}(r)) \dd r+\int_0^t \langle S_{k}(r),u_{k}'(r) \rangle_{V^*\times V}\dd r\\
& \leq \frac{1}{2}\vert v_0^{k} \vert^2 +\calE_0(u^k_0) + C_1\int_0^t \calE_r(u_{k}(r)) \dd r\\
&\quad +\int_0^t \langle f(r)-B(r,u_{k}(r),u_{k}'(r)),u_{k}'(r) \rangle_{V^*\times V}\dd r\\
& \leq \frac{1}{2}\vert v_0^{k} \vert^2 +\calE_0(u^k_0) + C_1\int_0^t \calE_r(u_{k}(r)) \dd r+\int_0^t \left( \frac{1}{2} \vert f(r)\vert^2+\frac{1}{2}\vert u_{k}'(r)\vert^2\right) \dd r\\ 
&\quad +\int_0^t \left( c\,\Psi(u_{k}'(r))+c\,\Psi^*(-B(r,u_{k}(r),u_{k}'(r))/c)\right) \dd r\\
& \leq \frac{1}{2}\vert v_0^{k} \vert^2 +\calE_0(u^k_0)+\frac{1}{2}\Vert f\Vert^2_{\rmL^2(0,T;H)} + C_1\int_0^t \calE_r(u_{k}(r)) \dd r+\frac{1}{2}\int_0^t \vert u_{k}'(r)\vert^2 \dd r\\ 
&\quad +\int_0^t \left( c\Psi(u_{k}'(r))+\beta(1+\calE_r(u_{k}(r))+\vert u_{k}'(r)\vert^2 +\Psi^\nu(u_{k}'(r)))\right) \dd r\\
&\leq \frac{1}{2}\vert v_0^{k} \vert^2 +\calE_0(u^k_0)+\frac{1}{2}\Vert f\Vert^2_{\rmL^2(0,T;H)} +C T+ C\int_0^t \left( \calE_r(u_{k}(r)) +\frac{1}{2}\vert u_{k}'(r)\vert^2 \right) \dd r\\ 
&\quad +(c+\tilde{c})\int_0^t \Psi(u_{k}'(r))\dd r
\end{align*} for a constant $C=C(\nu,C_1,\beta)>0$, where $S_{k}(r):=f(r)-B(r,u_{k}(r),u_{k}'(r)),\, r\in[0,T]$ and $\beta\geq0, c\in (0,1)$, and $\tilde{c}>0$ such that $c+\tilde{c}\in (0,1)$. Taking into account the non-negativity of $\Psi,\Psi^*$, by the lemma of \textsc{Gronwall}, there exists a constant $C_B>0$ such that
\begin{align}\label{eq:apriori.bound}
&\frac{1}{2}\vert u_{k}'(t) \vert^2 +\calE_t(u_{k}(t)) + \int_0^t \left(\Psi(u_{k}'(r))+\Psi^*(S_{k}(r)-\xi_{k}(r)-u_{k}''(r)) \right)\dd r \leq C_B
\end{align} for all $t\in [0,T]$. \\
Ad 2. With the same reasoning as in Lemma \ref{le:LimitPass.1}, we find (up to a subsequence) the following convergences
\begin{subequations}
\label{eq:LP.2}
\begin{align}
\label{eq:LP.uhu.weak.2}
u_{k} \overset{*}{\rightharpoonup} u \quad \text{in } &\rmL^{\infty}(0,T;U),\\
\label{eq:LP.u.2}
u_{k}-u_0^k \overset{*}{\rightharpoonup} u-u_0 \quad \text{in } &\rmL^{\infty}(0,T;V),\\
\label{eq:LP.u.strong.2}
u_{k}-u_0^k \rightarrow u-u_0 \quad \text{in } &\rmL^{2}(0,T;V),\\
\label{eq:C.uhu.weak.2}
u_{k}(t) \rightharpoonup u(t) \quad \text{in } & U \,\, \text{ for all }t\in[0,T],\\
\label{eq:B.uhh.2}
u_k \rightarrow u \quad \text{in } &\rmL^r(0,T;\WW) \quad \text{for any } r\geq 1,\\
\label{eq:B.uhh.ptw.2}
u_k(t) \rightarrow u(t) \quad \text{in } &\WW \,\, \text{ for all }t\in[0,T],\\
\label{eq:LP.uhh.weak.2}
u_{k}' \overset{*}{\rightharpoonup} u' \quad \text{in } &\rmL^2(0,T;V)\cap\rmL^{\infty}(0,T;H),\\
\label{eq:LP.vo.strong.2}
u_{k}'\rightarrow u' \quad \text{in } &\rmL^p(0,T;H) \quad \text{for all }p\geq 1,\\
\label{eq:vo.ptw.2}
u_{k}'(t) \rightarrow u'(t) \quad \text{in } & H\,\, \text{ for a.e. }t\in(0,T),\\
\label{eq:vh.ptw.weak.2}
u_{k}'(t) \rightharpoonup u'(t) \quad \text{in } & H \,\, \text{ for all }t\in[0,T],\\
\label{eq:LP.xi.2}
\xi_{k} \overset{*}{\rightharpoonup} \xi \quad \text{in } &\rmL^\infty(0,T;U^*+V^*),\\
 \text{and in Case  \textbf{(a)}}\notag\\
\label{eq:LP.vhn'.(a).2}
u_{k}'' \rightharpoonup u'' \quad \text{in } &\rmL^2(0,T;U^*+V^*),
\\
\label{eq:LP.B.(a).2}
B(\cdot,u_{k},u_{k}')\rightarrow B(\cdot,u,u')\quad \text{in } &\rmL^{2}(0,T;V^*),\\
 \text{in Case  \textbf{(b)}}\notag\\
 \label{eq:LP.weak.von.(b).2}
u_{k}' \rightharpoonup u' \quad \text{in } &\rmL^q(0,T;W),\\
\label{eq:LP.strong.von.(b).2}
u_{k}' \rightarrow u' \quad \text{in } &\rmL^{\max\lbrace2, q-\varepsilon \rbrace}(0,T;W) \quad \text{for any }\varepsilon\in [1,q),\\
\label{eq:L1.strong.DPsi2.2}
D_G\Psi_2(u_{k}') \rightarrow D_G\Psi_2(u') \quad \text{in } &\rmL^{r}(0,T;W^*) \quad \text{for any }r\in [1,q^*),\\
\label{eq:LP.vhn'.(b).2}
u_{k}'' \rightharpoonup u'' \quad \text{in } &\rmL^{\min\lbrace 2,q^*\rbrace}(0,T;U^*+V^*),\\
\label{eq:LP.B.(b).2}
B(\cdot,u_{k},u_{k}')\rightarrow B(\cdot,u,u')\quad \text{in } &\rmL^{2}(0,T;V^*)+\rmL^{q^*}(0,T;W^*),
\end{align}
\end{subequations} except from the strong convergence \eqref{eq:LP.u.strong.2}, which needs to be proven.
Thus, we show that $(u_k-u_0^k)_{k\in \mathbb{N}}$ is a \textsc{Cauchy} sequence in $\rmL^2(0,T;V)$. To do so, we consider 
\begin{align*}
&\frac{\dd}{\dd t}\Psi_1(u_l(t)-u_0^l-u_m(t)+u_0^m)\\
&=\langle A(u_l(t)-u_0^l-u_m(t)+u_0^m),u_l'(t)-u_m'(t)\rangle_{V^*\times V}\\
&=\langle A(u_l'(t)-u_m'(t)),u_l(t)-u_0^l-u_m(t)+u_0^m\rangle_{V^*\times V}\\
&=\langle \xi_m(t)-\xi_l(t)+u''_m(t)-u''_l(t)+S_m(t)-S_l(t),u_l(t)-u_0^l-u_m(t)+u_0^m \rangle_{V^*\times V}\\
&= \langle \xi_m(t)-\xi_l(t),u_l(t)-u_0^l-u_m(t)+u_0^m)\rangle_{U^*\times U}\\
&\quad+\langle u''_m(t)-u''_l(t),u_l(t)-u_0^l-u_m(t)+u_0^m\rangle_{(U^*+V^*)\times (U\cap V)}\\
&\quad +\langle S_m(t)-S_l(t),u_l(t)-u_0^l-u_m(t)+u_0^m\rangle_{V^*\times V}\\
&\leq \langle \xi_m(t)-\xi_l(t),u_0^m-u_0^l\rangle_{U^*\times U}\\
&\quad+\langle u''_m(t)-u''_l(t),u_l(t)-u_0^l-u_m(t)+u_0^m\rangle_{(U^*+V^*)\times (U\cap V)}\\
&\quad +\langle S_m(t)-S_l(t),u_l(t)-u_0^l-u_m(t)+u_0^m\rangle_{V^*\times V}.
\end{align*} where we have taken into account that $u_k$ is a solution of \eqref{eq:I.2} and that the subdifferential operator $\partial \calE_t$ is monotone. Integrating the latter inequality and using the integration by parts rule yields 
\begin{align*}
&\mu \Vert u_l(t)-u_0^l-u_m(t)+u_0^m\Vert^2_V\\
&\leq \Psi_1(u_l(t)-u_0^l-u_m(t)+u_0^m)\\
&\leq \int_0^t\langle \xi_m(r)-\xi_l(r),u_0^m-u_0^l\rangle_{U^*\times U}\dd r\\
&\quad+\int_0^t \langle u''_m(r)-u''_l(r),u_l(r)-u_0^l-u_m(r)+u_0^m\rangle_{(U^*+V^*)\times (U\cap V)}\dd r\\
&\quad +\int_0^t\langle S_m(r)-S_l(r),u_l(r)-u_0^l-u_m(r)+u_0^m\rangle_{V^*\times V}\dd r\\
&= \int_0^t\langle \xi_m(r)-\xi_l(r),u_0^m-u_0^l\rangle_{U^*\times U}\dd r\\
&\quad+\int_0^t \vert u'_m(r)-u'_l(r)\vert^2 \dd r +(u'_m(t)-u'_l(t),u_l(t)-u_0^l-u_m(t)+u_0^m)\\
&\quad +\int_0^t\langle S_m(r)-S_l(r),u_l(r)-u_0^l-u_m(r)+u_0^m\rangle_{V^*\times V}\dd r.
\end{align*} From the strong convergence $u_0^k\rightarrow u_0$ in $U$ as $k\rightarrow \infty$ and in view of the convergences \eqref{eq:LP.2} and the a priori bound \eqref{eq:apriori.bound}, the right-hand side is uniformly bounded and convergent to zero for every $t\in [0,T]$ as $m,l\rightarrow \infty$. Thus, by the dominated convergence theorem, we conclude that $(u_k-u_0^k)_{k\in \mathbb{N}}$ is a \textsc{Cauchy} sequence in $\rmL^2(0,T;V)$, and therefore strongly convergent in $\rmL^2(0,T;V)$ with the limit $u-u_0$.\\
Ad 3. With the same argument as before, we show that the equation \eqref{eq:LDIS.Inclusion} is fulfilled. However, it remains to identify $\xi(t)\in \partial_U\calE_t(u(t))$ a.e. in $(0,T)$. But this follows from the following limsup estimate and the closedness condition \ref{eq:cond.E1.6}
\begin{align*}
&\limsup_{k\rightarrow \infty}\int_0^t \langle \xi_k(r)-\xi(t),u_k(t)-u(t)\rangle_{V_\lambda^*\times V_\lambda}\dd r\\
&=\limsup_{k\rightarrow \infty}\int_0^t \langle \xi_k(r)-\xi(t),u_k(t)-u_0^k-u(t)+u_0\rangle_{V_\lambda^*\times V_\lambda}\dd r\\
&=\limsup_{k\rightarrow \infty} \bigg( \int_0^t \langle u''(r)-u''_k(r),u_k(t)-u_0^k-u(t)+u_0\rangle_{(U^*+V^*)\times (U\cap V)}\dd r\\
&\quad+\int_0^t \langle B(r,u(r),u'(r))-B(r,u_k(r),u_k'(r)),u_k(t)-u_0^k-u(t)+u_0\rangle_{V^*\times V}\dd r\\
&\quad+\int_0^t \langle Au'(r)-Au_k'(r),u_k(t)-u_0^k-u(t)+u_0\rangle_{V^*\times V}\dd r\bigg)\\
&=\limsup_{k\rightarrow \infty} \bigg( \int_0^t  \vert u'(r)-u'_k(r)\vert^2 \dd r +(u'(t)-u'_k(t),u_k(t)-u_0^k-u(t)+u_0)\\
&\quad+\int_0^t \langle B(r,u(r),u'(r))-B(r,u_k(r),u_k'(r)),u_k(t)-u_0^k-u(t)+u_0\rangle_{V^*\times V}\dd r\\
&\quad+\int_0^t \langle A(u_k(t)-u_0^k-u(t)+u_0),u'(r)-u_k'(r)\rangle_{V^*\times V}\dd r\bigg)=0
\end{align*} which again follows from the convergences \eqref{eq:LP.2}. Thus, there holds $\xi(t)\in \partial_{V_\lambda}\calE_t(u(t))$ a.e. in $(0,T)$, which shows the existence of solutions to the regularized problem \eqref{eq:I.2.reg}.\\\\ Now, let $(\varepsilon_l)_{l\in \mathbb{N}}$ be a sequence of regularization parameters such that $\varepsilon_l\searrow 0$ as $l\rightarrow \infty$, and let to each $l\in \mathbb{N}$, $u_l$ be the associated solution to \eqref{eq:I.2.reg}. It is easy to see that by the same argumentation as before, we can extract a subsequence (denoted as before) such that the convergences \eqref{eq:LP.uhu.weak.2}-\eqref{eq:LP.strong.von.(b).2}, \eqref{eq:LP.vhn'.(b).2}, and \eqref{eq:LP.B.(b).2}. If we show that there exists a subsequence (denoted as before) and a function $\eta\in \rmL^{q^*}(0,T;W^*)$ such that 
\begin{align*}
D_G\Psi^{\varepsilon_l}_2(u_{l}') \rightharpoonup \eta \quad \text{in } \rmL^{r}(0,T;W^*),
\end{align*} and $\eta\in \partial \Psi_2(u') \text{ a.e. in } (0,T)$, then we are done since all the other steps can be shown in the same manner as before. First, we observe that by Lemma \ref{le:Leg.Fen} and the growth condition \eqref{eq:Psi.growth}, there exists constants $c_1,c_2,c_3>0$ independent of $\varepsilon$ such that 
\begin{align*}
c_1\Vert u'(t)\Vert_{W}^q +c_2\Vert \rmD_G \Psi_2^{\varepsilon_l}(u_l'(t))\Vert_{W^*}^{q^*}-c_3 &\leq \Psi_2^{\varepsilon_l}(u_l'(t))+\Psi_2^{\varepsilon_l,*}(\rmD_G \Psi_2^{\varepsilon_l}(u_l'(t)))\\
&= \langle \rmD_G \Psi_2^{\varepsilon_l}(u_l'(t)),u_l'(t)\rangle_{W^*\times W}\\
&\leq C(\delta)\Vert u'(t)\Vert_{W}^q +\delta \Vert \rmD_G \Psi_2^{\varepsilon_l}(u_l'(t))\Vert_{W^*}^{q^*}
\end{align*} for all $\delta>0$, where we have used \textsc{Young}'s inequality in the last step. For $\delta<c_2$, we infer that $(\rmD_G \Psi_2^{\varepsilon_l}(u_l'(t)))_{l\in \mathbb{N}}$ is bounded in $\rmL^{q^*}(0,T;W^*)$. Hence, there exists a weak limit $\eta\in \rmL^{q^*}(0,T;W^*)$ such that (up to a subsequence)
\begin{align*}
D_G\Psi^{\varepsilon_l}_2(u_{l}') \rightharpoonup \eta \quad \text{in } &\rmL^{r}(0,T;W^*).
\end{align*} We want to show that $\eta$ is a subgradient of $\Psi_2$. As we mentioned in Section \ref{se:regularized.problem}, $\Psi_2^{\varepsilon}\Mto \Psi_2$, i.e., for all $v\in W$
\begin{align*}
\begin{cases}
a)\quad  \Psi_2(v)\leq \liminf_{l\to \infty} \Psi_2^{\varepsilon_l}(v_l) \quad \text{for all }v_n\rightharpoonup v \text{ in } W,\\
b)\quad \exists \hat{v}_l\rightarrow v \text{ in $W$ such that } \Psi_2(v) \geq \limsup_{l\to \infty} \Psi_2^{\varepsilon_l}(\hat{v}_l).
\end{cases}
\end{align*} In fact, by \textsc{B.} \cite[Theorem 2.2]{Bach22GMYR}, the assertion $b)$ can be replaced by the even stronger statement
\begin{align*}
\forall \hat{v}_l\rightarrow v \text{ in $W$ there holds } \Psi_2(v)= \lim_{l\to \infty} \Psi_2^{\varepsilon_l}(\hat{v}_l).
\end{align*} By \textsc{Fatou}'s lemma and the dominated convergence theorem, we obtain
\begin{align*}
\int_0^T \Psi_2(u'(t))\dd t-\int_0^T \Psi_2(v(t))\dd t&\leq  \int_0^T\liminf_{l\rightarrow \infty}\Psi^{\varepsilon_l}_2(u_l'(t))\dd t-\int_0^T \lim_{l\rightarrow \infty} \Psi^{\varepsilon_l}_2(v(t))\dd t\\
&\leq \liminf_{l\rightarrow \infty}\left( \int_0^T \Psi^{\varepsilon_l}_2(u_l'(t))\dd t-\int_0^T \Psi^{\varepsilon_l}_2(v(t))\dd t\right)\\
&=\liminf_{l\rightarrow \infty}\int_0^T \langle \rmD_G \Psi^{\varepsilon_l}_2(u_l'(t)),u_l'(t)\rangle_{W^*\times W}\dd t\\
&=\int_0^T \langle \eta(t),u'(t)\rangle_{W^*\times W}\dd t.
\end{align*} By \textsc{Kenmochi} \cite[Proposition 1.1]{Kenm75NPVI}, it follows that $\eta(t)\in \partial \Psi_2(u'(t))$ a.e. in $(0,T)$ which completes the proof of Theorem \ref{th:MainExist1}.

\end{proof}

\begin{rem}\label{relax} If we take a closer look into the proof, we see that the assumption that $\calE_t$ is sequentially weakly lower semicontinuous has only been used to show the existence of solutions to the discrete problem and to show the energy-dissipation inequality. If we only address the existence of solutions without the energy-dissipation inequality, we can relax the condition by assuming (in both cases) that there exists $r_0>0$ such that $u\mapsto \frac{1}{r_0} a(u,u)+\calE_t(u)$ is sequentially weakly lower semicontinuous. The existence of discrete solutions under this assumption follows from the fact that
\begin{align*}
\frac{1}{\tau}a(u-u_0,u-u_0)+\calE_t(u)=\frac{1}{\tau}a(u,u)+\calE_t(u)-\frac{2}{\tau}a(u,u_0)+\frac{1}{\tau}a(u_0,u_0),
\end{align*} so that the first two terms are sequentially weakly lower semicontinuous and that the remaining terms are weak-to-weak continuous. 
\end{rem}

\begin{section}{Applications}\label{se:App}
Now, we want to provide some applications of the abstract theory. First, we discuss in detail a mathematical example to illustrate the strength of the theory and highlight the case (a) and (b). Then, we also discuss some physically meaningful examples.   

\subsection{Differential inclusion I}\label{se:app.DI1A} In the first example, we consider a system which can be treated in the Case \textbf{(a)} of the linearly damped inertial system, where the dissipation potential is given by the \textsc{Dirichlet} energy and the energy functional is a nonsmooth $\lambda$-convex function which to the best of the authors' knowledge can not be treated with the abstract results known thus far. More precisely, we consider the initial-boundary value problem 
\begin{align*} 
\text{(P1)}
\begin{cases}
\partial_{tt}\uu-\Delta \partial_t \uu-\Delta_p \uu+(\vert \uu\vert ^2-1)\uu - \nabla \cdot \mathbold{p}+\bb\left(\xx,t,\uu,\partial_t \uu\right) = \ff \quad \text{in } \Omega_T,\\
\mathbold{p}(\xx,t)\in \mathbf{Sgn}\left(\nabla \uu(\xx,t)\right)\quad\text{a.e. in } \Omega_T,\\
\uu(\xx,0)\,\,=\uu_0(\xx) \quad \text{on } \Omega,\\
\uu'(\xx,0)=\vv_0(\xx) \quad \,\text{on } \Omega, \\
\uu(\xx,t)\,\,=0 \quad  \qquad \text{on } \partial \Omega\times[0,T],
\end{cases}
\end{align*} where $\mathbf{{Sgn}}:\mathbb{R}^{d\times m} \rightrightarrows \mathbb{R}^{d\times m}$ is the multi-valued and multi-dimensional sign function defined by
\begin{align}\label{sign}
\mathbf{Sgn}(\AB)=
\begin{cases} B_{\mathbb{R}^{d\times m}}(0,1) \quad &\text{if }\AB=0 \\
\frac{\AB}{\vert \AB\vert } \quad &\text{otherwise,}
\end{cases}
\end{align}  $\ff:\Omega\rightarrow \mathbb{R}^m$ is an external force, $\bb:\Omega\times [0,T]\times \mathbb{R}^m\times \mathbb{R}^m\rightarrow \mathbb{R}^m$ a \textsc{Carath\'{e}odory} function in the sense that $\bb(\xx,\cdot,\cdot,\cdot)$ is continuous for almost every $\xx\in \Omega$ and $\bb(\cdot,t,\yy,\mathbold{z})$ is measurable for all $t\in [0,T]$ and $\yy,\mathbold{z}\in \mathbb{R}^m$. Furthermore, $\bb$ is assumed to satisfy the following growth condition: there exists a constant $C_b>0$ and numbers $q,r>1$ such that 
\begin{align*}
\vert \bb(\xx, t,\uu,\vv)\vert\leq C_b(1+ \vert \uu\vert^{q-1}+ \vert \vv\vert^{r-1}) \quad \text{for a.e. }\xx\in \Omega, t\in [0,T] \text{ and all }\uu,\vv\in \mathbb{R}^m.
\end{align*} Here, $p,q,r\geq 1$ are to be chosen in accordance with the assumptions.

Choosing the spaces $U=\rmW_0^{1,p}(\Omega)^m\cap \rmL^4(\Omega)^m,V=\rmH_0^{1}(\Omega)^m,\, W=\rmL^{\max\lbrace 2,q\rbrace}(\Omega)^m$ and $H= \rmL^2(\Omega)^m$  equipped with the standard norms, we assume $\mathbold{f}\in \rmL^2(0,T;V^*)$. The energy functional $\calE: V\rightarrow (-\infty,+\infty]$ and the dissipation potential $\Psi: V\rightarrow \mathbb{R}$ are given by
\begin{align*}
\calE(u)=
\begin{cases}
\int_{\Omega}\left( \frac{1}{p}\vert \nabla \uu(\xx)\vert^{p} +\vert \nabla \uu(\xx)\vert+\frac{1}{4}(\vert \uu(\xx)\vert^2-1)^2\right)\dd \xx \quad \text{if } \uu \in \DOM(\calE),\\
+\infty \quad \text{otherwise},
\end{cases}
\end{align*} and
\begin{align*}
\Psi(v)=\frac{1}{2}\int_{\Omega}\vert \nabla \vv(\xx)\vert^2 \dd \xx,
\end{align*}
 respectively, whereas the perturbation $B:[0,T]\times W\times H\rightarrow V^*$ is defined by
\begin{align*}
\langle B(t,\uu,\vv),\ww \rangle_{V^*\times V} &=\langle B(t,\uu,\vv),\ww \rangle_{W^*\times W}= \int_{\Omega} \bb(\xx,t,\uu(\xx),\vv(\xx)) \cdot \ww(\xx) \dd \xx.
\end{align*} The \textsc{Legendre--Fenchel} transformation   $\Psi^*:\rmH^{-1}(\Omega)^m \rightarrow \mathbb{R}$ of $\Psi$ is obviously given by $\Psi^*(\xii)=\frac{1}{2}\Vert \xii\Vert^2_{-1,2}$. Furthermore, it is readily seen that the energy functional is not \textsc{G\^{a}teaux} differentiable, and its effective domain is given by $\DOM(\calE)= \rmW_0^{1,p}(\Omega)^m\cap \rmL^4(\Omega)^m$. The values $p,q,r\geq 1$ are to be chosen such that all assumptions are fulfilled. We can choose, e.g.,    
 \begin{align*}
 &d=1,\, p\in (1,+\infty), r\in [1,2], q\in [1,p/2+1],\\
 &d=2,\, p\in (1,+\infty), r\in [1,2], q\in 
 \begin{cases} [1,pd/(p-d))\cap [1,p/2+1] \quad \text{if }p\in (1,2),\\
[1,p/2+1]\quad \text{if }p\geq 2,
 \end{cases}\\ 
&d\geq 3,\, p\in (1,+\infty),r\in [1,2], q\in 
 \begin{cases} [1,q^*) \quad \text{if }p\in (1,2),\\
[1,p/2+1]\quad \text{if }p\geq 3,
 \end{cases}
 \end{align*} where $q^*=\min\lbrace \frac{d(p+2)}{2(d-p)},\frac{3d+4}{d},p+1\rbrace$. Then, by the \textsc{Sobolev} embedding theorem and the \textsc{Rellich--Kondrachov} theorem, $U$ and $V$ are densely, continuously, and compactly embedded in $W$ and $H$, respectively. We will verify for illustration the assumptions for the case $d\geq 3$. Since the dissipation potential is state-independent, it is induced by the bilinear form $a:V\times V\rightarrow \mathbb{R}$,
 \begin{align*}
 a(\vv,\ww): =\frac{1}{2}\int_\Omega \nabla \vv\ccdot \nabla \ww \dd \xx,
 \end{align*} and therefore satisfies all conditions. The conditions \ref{eq:cond.E1.2}-\ref{eq:cond.E1.4} are obviously fulfilled by the energy functional. In order to verify \ref{eq:cond.E1.1}, we note that every convex and lower semicontinuous functional on a \textsc{Banach} space is weakly lower semicontinuous. Taking the latter into account, we observe that for $u\in \DOM(\calE)$, the energy functional
 \begin{align*}
 \calE(\uu)&=\int_{\Omega}\left( \frac{1}{p}\vert \nabla \uu(\xx)\vert^{p} +\vert \nabla \uu(\xx)\vert+\frac{1}{4}(\vert\uu(\xx)\vert^2-1)^2\right)\dd \xx \notag\\
&=\int_{\Omega}\left( \frac{1}{p}\vert \nabla \uu(\xx)\vert^{p}+\vert \nabla \uu(\xx)\vert +\frac{1}{4}(\vert \uu(\xx)\vert^4-2\vert \uu(\xx)\vert^2+1)\right)\dd \xx \notag \\
&=\int_{\Omega}\left( \frac{1}{p}\vert \nabla \uu(\xx)\vert^{p} +\vert \nabla \uu(\xx)\vert+\frac{1}{4}(\vert \uu(\xx)\vert^4+1)\right)\dd \xx-\frac{1}{2}\int_\Omega \vert \uu(\xx)\vert^2\dd \xx\notag \\
&=\calW(\uu)-\frac{1}{2}\int_\Omega \vert \uu(\xx)\vert^2\dd \xx
\end{align*} is the sum of a convex function $\calW$ and a concave function $u\mapsto -\frac{1}{2}\int_\Omega \vert \uu(\xx)\vert^2\dd \xx $ on $V$. The lower semicontinuity of $\calW$ on $V$ follows immediately from the converse of the dominated convergence theorem (see, e.g., \textsc{Br\'{e}zis} \cite[Theorem 4.9, p. 94]{Brez11FASS}) and \textsc{Fatou}'s lemma. Further, due to the compact embedding of $V$ in $H$, the concave function is continuous on $V$ with respect to the weak topology. This implies $\calE$ to be weakly lower semicontinuous on $V$. In fact, the convex part of the energy is perturbed by the negative \textsc{Hilbert} space norm of $H$ squared, which by the parallelogram law and the embedding $V\hookrightarrow H$ leads to the $\lambda$-convexity of $\calE$ with $\lambda:=C$ being the constant of the very same embedding. Now, we show the closedness property \ref{eq:cond.E1.6}. First, we note that for each $\uu\in D(\partial \calE)$, there holds $\xii\in \partial_U \calE(\uu)=\partial_U \calW(\uu)-\uu$ if and only if $\xii= -\Delta_p\uu+\nabla \cdot \mathbold{p}+(\vert \uu\vert^2-1)\uu\in U^*$ for a measurable selection $\mathbold{p}\in \rmL^\infty(\Omega)^{d\times m}$ satisfying $\mathbold{p}\in \mathbf{Sgn}(\nabla \uu)$ a.e. in $\Omega$. This can be seen as follows: we define the functionals $\calW_1:W_0^{1,p}(\Omega)^m \rightarrow [0,+\infty]$ and $\calW_2:\rmL^1(\Omega)^{d\times m}\rightarrow \mathbb{R}$ as well as the operator $\Lambda: W_0^{1,p}(\Omega)^m\rightarrow \rmL^1(\Omega)^{d\times m}$ via
\begin{align*}
\calW_1(u)=
\begin{cases}
\int_{\Omega}\left( \frac{1}{p}\vert \nabla \uu(\xx)\vert^{p} +\frac{1}{4}(\vert \uu(\xx)\vert^2-1)^2\right)\dd \xx \quad \text{if } \uu \in \DOM(\calE),\\
+\infty \quad \text{otherwise},
\end{cases}
\end{align*} 
\begin{align*}
\calW_2(\AB)=\int_\Omega \vert  \AB(x)\vert \dd \xx
\end{align*} and $\Lambda u=\nabla u$. We note that $\Lambda$ is linear and bounded and has as adjoint operator  $\Lambda^*: \rmL^\infty(\Omega)^{d\times m}\rightarrow W^{-1,p^*}(\Omega)^m, \AB\mapsto -\nabla\cdot \AB$ the divergence operator. Let $\uu\in \DOM(\partial_U\calW)$, then by Lemma \ref{le:Subdif2} and Lemma \ref{le:chain.rule}, there holds
\begin{align*}
\xii&\in \partial_U\left( \calW_1(\uu)+\calW_2(\Lambda \uu)\right)\\
&=\partial_U \calW_1(\uu)+\partial_U\calW_2(\Lambda \uu)\\&= \partial_U \calW_1(\uu)+\Lambda^*\partial_X\calW_2(\Lambda \uu),
\end{align*} where $X=\rmL^1(\Omega)^{m \times d}$. Thus, there exists $\xii_1\in\partial_U \calW_1(\uu) $ and $\xii_2\in +\Lambda^*\partial_X\calW_2(\Lambda \uu)$ such that $\xii=\xii_1+\xii_2$. Now, we shall determine $\xii_1$ and $\xii_2$. Since $\calW_1$ is \textsc{G\^{a}teaux} differentiable on $\rmW^{1,r}(\Omega)^m\cap \rmL^4(\Omega)^m$, we deduce immediately $\xii_1=-\Delta_p \uu+(\vert \uu\vert^2-1)\uu$ a.e. in $\Omega$. In order to determine $\xi_2$, we note first that $\xii_2=\nabla \cdot \mathbold{p}$ with $\mathbold{p}\in \partial_X \calW_2(\Lambda \uu)$. Second, we express $\partial_X \calW_2(\Lambda \uu)$ with the aid of Lemma \ref{le:Leg.Fen} equivalently through the equation
\begin{align}\label{indi.conj}
\langle \mathbold{p},\Lambda \uu\rangle_{X^*\times X}=\calW_2(\Lambda \uu)+\calW_2^*(\mathbold{p}).
\end{align} Third, by \textsc{Ekeland \& Temam} \cite[Proposition 1.2, p. 87]{EkeTem76CAVP}, the conjugate $\calW_2^*$ is given by 
\begin{align*}
\calW_2^*(\BB)=\int_\Omega \imath_{\overline{B}_{\mathbb{R}^{m\times d}}(0,1)}(\BB(\xx))\dd \xx,
\end{align*} with the indicator function $\imath_{\overline{B}_{\mathbb{R}^{m\times d}}(0,1)}\rightarrow \lbrace 0,+\infty\rbrace$ defined by \begin{align*}
\imath_{\overline{B}_{\mathbb{R}^{m\times d}}(0,1)}(\AB)=\begin{cases} 0 \quad &\text{if }\vert  \AB\vert \leq 1 \,\\
+\infty \quad &\text{otherweise.}
\end{cases}
\end{align*} This implies
\begin{align*}
\calW_2^*(\BB)=\begin{cases} 0 \quad &\text{if }\vert  \BB(\xx)\vert \leq 1 \quad \text{a.e. in }\Omega\,\\
+\infty \quad &\text{otherweise.}
\end{cases}
\end{align*} Inserting the latter expression into the equality \eqref{indi.conj}, we obtain
\begin{align*}
\int_{\Omega}\mathbold{p}(\xx):\nabla \uu(\xx)\dd \xx=\int_{\Omega}\vert \mathbold{p}(\xx)\vert \dd \xx
\end{align*} and $\vert \mathbold{p}(\xx)\vert\leq 1$ a.e. in $\Omega$. Since $\mathbold{p}(\xx):\nabla \uu(\xx)\leq \vert \mathbold{p}(\xx)\vert$ by the \textsc{Fenchel--Young} (or \textsc{Cauchy--Schwarz}) inequality, we deduce
\begin{align*}
\mathbold{p}(\xx):\nabla \uu(\xx)=\vert \mathbold{p}(\xx)\vert\quad \text{a.e. in }\Omega.
\end{align*} Therefore, $\mathbold{p}(\xx)\in \overline{B}_{\mathbb{R}^{m\times d}}(0,1)$ if $\vert \nabla \uu(\xx)\vert=0$  and $\mathbold{p}(\xx)=\frac{ \nabla \uu(\xx)}{\vert \nabla \uu(\xx)\vert}$ otherwise. We obtain $\mathbold{p}(\xx)\in \mathbf{Sgn}(\nabla \uu(\xx))$ a.e. in $\Omega$. Now let $\uu_n\overset{*}{\rightharpoonup} \uu$ in $\rmL^{\infty}(0,T;U)\cap \rmH^1(0,T;V)$, $\uu_n \rightarrow \uu$ in $\rmL^{2}(0,T;V)$, and $\xii_n\overset{*}{\rightharpoonup} \xii $ in $\rmL^{2}(0,T;U^*)$ as $n\rightarrow \infty$ such that $\xii_n(t)\in \partial \calE(\uu_n(t))$ for almost every $t\in (0,T)$, $\sup_{n\in\mathbb{N},t\in[0,T]}\calE(\uu_n(t))\leq C_2$, and
 \begin{align}\label{limsup}
 \limsup_{n \rightarrow \infty}\int_0^T \langle \xii_n(t),\uu_n(t)\rangle_{U^*\times U} \dd t \leq \int_0^T \langle \xii(t),\uu(t)\rangle_{U^*\times U} \dd t.
 \end{align} We note that we can decompose $\xii=\zzeta-\uu\in V^*$ with $\zzeta \in \partial \calW(\uu)$. Then, defining $\zzeta_n:=\xii_n+\uu_n$, there holds $\zzeta_n\in \partial W(\uu_n)$ and $\zzeta_n\rightharpoonup \zzeta:=\xii+\uu$ in $\rmL^2(0,T;U^*)$. By the \textsc{Lions--Aubin} lemma, we obtain the strong convergence of $\uu_n\rightarrow \uu$ in $\rmC([0,T];H)$. Thus, in view of \eqref{limsup}, we deduce
 \begin{align*}
 \limsup_{n \rightarrow \infty}\int_0^T \langle \zzeta_n(t),\uu_n(t)\rangle_{U^*\times U} \dd t \leq \int_0^T \langle \zzeta(t),\uu(t)\rangle_{U^*\times U} \dd t.
 \end{align*} Since $W$ is convex, by \cite[Lemma 1.2]{BrCrPa70PNMB} and \cite[Theorem A.3]{Bach22NDIS}, there holds $\zzeta (t)\in \partial_U W(\uu(t))$ in $U^*$ and $W(\uu_n(t))\rightarrow W(\uu(t))$ as $n\rightarrow \infty$ a.e. in $(0,T)$, whence $\xii(t)\in\partial_U\calE(\uu(t))$ a.e. in $(0,T)$ and $\calE(\uu_n(t))\rightarrow \calE(\uu(t))$ as $n\rightarrow \infty$ a.e. in $(0,T)$. We proceed with showing the control of the subgradient of $\calE$, i.e., Condition \ref{eq:cond.E1.8}. Let $\uu\in D(\partial \calE)$ and $\xii\in \partial_U \calE(\uu)$. Then, by \textsc{H\"{o}lder}'s and \textsc{Young}'s inequality, the \textsc{Sobolev} embedding theorem, we obtain
 \begin{align*}
 \langle \xii,\vv \rangle_{U^*\times U} &= \int_\Omega \left(\vert \nabla \uu \vert^{p-2}\nabla \uu : \nabla \vv+\mathbold{p}:\nabla \vv+(\vert\uu\vert^2-1)\uu\cdot \vv\right) \dd \xx\\
 &\leq C \left(\Vert \uu\Vert_{\rmW^{1,p}_0(\Omega)^m}^{p-1}+\Vert \mathbold{p} \Vert_{\rmL^{\infty}(\Omega)^{m\times d}}\right)\Vert \vv\Vert_{\rmW^{1,p}_0(\Omega)^m}\\
 &\quad +\Vert \uu\Vert_{\rmL^{4}(\Omega)^m}^3\Vert \vv\Vert_{\rmL^{4}(\Omega)^m}+ \Vert \uu\Vert_{\rmL^{2}(\Omega)^m}\Vert \vv\Vert_{\rmL^{2}(\Omega)^m}\\
 &\leq C \left(\Vert \uu\Vert_{\rmW^{1,p}_0(\Omega)^m}^{p-1}+\Vert \mathbold{p} \Vert_{\rmL^{\infty}(\Omega)^{m\times d}}+\Vert \uu\Vert_{\rmL^{4}(\Omega)^m}^3+ \Vert \uu\Vert_{\rmL^{2}(\Omega)^m}\right)\Vert \vv\Vert_{U}\\
 &\leq C\left(1+\frac{1}{p}\Vert \uu\Vert_{\rmW^{1,p}_0(\Omega)^m}^{p}+\frac{1}{4}\Vert \uu\Vert_{\rmL^{4}(\Omega)^m}^4 +\frac{1}{2}\Vert \uu\Vert^2_{\rmL^{2}(\Omega)^m} \right) \Vert \vv\Vert_{U}\\
  &\leq C\left(1+\frac{1}{p}\Vert \uu\Vert_{\rmW^{1,p}_0(\Omega)^m}^{p}+\frac{1}{4}\Vert \uu\Vert_{\rmL^{4}(\Omega)^m}^4-\frac{1}{2}\Vert \uu\Vert_{\rmL^{2}(\Omega)^m}^2 +\Vert \uu\Vert_{\rmL^{2}(\Omega)^m}^2\right) \Vert \vv\Vert_{U}\\
 &\leq C\left(1+\calE(\uu)+\Vert \uu\Vert_{\rmL^{2}(\Omega)^m}\right)\Vert \vv\Vert_{U}\\
 &\leq C\left(1+\calE(\uu)+\Vert \vv\Vert_{\rmW^{1,p}_0(\Omega)^m}\right)\Vert \vv\Vert_U 
\end{align*} for all $\vv \in U=\rmW^{1,p}_0(\Omega)^m\cap \rmL^4(\Omega)^m$, where we also used the fact that $\mathbf{Sgn}$ is uniformly bounded, from which \ref{eq:cond.E1.8} follows. Finally, we verify the assumptions on the perturbation $B$. The continuity condition \ref{eq:B1.1} can be easily checked with the dominated convergence theorem. \\
Ad \ref{eq:B1.2}. Let $\uu\in \DOM(\calE) $ and $v,w\in V$. Then, by the \textsc{H\"{o}lder \& Young} inequalities as well as the \textsc{Sobolev} embedding theorem, there holds
\begin{align}\label{eq:App.est.B}
\langle B(\uu,\vv),\ww\rangle_{V^*\times V}&= 
\int_{\Omega} \bb(\xx,t,\uu(\xx),\vv(\xx))\cdot\ww(\xx) \dd \xx\notag \\&\leq C_b\int_{\Omega} ( \vert \uu(\xx)\vert^{q-1}+ \vert \vv(\xx)\vert^{r-1})\vert \ww(\xx)\vert \dd x\notag \\
&\leq C \left(\Vert \uu \Vert_{\rmL^{(q-1)2d/(d+2)}(\Omega)^m}^{q-1}+\Vert \vv \Vert_{\rmL^{(r-1)2d/(d+2)}(\Omega)^m}^{r-1}\right)\Vert \ww \Vert_{\rmL^{(2d/(d-2)}(\Omega)^m}\notag \\
&\leq C \left(\Vert \uu \Vert_{\rmW_0^{1,p}(\Omega)^m}^{q-1}+\Vert \vv \Vert_{\rmL^2(\Omega)^m}^{r-1}\right)\Vert \ww \Vert_{\rmH^1_0(\Omega)^m} \notag \\
&\leq C\left( (1+\calE(\uu))+\Vert \vv \Vert_{\rmL^2(\Omega)^m}^{2}\right)^\frac{1}{2}\Vert \ww \Vert_{\rmH^1_0(\Omega)^m},
\end{align} where $\tilde{c}\in (0,1)$. Recalling that the conjugate is given by $\Psi^*(\xii)=\frac{1}{2}\Vert \xii\Vert^2_{V^*}$ for all $\xii\in V^*=\rmH^{-1}(\Omega)^m$, we conclude \ref{eq:B1.2}. Therefore, for every initial datum $\uu_0\in \DOM(\calE)= \rmW_0^{1,p}(\Omega)^m\cap \rmL^4(\Omega)^m$ and $\vv_0\in \rmL^2(\Omega)^m$, there exists a weak solution
  \begin{align*}
  \uu\in \rmL^{\infty}(0,T;U)\cap\rmH^1(0,T;V)\cap \rmW^{1,\infty}(0,T;H)\cap \rmH^{2}(0,T;U^*)
\end{align*} to (P1) in the sense that 
\begin{align*}
&\int_0^T \Big( \langle \uu'', \vv \rangle_{U^*\times U}+\int_\Omega \big(  \nabla \partial_t\uu:\nabla \vv + \vert\nabla \uu\vert^{p-2}\nabla \uu :\nabla \vv+ (\vert \uu\vert^2-1)\uu\cdot \vv +\mathbold{p}:\nabla \vv \\
&\quad +\bb(\xx,t,\partial_t\uu,\uu) \big) \dd x+  \langle \ff,\vv\rangle_{U^*\times U} \Big) \dd t \quad \text{for all }\vv \in \rmL^2(0,T;U)
\end{align*} with $\mathbold{p}(\xx,t)\in \mathbf{Sgn}(\nabla \uu(\xx,t))$ a.e. in $\Omega_T$, and the energy-dissipation inequality 
\begin{align*}
&\frac{1}{2}\Vert \uu'(t)\Vert^2_{\rmL^2(\Omega)}+\calE(\uu(t)) +\int_s^t \Psi(\uu'(r))\dd r  \\
& +\int_s^t\Psi^*(\mathbold{f}(r)-\bb(r,\uu(r),\uu'(r))-\uu''(r)-\xii(r)) \dd r\\
&\leq \frac{1}{2}\Vert \uu'(s)\Vert^2_{\rmL^2(\Omega)}+\calE(\uu(s)) +\int_s^t \langle \mathbold{f}(r)-\bb(r,\uu(r),\uu'(r)),\uu'(r)\rangle_{V^*\times V}\dd r
\end{align*} holds for all $t\in [0,T]$ if $s=0$ and a.e. $s\in(0,t)$, where $\xii\in \rmL^{\infty}(0,T;U^*)$ and $\xii(t)=-\Delta_p \uu(t)-(\vert \uu(t)\vert^2-1)\uu(t)$ in $U^*=\rmW^{-1,p^*}(\Omega)^m+\rmL^{\frac{4}{3}}(\Omega)^m$ a.e. in $t\in(0,T)$. 
\subsection{Differential inclusion II} \label{se:app.DI1B}
In the following example, we will cover Case \textbf{(b)} while also highlighting the difference between Case \textbf{(a)} and \textbf{(b)}. Hence, we consider the initial-boundary value problem
\begin{align*} 
\text{(P2)}
\begin{cases}
\partial_{tt}\uu-\!\Delta \partial_t \uu+\mathbold{\eta}-\Delta_p \uu+(\vert \uu\vert ^2-\!1)\uu -\! \nabla \cdot \mathbold{p}+\bb\left(\xx,t,\uu\right) = \ff\, \text{ in } \Omega_T,\\
\mathbold{p}(\xx,t)\in \mathbf{Sgn}\left(\nabla \uu(\xx,t)\right)\quad\text{a.e. in } \Omega_T,\\
\mathbold{\eta}(\xx,t)\in \partial \psi(\xx,\partial_t \uu(\xx,t))\quad\text{a.e. in } \Omega_T,\\
\uu(\xx,0)\,\,=\uu_0(\xx) \quad \text{on } \Omega,\\
\uu'(\xx,0)=\vv_0(\xx) \quad \,\text{on } \Omega, \\
\uu(\xx,t)\,\,=0 \quad  \qquad \text{on } \partial \Omega\times[0,T],
\end{cases}
\end{align*} where we assumed that the contribution of $\partial_t \uu$ in the perturbation from the first example is variational, i.e., it has a potential $\psi:\Omega \times\mathbb{R}^m\rightarrow \mathbb{R}$ so that the perturbation $\bb$ has only a contribution from $\uu$. Here, $\psi(\xx,\cdot)$ is for almost every $\xx\in \Omega$ a proper and convex function with subdifferential $\partial \psi:\mathbb{R}^{m}\rightrightarrows \mathbb{R}^{m}$ such that $\psi'(\xx,\cdot)$ is a \textsc{Carath\'{e}odory} function and satisfies the following growth conditions: there exists a number $r>1$ and constants $c_1,c_2,\tilde{C_1}>0$ such that 
\begin{align*}
c_1\left(\vert \zz\vert^r-1\right)\leq &\psi(\xx,\zz)\leq \tilde{C_1}\left( 1+\vert \zz\vert^r\right),\\
&\vert \eta \vert \leq c_1\left(1+\vert \zz\vert^{r-1}\right) \quad \text{for all }\eta \in \partial \psi(\xx,\zz)
\end{align*} for almost every $\xx \in \Omega$ and all $\zz\in \mathbb{R}^m$. As we discussed in Sections \ref{su:Diss.LDS} a prototypical example is $\psi(\zz)=\frac{1}{r}\vert \zz\vert^r+\vert \zz\vert$. 
We choose the same function spaces for $U=\rmW_0^{1,p}(\Omega)^m\cap \rmL^4(\Omega)^m,V=\rmH_0^{1}(\Omega)^m, W=\rmL^{\max\lbrace 2,r\rbrace}(\Omega)^m, \WW=\rmL^{\max\lbrace 2,q\rbrace}(\Omega)^m,$ and $H=\rmL^2(\Omega)^m$ as above. Again $p,q,r\geq 1$ are to be chosen suitably. The external force is assumed to satisfy the weaker assumption  $\mathbold{f}\in \rmL^{r^*}(0,T;W^*)+ \rmL^2(0,T;V^*)$. Further, we assume $u_0\in U$ and $v_0\in H$. In this case, the dissipation potential $\Psi: V\rightarrow \mathbb{R}$ is given 
\begin{align*}
\Psi(\vv)=\frac{1}{2}\int_{\Omega}\vert \nabla \vv(\xx)\vert^2 \dd \xx +\int_{\Omega}\psi(\xx, \vv(\xx)) \dd \xx.
\end{align*} The conjugate $\Psi^*:V^*\rightarrow \mathbb{R}$ is by Lemma \ref{le:Bre.Att} given by the expression
\begin{align*}
\Psi^*(\xii)= \min_{\mathbold{\eta}\in W^*}\left( \frac{1}{2}\Vert \xii-\mathbold{\eta}\Vert^2_{-1,2}+\int_\Omega \psi^*(\xx,\mathbold{\eta}(\xx))\dd \xx\right),
\end{align*} where $\psi^*:\Omega \times \mathbb{R}^{m}\rightarrow \mathbb{R}$ is the conjugate of $\psi$.  The energy functional $\calE: V\rightarrow [0,+\infty]$ is given as in the previous example. The perturbation $B:[0,T]\times W\rightarrow V^*$ is consequently given by
\begin{align*}
\langle B(t,\uu),\ww \rangle_{V^*\times V}= \int_{\Omega} \bb(\xx,t,\uu(\xx))\cdot \ww(\xx) \dd \xx.
\end{align*} It is readily seen that the assumptions for the dissipation potential follow from the assumptions on $\psi$. From $\frac{1}{2}\Vert \vv\Vert_{\rmH_0^1}^2\leq \Psi(\vv)$, we find $\Psi^*(\xii)\leq \frac{1}{2}\Vert \xii\Vert_{-1,2}^2$, and if we choose, e.g., 
\begin{align*}
d\geq 3,\, p\in (1,+\infty),\,q\in 
 \begin{cases} [1,q^*) \quad \text{if }p\in (1,3),\\
[1,1+p)\cap (1,2d/(d-2)]\quad \text{if }p\geq 3,
 \end{cases} 
 \end{align*} where $q^*=\min\lbrace \frac{d(p+2)}{2(d-p)},\frac{3d+4}{d},\frac{2d}{(d-2)},p/2+1\rbrace$, $d\geq 3,\, p\in (1,+\infty)$, we obtain again the estimate \eqref{eq:App.est.B} without any restriction on $r$. However, from the condition $V \overset{c}{\hookrightarrow} W$, we obtain the restriction $r \in [1,1+p)\cap (1,2d/(d-2)]$, which is now a larger range as opposed to the previous case. Simple calculations show that the same restrictions for the exponents $p,q$ and $r$ hold in the dimensions $d=1$ and $d=2$. Again, we obtain for every initial values $u_0\in \DOM(\calE)$ and $v_0\in H$, the existence of a weak solution 
 \begin{align*}
   u\in\rmC_{w}([0,T];U)\cap \rmH^1(0,T;V) \cap\rmW^{1,\infty}(0,T;H)\cap\rmW^{1,r}(0,T;W)\cap \rmW^{2,r^*}(0,T;U^*)
\end{align*} with $r^*=\min\lbrace 2,r^*\rbrace$ to (P2) such that the initial conditions $u(0)=u_0$ and $u'(0)=v_0$ are satisfied, the integral equation
 \begin{align*}
&\int_0^T \Big( \langle \uu'', \vv \rangle_{U^*\times U}+\int_\Omega \big(  \nabla \partial_t\uu:\nabla \vv +\mathbold{\eta} \cdot \vv+ \vert\nabla \uu\vert^{p-2}\nabla \uu :\nabla \vv +\mathbold{p}:\nabla \vv \\
&\quad + (\vert \uu\vert^2-1)\uu\cdot \vv +\bb(\xx,t,\uu) \big) \dd x+  \langle \ff,\vv\rangle_{U^*\times U} \Big) \dd t \quad \text{for all }\vv \in \rmL^{\max\lbrace 2,r\rbrace}(0,T;U)
\end{align*} with $\mathbold{p}(\xx,t)\in \mathbf{Sgn}(\nabla \uu(\xx,t))$ and $\mathbold{\eta}(\xx,t)\in \partial \psi(\xx,\partial_t \uu(\xx,t))$ a.e. in $\Omega_T$ is fulfilled, and the energy-dissipation inequality \eqref{sol:EDI} holds. 

\end{section}

\section{Conclusion}

In this work, we have for the first time demonstrated the existence of strong solutions to doubly nonlinear and multivalued second-order evolution equations under general conditions. The set-valued and nonlinear operators involved are the subdifferentials of nonlinear functionals, namely the dissipation potential $\Psi
\rightarrow \mathbb{R}$, and the energy functional $\calE
\rightarrow (-\infty,+\infty]$. We assumed that the domains $V$ and $U$ generally do not satisfy $U\hookrightarrow V$ or $V \hookrightarrow U$. The existence proof is based on a regularization of the problem through a generalized Moreau-Yosida regularization of the dissipation potential $\Psi$ and a semi-implicit time discretization. Based on the time discretization, which is suitable for computing numerical approximation solutions, a variational approximation scheme was derived, and the existence of discrete solutions was proven. Subsequently, the compactness of the approximate solutions in suitable Bochner spaces was demonstrated, and for the limit $\tau\rightarrow 0$ of vanishing time step size, the solvability of the regularizing problem was proven. Following this, the solvability of the original problem was demonstrated for the case when the regularization parameter tends to zero. Finally, applications of the abstract theory were presented.

\section{Author’s Statements}
All aspects of the present work were solely conducted by the author. The author declares no conflicts of interest.

\bibliographystyle{acm}
\bibliography{alex_pub,bib_aras}

%
\end{document}